\pdfoutput=1
\documentclass[3p,times]{elsarticle}

\usepackage{latexsym}
\usepackage{macro}
\usepackage[babel]{csquotes}
\usepackage{indentfirst}
\usepackage[english]{babel}
\usepackage{cancel}
\usepackage{array}
\usepackage{color}
\usepackage{tabularx}
\usepackage{graphicx}
\usepackage{amsfonts}
\usepackage{amsmath}
\usepackage{amsthm}
\usepackage{amssymb}
\usepackage{dsfont}
\usepackage{upgreek}
\usepackage{bm}
\usepackage{enumitem}
\usepackage{multirow}
\usepackage[dvipsnames]{xcolor}
\usepackage{mathtools}
\usepackage{caption}
\usepackage{subcaption}
\usepackage{booktabs}
\usepackage{relsize}
\usepackage{tikz}
\usepackage{placeins}
\newcommand{\figscale}{0.588}
\newcommand{\figB}[1]{\includegraphics[scale=0.445,trim=4.2bp 2bp 32bp 16bp,clip]{#1}}
\allowdisplaybreaks
\numberwithin{equation}{section}

\newtheorem{theorem}{Theorem}[section]
\newtheorem{proposition}[theorem]{Proposition}
\newtheorem{lemma}[theorem]{Lemma}

\theoremstyle{definition}
\newtheorem{definition}[theorem]{Definition}
\theoremstyle{remark}
\newtheorem{remark}[theorem]{Remark}

\usepackage{hyperref}
\hypersetup{colorlinks=true, linkcolor=blue, citecolor=red, urlcolor=blue}

\makeatletter
\def\ps@pprintTitle{%
	\let\@oddhead\@empty
	\let\@evenhead\@empty
	\def\@oddfoot{\footnotesize\itshape Preprint\hfill\today}%
	\let\@evenfoot\@oddfoot}
\makeatother

\begin{document}

\begin{frontmatter}

\title{A configurational numerical viscosity for thermodynamically compatible finite volume schemes of compressible multiphase flows}

\author[Inria,UniTN]{Davide Ferrari\corref{cor1}}
\ead{davide.ferrari@unitn.it}
\cortext[cor1]{Corresponding author}
\author[UniTN]{Ilya Peshkov}
\ead{ilya.peshkov@unitn.it}
\author[UniTN]{Michael Dumbser}
\ead{michael.dumbser@unitn.it}
\address[Inria]{Inria, Universit\'e de Strasbourg, 67000 Strasbourg, France}
\address[UniTN]{Laboratory of Applied Mathematics, DICAM, University of Trento, via Mesiano 77, 38123 Trento, Italy}

\begin{abstract}
We consider the fluid subsystem of a variationally derived hyperbolic thermodynamically compatible (HTC) model of compressible multiphase flows. Its total energy is in general not convex in the natural state variables, so a scalar numerical viscosity acting on these variables need not produce non-negative entropy. We reformulate the equations in terms of the volume fractions, making explicit the distinction between changes in the amount of each phase and changes in its intrinsic state. This structure guides the construction of a configurational numerical viscosity in phase variables in which the energy is convex, while the scheme retains the original mathematical structure of the HTC model. At the semi-discrete level, the new schemes satisfy the following properties: I)~local conservation of the phase masses, of the mixture momentum and of the total energy, which is not evolved, with a single-valued numerical energy flux; II)~a cell entropy inequality whose production is non-negative face by face and phase by phase, without any convexity assumption in the evolved variables and without clipping; III)~positivity of the partial densities and of the volume fractions under explicit time-step conditions; IV)~exact preservation of contact equilibria, across which the volume fractions are transported by the upwind scheme with the material velocity; V)~equivalence with a scheme for the density-weighted volume fractions of the original model. The key ingredients of the new schemes are: i)~an exact decomposition of the jump of each phase across a face into a configurational part, dissipated with the material speed, and an intrinsic part, dissipated with the acoustic speed; ii)~the form of this dissipation as an exchange of each phase between neighbouring cells, with non-negative weights, from which positivity follows; iii)~a side-wise transfer of the dissipation to the evolved variables through the chain rule of each cell, which turns the energy defects into Bregman divergences of the single-fluid energies; iv)~a direction of the thermodynamic correction of the reversible part built in the phase variables, whose denominator is non-negative. Properties I)--V) are first proved mathematically and are then verified by suitable numerical tests, which include a shock tube in which the Rusanov viscosity produces negative entropy and the interaction of a shock with a water column.
\end{abstract}


\end{frontmatter}

\section{Introduction}\label{sec:intro}
Compressible flows of several immiscible fluids separated by material interfaces, such as bubbles or droplets in a liquid, shocks interacting with liquid columns or cavitating flows, are often computed with diffuse-interface models \cite{baernunziato,SaurelAbgrall,AbgrallSaurel,SaurelCavitation,Saurel2018}. For models of Baer--Nunziato type, the literature ranges from their mathematical analysis and the solution of the Riemann problem \cite{Embid,AndrianovWarnecke,Schwendeman,DeledicquePapalexandris} to path-conservative and pressure-based schemes \cite{USFORCE2,ReAbgrall2022} and to models with more than two phases \cite{Herard1,Herard2,Pelanti3phase,FerrariPeshkovDumbser}. Each phase is described as a continuum with a volume fraction, and an interface is the thin region in which the volume fractions change. A numerical scheme for such models has to regularize shocks without diffusing the interfaces more than necessary, while preserving the conservation properties and the thermodynamic structure of the model. With a scalar numerical viscosity applied uniformly to all variables, these goals can conflict. In Rusanov-type schemes, and more generally in scalar parabolic numerical regularizations based on characteristic speeds, the dissipative scale is set by the fastest waves. In a liquid, the acoustic speed can exceed the material velocity by an order of magnitude or more. The same acoustic scale then controls the numerical diffusion of an interface, even though the interface is transported with the material velocity: in water-air computations, this can lead to substantial interface smearing.

\paragraph{Hyperbolic thermodynamically compatible models of mixtures}
The models considered here belong to the class of symmetric hyperbolic thermodynamically compatible (SHTC) systems introduced by Godunov and Romenski \cite{God1961,GodRom2003}. Their two-phase and multiphase versions \cite{Romenski2004,RomDrikToro2010,Romenski2016,Thein2022} evolve the relative velocities of the phases as state variables and combine conservative fluxes with non-conservative terms reflecting the geometric structure of the model. Interfacial friction and pressure relaxation are described by algebraic source terms with non-negative entropy production. Numerical schemes for these models include high-order ADER discontinuous Galerkin, all-Mach-number and implicit--explicit finite volume methods, and exactly curl-free schemes \cite{Rio-Martin2023,Thomann2022,Thomann2024,Chiocchetti2023,RioMartinDhaouadiDumbser2025}. The unified HTC multiphase model of \cite{FerrariPeshkovDumbser} extends this description to elastic and heat-conducting media, in the spirit of the unified first-order formulation of continuum mechanics \cite{DUMBSER2016824}, for which structure-preserving semi-implicit schemes have recently been developed \cite{DumbserBoscheriTavelliThomann2026}. In this work we consider its fluid subsystem, with relative velocities, interfacial friction and pressure relaxation, but without elastic and heat-conducting effects.

\paragraph{Thermodynamically compatible schemes}
The thermodynamic structure of the continuous model should also be reflected in its numerical approximation: the discrete evolution should satisfy both an energy balance and an entropy inequality. Entropy-conservative and entropy-stable schemes provide a systematic framework for enforcing discrete entropy relations, following the work of Tadmor \cite{Tadmor1987,Tadmor2003}, also for non-conservative systems \cite{CastroFjordholmMishraPares2013}; recent developments include discontinuous Galerkin formulations \cite{Chan2018,Chan2025}.
For fluid systems, these methods typically discretize the conservation laws, including total energy, and enforce an additional discrete entropy relation. Hyperbolic thermodynamically compatible (HTC) schemes \cite{ABGRALL2018640,ABGRALL2023127629,HTC2022,Busto2021b,Busto2022a,BustoDumbserMHD2023,BustoDumbser2024} reverse this choice: they evolve the entropy densities and construct the discretization so that the total-energy balance follows as a consequence. Their convergence towards dissipative weak solutions has been proved for the Euler equations \cite{DumbserLukacovaThomann2026}, semi-implicit variants whose time step is restricted only by the material velocity have recently been proposed \cite{DumbserPuppoRinaldiThomann2026}, and structure-preserving HTC schemes for nonlinear SHTC systems with involutions have been developed in \cite{LuccaDumbser2025}, together with exactly div--curl--grad compatible and energy-conservative discontinuous Galerkin schemes \cite{AbgrallDumbserMaire2025}. They have been applied, among others, to a two-fluid model with two entropy inequalities \cite{Thomann2023} and to Lagrangian gas dynamics and continuum mechanics \cite{boscheri2023,boscheri2025}, and a space--time discontinuous Galerkin method that uses the state and the dual variables together has recently been proposed in \cite{DumbserPerugiaZampa2026}; see also \cite{AbgrallDumbserMaireZampa2026} for structure-preserving continuous--discontinuous Galerkin schemes.
In the finite volume HTC framework mentioned so far the numerical flux has two parts. The reversible part is made exactly compatible with the energy conservation law, either through a path integral in the dual variables or by correcting a central flux along a prescribed direction with one scalar coefficient per face, the so-called Abgrall correction \cite{ABGRALL2018640}, which has also been reinterpreted and extended to residual distribution and discontinuous Galerkin schemes \cite{AbgrallOffnerRanocha2022,GaburroOffnerRicchiutoTorlo2023}. The dissipative part is a numerical viscosity with non-negative entropy production. A common choice is a scalar viscosity acting on the evolved variables, whose production is proportional to the pairing between the jumps of the state and of the thermodynamic forces. This pairing is non-negative when the energy is convex in the evolved variables \cite{ABGRALL2023127629}. A viscosity written directly in the forces, in the spirit of \cite{Tadmor1987}, has a non-negative production by construction.

\paragraph{Numerical dissipation for multiphase flows}

Two observations motivate the present work. The first is that, for the model considered here, the convexity on which the classical scalar viscosity relies is not available. The energy is convex in suitable phase variables, but not in the variables in which the scheme is written as soon as the phases slip with respect to each other (Section~\ref{subsec:three_representations}); the pairing between the jumps of the state and of the thermodynamic forces, and with it the entropy production of the scalar viscosity, can therefore be negative. Section~\ref{subsec:test_shocktube} shows that this happens in a shock tube, in the regions where one of the phases is present only in traces, and that clipping the negative production breaks the conservation of the total energy.

A dissipation written in the forces, \(\mathbf K\Delta\pp\) with a symmetric positive semi-definite mobility \(\mathbf K\), guarantees the sign of the production, but does not by itself determine how the viscosity acts on the evolved variables. For small jumps, the jumps of the forces and of the state are related by the Hessian of the energy, \(\Delta\pp\approx\mathbf H\,\Delta\qq\) with \(\mathbf H=\nabla^2_{\qq}\rhoE\), so that such a viscosity acts on perturbations of the state through \(\mathbf K\mathbf H\). Reproducing a scalar diffusion in all directions of the state would require \(\mathbf K\) proportional to \(\mathbf H^{-1}\), which is indefinite where the phases slip and does not exist where they do not, since the configurational direction then lies in the kernel of \(\mathbf H\). In the phase variables, where the energy is convex, a different obstruction remains: the energy has flat configurational directions, along which the phase variables change while the forces do not. Across a contact between phases with uniform intrinsic states, a common velocity and a common pressure, the forces are therefore continuous while the volume fractions jump, and any dissipation of the form \(\mathbf K\Delta\pp\) vanishes. More generally, a non-negative production says nothing about how the phase masses are exchanged between the cells, nor about the products, such as the partial densities \(\crho_a=\alpha_a\rho_a\), of which the state is made. These are requirements of the multiphase problem rather than shortcomings of the general frameworks, which were designed for systems in which they do not arise; we formulate them in Section~\ref{subsec:viscous_requirements}.

The second observation concerns the interfaces. The volume fractions carry the geometry of the phases and are transported with the material velocity, so that their dissipation could be smaller than the acoustic one. The volume fraction of a phase, however, cannot be dissipated independently of the mass, momentum and entropy that the phase carries: unrelated changes of the amount of a phase and of its extensive variables change its intrinsic state by an amount inversely proportional to the volume fraction, which is large precisely where the phase is present only in traces. The central question of this work is therefore how to reduce the dissipation of the volume fractions with respect to the acoustic scale without decoupling them from the extensive variables, while keeping the conservation of the total energy and a non-negative entropy production.

\paragraph{Three sets of variables}
Our answer rests on a simple observation: the variables that are best suited to \emph{derive} a model are not necessarily those in which its equations are best \emph{solved}, and neither need be those in which its thermodynamic structure is most \emph{transparent}. We therefore use three representations of the same physical state, each with a distinct role. The variational derivation of the multiphase HTC model \cite{FerrariPeshkovDumbser} dictates the density-weighted volume fractions and the relative velocities. The scheme evolves a reformulation in terms of the volume fractions: in these \emph{evolved variables} the admissible set is simple, and the fluxes, the non-conservative products and the energy identities are naturally formulated.
The energy, finally, is convex and positively homogeneous of degree one in the phase variables \(Z_a=\alpha_a(1,Y_a)\), where \(Y_a\) collects the intrinsic density, momentum and entropy of phase \(a\): the energy of each phase is the perspective of the energy of a single fluid, and its flat direction, along which the amount of a phase changes at fixed intrinsic state, is explicit. The guiding principle of the method is to build the scheme in the evolved variables and to construct in the phase variables the two ingredients that need a sign, the numerical dissipation and the denominator of the thermodynamic correction. The joint use of different sets of variables within one scheme is also at the heart of \cite{DumbserPerugiaZampa2026}, where the state and the dual variables are discretized in the same space.

\paragraph{Contributions}
To address these points we propose a \emph{configurational viscosity} for finite volume HTC schemes of multiphase flows. The jump of the phase variables across a face is split exactly into a configurational part, which changes the amount of a phase at an averaged intrinsic state, and an intrinsic part, weighted by the harmonic mean of the volume fractions. The first part is dissipated with the speed at which material crosses the face, the second with the acoustic speed. The resulting exchange has non-negative weights: each phase receives only the same phase from the neighbouring cell, and its mass, momentum, entropy and volume are exchanged in the same ratio. The exchange is transferred to the evolved variables side by side, each cell receiving it through the chain rule at its own state. The total energy is then moved between the cells by a single-valued flux, and the energy defect of each cell is a sum of Bregman divergences of the single-fluid energies, which is converted into the entropy of the phase it comes from. More precisely, the main contributions of this work are the following.
\begin{itemize}
	\item A numerical dissipation for multiphase HTC schemes whose entropy production is non-negative face by face and phase by phase without any convexity of the energy in the evolved variables, which preserves the positivity of the phase masses and of the volume fractions, and which preserves contact equilibria, across which it reduces to the upwind scheme driven by the material velocity.
	\item A direction of the Abgrall correction built in the phase variables and transferred side by side, whose denominator is a sum of Bregman divergences and is therefore non-negative. With the original direction, the jump of the evolved variables, the denominator may vanish across interfaces where the phases slip; in the shock--water column computation of Section~\ref{subsec:test_shockcolumn} the time step then collapses by orders of magnitude.
	\item A single principle, the side-wise chain rule, which also yields the discrete product rule between the volume fractions and the density-weighted volume fractions, and hence the semi-discrete equivalence of the two formulations; the conservative flux of the density-weighted volume fractions and the averaged velocity at the faces follow from it instead of being imposed.
\end{itemize}

The rest of the paper is organized as follows. Section~\ref{sec:model} recalls the fluid subsystem of the unified HTC model of mixtures and its formulation in terms of the volume fractions. Section~\ref{sec:geometry} introduces the phase variables and the perspective energy, and explains why the energy is convex in these variables but not in the evolved ones. Section~\ref{sec:framework} recalls the semi-discrete HTC framework, introduces the thermodynamic correction with the direction built in the phase variables and states the requirements on the numerical dissipation. Section~\ref{sec:viscosity} constructs the configurational viscosity: the decomposition of the phase jump, the two dissipative scales, the exchange form, the side-wise transfer and the energy defects. Section~\ref{sec:properties} collects the properties of the resulting scheme: conservation of the total energy, non-negative entropy production, positivity of phase masses and volume fractions, preservation of contact equilibria, and the discrete product rule with the semi-discrete equivalence of the two formulations. Section~\ref{sec:results} presents the numerical results, and Section~\ref{sec:conclusions} concludes the paper. The appendices collect the convexity properties of the energy in the different sets of variables (\ref{app:convexity}), the characteristic structure and the asymptotic limits of the model (\ref{app:characteristics}), and the evaluation of the Bregman divergences for the equations of state used here (\ref{app:bregman}).

\section{The fluid multiphase HTC model and its volume fraction formulation}\label{sec:model}

In this section we recall the fluid subsystem of the unified HTC model for mixtures introduced in \cite{FerrariPeshkovDumbser}, and we rewrite it in terms of the volume fractions, which are the variables used by the numerical method. We keep only what is needed later: the state variables, the energy, the structure of the balance laws and the two relaxation mechanisms that act on the fluid subsystem, namely interfacial friction and pressure relaxation. Elastic and heat-conducting effects, which in the complete model are carried by the distortion field and by the thermal impulses, are not considered.

\subsection{Variables, kinematics and energy}\label{subsec:variables}

\paragraph{Masses and volumes}
We consider a mixture of \(\Nph\) constituents, or phases, labelled by \(a,b=1,\ldots,\Nph\). An infinitesimal mixture element of volume \(V\) and mass \(M\) contains the volume \(\nu_a\) and the mass \(m_a\) of each phase, with \(V=\sum_a\nu_a\) and \(M=\sum_a m_a\). Three densities and a volume fraction follow from these quantities, and it is useful to keep them clearly apart from the beginning:
\begin{equation}\label{eq:densities}
	\rho=\frac{M}{V},
	\qquad
	\rho_a=\frac{m_a}{\nu_a},
	\qquad
	\crho_a=\frac{m_a}{V}=\alpha_a\rho_a,
	\qquad
	\alpha_a=\frac{\nu_a}{V}.
\end{equation}
Here \(\rho\) is the mixture density, \(\rho_a\) the intrinsic density of phase \(a\) (its mass per unit volume of the phase itself), \(\crho_a\) its partial density (its mass per unit volume of the mixture) and \(\alpha_a\) its volume fraction. The volume and mass fractions satisfy
\begin{equation}\label{eq:saturation}
	\sum_{a=1}^{\Nph}\alpha_a=1,
	\qquad
	\sum_{a=1}^{\Nph}c_a=1,
	\qquad
	c_a:=\frac{\crho_a}{\rho},
	\qquad
	\rho=\sum_{a=1}^{\Nph}\crho_a .
\end{equation}
The original HTC formulation does not evolve \(\alpha_a\) but the density-weighted volume fractions
\begin{equation}\label{eq:phi_definition}
	\varphi_a:=\rho\,\alpha_a .
\end{equation}
Although \(\varphi_a\) and \(\crho_a\) are both densities, they describe different things. The partial density \(\crho_a=\alpha_a\rho_a\) measures \emph{how much mass} of phase \(a\) is present; the quantity \(\varphi_a=\alpha_a\rho\) measures \emph{how much volume} phase \(a\) occupies, weighted by the density of the mixture. The two coincide only where \(\rho_a=\rho\). Across an air--water interface, for instance, \(\crho_a\) of the water changes by a factor of the order of the density ratio, whereas \(\varphi_a\) follows the volume fraction. This distinction will reappear several times: \(\crho_a\) is the natural carrier of the mass of a phase, \(\alpha_a\) of its geometry.

\paragraph{Velocities and entropies}
Each phase moves with its own velocity \(\vvv_a\) and has the entropy density \(\eta_a=\crho_a s_a\), where \(s_a\) is its specific entropy; in the model, the phase entropies are transported with the mixture velocity, not with the phase velocities (Section~\ref{subsec:balances}). The mixture momentum, the mixture velocity and the relative velocities are
\begin{equation}\label{eq:kinematics}
	\UU=\sum_{a=1}^{\Nph}\crho_a\vvv_a,
	\qquad
	\VV=\frac{\UU}{\rho},
	\qquad
	\ww_a=\vvv_a-\vvv_\Nph,
	\quad a=1,\ldots,\Nph-1,
\end{equation}
where phase \(\Nph\) is an arbitrarily chosen reference phase. The mixture velocity \(\VV\) is the velocity of the centre of mass; the phase velocities are recovered from \(\VV\) and the relative velocities, \(\vvv_\Nph=\VV-\sum_{b<\Nph}c_b\ww_b\) and \(\vvv_a=\vvv_\Nph+\ww_a\).

\paragraph{Energy}
Each phase is described by a standard single-phase equation of state, given through its specific internal energy \(e_a(\rho_a,s_a)\). Pressure, temperature and specific Gibbs energy of the phase are
\begin{equation}\label{eq:phase_thermo}
	p_a=\rho_a^2\frac{\partial e_a}{\partial\rho_a},
	\qquad
	T_a=\frac{\partial e_a}{\partial s_a},
	\qquad
	g_a=e_a+\frac{p_a}{\rho_a}-T_as_a .
\end{equation}
The total energy density of the mixture is the sum of the kinetic and internal energies of the phases,
\begin{equation}\label{eq:energy_phases}
	\rhoE=\sum_{a=1}^{\Nph}\crho_a\left(\tfrac12|\vvv_a|^2+e_a(\rho_a,s_a)\right).
\end{equation}
In the HTC variables the kinetic part is split into the kinetic energy of the centre of mass and the kinetic energy of the relative motion,
\begin{equation}\label{eq:energy_split}
	\rhoE
	=
	\frac{|\UU|^2}{2\rho}
	+
	\underbrace{\frac12\sum_{a=1}^{\Nph-1}\crho_a|\ww_a|^2-\frac{1}{2\rho}\Bigl|\sum_{a=1}^{\Nph-1}\crho_a\ww_a\Bigr|^2}_{=\,\frac12\sum_{a=1}^{\Nph}\crho_a|\vvv_a-\VV|^2}
	+
	\sum_{a=1}^{\Nph}\crho_a\,e_a\!\left(\frac{\crho_a}{\alpha_a},\frac{\eta_a}{\crho_a}\right).
\end{equation}
The identity under the brace is the parallel-axis theorem for the phase momenta. Two remarks will guide the rest of the paper. First, \eqref{eq:energy_phases} is a sum over phases of one and the same kind of term: mass times the energy per unit mass of a single fluid. Second, the HTC variables mix quantities of different phases (\(\UU\), \(\ww_a\), \(\rho\)), so that this phase-by-phase structure is not visible in them. Section~\ref{sec:geometry} returns to this point.

\subsection{Thermodynamic balances and relaxation}\label{subsec:balances}

\paragraph{State vector and thermodynamic forces}
The original fluid HTC system \cite{FerrariPeshkovDumbser} evolves the state vector
\begin{equation}\label{eq:state_phi}
	\qq^\varphi
	=
	\left(
		\UU,\,
		\crho_1,\ldots,\crho_{\Nph-1},\,
		\ww_1,\ldots,\ww_{\Nph-1},\,
		\eta_1,\ldots,\eta_\Nph,\,
		\rho,\,
		\varphi_1,\ldots,\varphi_{\Nph-1}
	\right)^{\mathsmaller T},
\end{equation}
the variables of phase \(\Nph\) being recovered from \(\crho_\Nph=\rho-\sum_{b<\Nph}\crho_b\) and \(\varphi_\Nph=\rho-\sum_{b<\Nph}\varphi_b\). All fluxes are built from the derivatives of the energy with respect to the state variables, the \emph{thermodynamic forces} or main-field variables \cite{MullerRuggeri1998}. For a state variable \(\chi\) we write \(\rhoE_\chi=\partial\rhoE/\partial\chi\). The forces needed in this section are
\begin{subequations}\label{eq:forces_phi}
\begin{align}
	\rhoE_{\UU}&=\VV,
	\qquad
	\rhoE_{\ww_a}=\crho_a(\vvv_a-\VV),
	\qquad
	\rhoE_{\eta_a}=T_a,
	\qquad
	\rhoE_{\varphi_a}=\frac{p_\Nph-p_a}{\rho},
	\label{eq:forces_phi_a}
	\\
	\rhoE_{\crho_a}&=\left(g_a+\tfrac12|\vvv_a-\VV|^2\right)-\left(g_\Nph+\tfrac12|\vvv_\Nph-\VV|^2\right).
	\label{eq:forces_phi_b}
\end{align}
\end{subequations}
Each of them has a direct physical reading. The force conjugate to \(\ww_a\) is the relative momentum of phase \(a\) with respect to the centre of mass; the force conjugate to the partial density is the difference between the Gibbs energies of the phases, each augmented by its relative kinetic energy; the force conjugate to \(\varphi_a\) is the pressure difference. The mixture pressure, defined by the partial Legendre transform of the energy with respect to the extensive variables, is the volume average of the phase pressures (Dalton's law),
\begin{equation}\label{eq:mixture_pressure}
	P
	:=
	U_i\rhoE_{U_i}
	+\sum_{a<\Nph}\crho_a\rhoE_{\crho_a}
	+\sum_{a=1}^{\Nph}\eta_a\rhoE_{\eta_a}
	+\rho\rhoE_\rho
	+\sum_{a<\Nph}\varphi_a\rhoE_{\varphi_a}
	-\rhoE
	=
	\sum_{a=1}^{\Nph}\alpha_ap_a .
\end{equation}

\paragraph{Balance laws}
With summation over repeated spatial indices, and over \(a<\Nph\) where a phase index is repeated in the momentum and energy fluxes, the fluid HTC system reads
\begin{subequations}\label{eq:system_phi}
\begin{align}
	&\partial_t U_i+\partial_k\left(U_iV_k+P\delta_{ik}+w_{a,i}\,\rhoE_{w_{a,k}}\right)=0,
	\label{eq:system_phi_U}\\
	&\partial_t\crho_a+\partial_k\left(\crho_aV_k+\rhoE_{w_{a,k}}\right)=0,
	&&a=1,\ldots,\Nph-1,
	\label{eq:system_phi_rho_a}\\
	&\partial_t w_{a,k}+\partial_k\left(w_{a,l}V_l+\rhoE_{\crho_a}\right)+V_i\left(\partial_iw_{a,k}-\partial_kw_{a,i}\right)=\Lambda_{a,k},
	&&a=1,\ldots,\Nph-1,
	\label{eq:system_phi_w}\\
	&\partial_t\eta_a+\partial_k\left(\eta_aV_k\right)=\Pi_a,
	&&a=1,\ldots,\Nph,
	\label{eq:system_phi_eta}\\
	&\partial_t\rho+\partial_k\left(\rho V_k\right)=0,
	\label{eq:system_phi_rho}\\
	&\partial_t\varphi_a+\partial_k\left(\varphi_aV_k\right)=\Phi_a,
	&&a=1,\ldots,\Nph-1.
	\label{eq:system_phi_phi}
\end{align}
\end{subequations}
Using \eqref{eq:forces_phi}, the fluxes take a familiar form. The mass flux of phase \(a\) in \eqref{eq:system_phi_rho_a} is \(\crho_aV_k+\rhoE_{w_{a,k}}=\crho_av_{a,k}\): each phase is transported with its own velocity. The momentum flux in \eqref{eq:system_phi_U} is \(\sum_{a=1}^{\Nph}\crho_av_{a,i}v_{a,k}+P\delta_{ik}\): the sum of the convective fluxes of the phases plus the mixture pressure. The relative velocity in \eqref{eq:system_phi_w} is driven by the gradient of \(\rhoE_{\crho_a}\), that is, by the difference of the phase Gibbs energies. Volume fractions and phase entropies are transported with the mixture velocity.

\paragraph{Relaxation and the two laws of thermodynamics}
The source terms describe interfacial friction and pressure relaxation. They are gradient dynamics in state space \cite{PKGBook2018}:
\begin{equation}\label{eq:relaxation}
	\Lambda_{a,k}=-\sum_{b<\Nph}\lambda_{ab,k}\,\rhoE_{w_{b,k}},
	\qquad
	\Phi_a=-\sum_{b<\Nph}\phi_{ab}\,\rhoE_{\varphi_b},
\end{equation}
with symmetric positive semi-definite coefficient matrices, proportional to the inverse relaxation times \(1/\tau^{\mathrm w}\) and \(1/\tau^{\mathrm p}\). By \eqref{eq:forces_phi_a}, friction drives the relative momenta to zero and pressure relaxation drives the phase pressures to a common value. For smooth solutions, multiplying each equation of \eqref{eq:system_phi} by its thermodynamic force and summing gives the conservation of the total energy,
\begin{equation}\label{eq:energy_conservation}
	\partial_t\rhoE
	+
	\partial_k\left[
		V_k\rhoE
		+V_i\left(P\delta_{ik}+w_{a,i}\rhoE_{w_{a,k}}\right)
		+\rhoE_{\crho_a}\rhoE_{w_{a,k}}
	\right]
	=0,
\end{equation}
provided the phase entropy productions compensate the energy removed by the relaxation sources. This fixes
\begin{equation}\label{eq:entropy_production}
	\Pi_a=\frac{c_a}{T_a}\,\mathcal D,
	\qquad
	\mathcal D=\sum_{a,b<\Nph}\left[\sum_k\rhoE_{w_{a,k}}\lambda_{ab,k}\rhoE_{w_{b,k}}+\rhoE_{\varphi_a}\phi_{ab}\rhoE_{\varphi_b}\right]\ \ge0,
\end{equation}
so that \(\sum_aT_a\Pi_a=\mathcal D\). Since \(\Pi_a\ge0\), the mixture entropy \(\eta=\sum_a\eta_a\) satisfies \(\partial_t\eta+\partial_k(\eta V_k)=\sum_a\Pi_a\ge0\). This is a statement about the sources: a phase entropy density \(\eta_a\) can still decrease locally through transport, and nonnegative production does not mean that every local value of \(\eta_a\) must increase. The reversible part of the system conserves the energy by construction, and the irreversible part converts the dissipated mechanical energy into phase entropy: the model is compatible with the first and the second law of thermodynamics.

\paragraph{Asymptotic regimes}
Three limiting regimes clarify the hyperbolic structure of the model, its characteristic speeds and the effect of the relaxation processes on the acoustic scales, which range from the frozen speeds of sound to Wood's speed of the relaxed mixture. The details, including a dispersion analysis that connects the two limits, are collected in \ref{app:characteristics}.
\begin{itemize}
	\item \emph{Frozen regime} (\(\tau^{\mathrm w},\tau^{\mathrm p}\to\infty\)). For two phases with equal velocities (\(\ww_1=\mathbf 0\)), the characteristic speeds of the frozen system are computed explicitly in \ref{app:characteristics}: contact waves moving with \(\VV\) and four thermo-acoustic waves whose speeds involve both phases. For \(\ww_1\neq\mathbf 0\) no closed form is available, and we do not attempt a general characteristic analysis. The numerical dissipation does not rely on one: the bound \(\lambda_{ij}\) of Section~\ref{subsec:two_scales} is computed from the eigenvalues of the quasi-linear matrix of the fluid system, evaluated numerically in each cell, and no complex eigenvalue was detected in the computations of Section~\ref{sec:results}.
	\item \emph{Stiff relaxation} (\(\tau^{\mathrm w},\tau^{\mathrm p}\to0\)). The sources enforce \(\ww_a=0\) and \(p_a=p_\Nph\), and long acoustic waves propagate with Wood's speed of a mechanically equilibrated, phase-wise adiabatic mixture \cite{wood,Bdzil}.
	\item \emph{Pure-phase limit} (\(\alpha_a\to1\) for one phase). The system reduces to the compressible Euler equations of that phase. Diffuse-interface computations keep the other phases at a small residual volume fraction; such \emph{trace} amounts are a regularization of this limit and play an important role in the construction of the numerical dissipation.
\end{itemize}
The test problems of Section~\ref{sec:results} cover the frozen regime (vortex, transported contact and shock tube) and stiff relaxation (shock--water column), and the last three contain phases in traces.

\subsection{Why we solve for the volume fractions}\label{subsec:why_alpha}

\paragraph{The volume fraction formulation}
For \(\rho>0\), combining the mixture mass balance \eqref{eq:system_phi_rho} with \eqref{eq:system_phi_phi} gives
\begin{equation}\label{eq:alpha_equation}
	\partial_t\alpha_a+V_k\,\partial_k\alpha_a=\Psi_a,
	\qquad
	\Psi_a:=\frac{\Phi_a}{\rho},
	\qquad a=1,\ldots,\Nph-1,
\end{equation}
which replaces \eqref{eq:system_phi_phi}. The resulting system evolves
\begin{equation}\label{eq:state_alpha}
	\qq
	=
	\left(
		\UU,\,
		\crho_1,\ldots,\crho_{\Nph-1},\,
		\ww_1,\ldots,\ww_{\Nph-1},\,
		\eta_1,\ldots,\eta_\Nph,\,
		\rho,\,
		\alpha_1,\ldots,\alpha_{\Nph-1}
	\right)^{\mathsmaller T},
\end{equation}
with \(\alpha_\Nph=1-\sum_{b<\Nph}\alpha_b\). From now on \(\qq\) denotes this vector; it is the state vector of the numerical method, and \(\pp=\partial_{\qq}\rhoE\) denotes the corresponding thermodynamic forces. The energy is the same function \eqref{eq:energy_split}, now read as a function of \(\qq\). By the chain rule, the forces conjugate to \(\UU\), \(\crho_a\), \(\ww_a\), \(\eta_a\) are unchanged, while
\begin{equation}\label{eq:forces_alpha}
	\rhoE_{\alpha_a}=\rho\,\rhoE_{\varphi_a}=p_\Nph-p_a,
	\qquad
	\rhoE_{\rho}\big|_{\qq}=\rhoE_\rho\big|_{\qq^\varphi}+\sum_{a<\Nph}\alpha_a\rhoE_{\varphi_a}
	=g_\Nph+\tfrac12|\vvv_\Nph-\VV|^2-\tfrac12|\VV|^2 .
\end{equation}
The mixture pressure is still \(P=\sum_a\alpha_ap_a\); in the Legendre transform \eqref{eq:mixture_pressure} written with \(\qq\), the term \(\sum_a\varphi_a\rhoE_{\varphi_a}\) is absorbed by \(\rho\,\rhoE_\rho|_{\qq}\), and no term \(\alpha_a\rhoE_{\alpha_a}\) appears, consistently with the fact that \(\alpha_a\), a ratio of volumes, does not change under a uniform compression of the mixture. Pressure relaxation becomes \(\Psi_a=-\rho^{-2}\sum_b\phi_{ab}\rhoE_{\alpha_b}\) and acts directly on the pressure differences.

\paragraph{Three reasons}
The two formulations describe the same model (Proposition~\ref{prop:continuous_equivalence} below). We nevertheless prefer \(\alpha_a\), for three reasons.
\begin{enumerate}
	\item \emph{The admissible set is simple.} Physically admissible volume fractions form the simplex \(\alpha_a\ge0\), \(\sum_a\alpha_a=1\), independently of all other variables. In terms of \(\varphi_a\) the same constraint reads \(0\le\varphi_a\le\rho\) and couples the composition to the mixture density. Positivity statements for the numerical scheme are therefore naturally formulated in \(\alpha_a\).
	\item \emph{The phase energies depend only on phase variables.} The internal energy of phase \(a\),
	\[
		\rhoE_{\mathrm{int},a}=\crho_a\,e_a\!\left(\frac{\crho_a}{\alpha_a},\frac{\eta_a}{\crho_a}\right),
	\]
	involves only the variables \((\alpha_a,\crho_a,\eta_a)\) of that phase. In terms of \(\varphi_a\), the intrinsic density reads \(\rho_a=\rho\crho_a/\varphi_a\) and brings in the density of the whole mixture. In \(\alpha_a\) the perspective structure of the internal energy becomes visible (Section~\ref{sec:geometry}). This alone does not make the total energy convex in \(\qq\), where the relative velocities remain; the complete convex structure also requires the phase momenta (Section~\ref{subsec:three_representations} and \ref{app:convexity}).
	\item \emph{It carries the geometry of the phases.} Without relaxation, \eqref{eq:alpha_equation} says that \(\alpha_a\) is a material label transported by the mixture flow. In particular, if every phase has a uniform intrinsic state (constant \(\rho_a\) and \(s_a\)) and all phases share the same velocity and the same pressure, the system reduces to the advection of the volume fractions with the common velocity: an interface between such phases is simply transported. The discrete counterpart of this equilibrium is one of the requirements on the numerical dissipation (Sections~\ref{sec:framework} and~\ref{sec:viscosity}).
\end{enumerate}

\paragraph{Continuous equivalence}
The equivalence is conveniently stated in terms of residuals. Let
\begin{equation}\label{eq:residuals}
	\mathcal R_\rho:=\partial_t\rho+\partial_k(\rho V_k),
	\qquad
	\mathcal R_{\alpha_a}:=\partial_t\alpha_a+V_k\partial_k\alpha_a-\Psi_a,
	\qquad
	\mathcal R_{\varphi_a}:=\partial_t\varphi_a+\partial_k(\varphi_aV_k)-\Phi_a .
\end{equation}

\begin{proposition}[Continuous equivalence]\label{prop:continuous_equivalence}
For smooth fields with \(\rho>0\), \(\varphi_a=\rho\alpha_a\) and \(\Phi_a=\rho\Psi_a\),
\begin{equation}\label{eq:residual_identity}
	\mathcal R_{\varphi_a}=\alpha_a\,\mathcal R_\rho+\rho\,\mathcal R_{\alpha_a},
	\qquad a=1,\ldots,\Nph-1 .
\end{equation}
Consequently, a smooth solution of \eqref{eq:system_phi_rho} and \eqref{eq:alpha_equation} solves \eqref{eq:system_phi_phi}, and conversely. Moreover, the energy residuals of the two formulations coincide,
\begin{equation}\label{eq:energy_residual_identity}
	\rhoE_\rho\big|_{\qq^\varphi}\mathcal R_\rho+\sum_{a<\Nph}\rhoE_{\varphi_a}\mathcal R_{\varphi_a}
	=
	\rhoE_\rho\big|_{\qq}\mathcal R_\rho+\sum_{a<\Nph}\rhoE_{\alpha_a}\mathcal R_{\alpha_a},
\end{equation}
so that the energy conservation law \eqref{eq:energy_conservation} and the entropy productions \eqref{eq:entropy_production} hold unchanged in the volume fraction formulation.
\end{proposition}

\begin{proof}
Expanding the derivatives of the product \(\varphi_a=\rho\alpha_a\),
\[
	\partial_t(\rho\alpha_a)+\partial_k(\rho\alpha_aV_k)-\rho\Psi_a
	=
	\alpha_a\left[\partial_t\rho+\partial_k(\rho V_k)\right]
	+
	\rho\left[\partial_t\alpha_a+V_k\partial_k\alpha_a-\Psi_a\right],
\]
which is \eqref{eq:residual_identity}; the equivalence of the solutions follows because \(\rho>0\). Inserting \eqref{eq:residual_identity} into the left-hand side of \eqref{eq:energy_residual_identity} and using the chain rule \eqref{eq:forces_alpha} gives the right-hand side. All the other terms of the energy identity involve forces that are the same in both representations. Finally, \(\sum_a\rhoE_{\alpha_a}\Psi_a=\sum_a\rhoE_{\varphi_a}\Phi_a\), so the relaxation sources dissipate the same energy in both formulations.
\end{proof}

The identity \eqref{eq:residual_identity} is a product rule. It holds automatically for smooth fields, but a numerical scheme written for \(\qq\) does not satisfy its discrete counterpart unless it is built to do so. Section~\ref{sec:equivalence} shows that the construction of the numerical dissipation developed in Section~\ref{sec:viscosity} provides it as a by-product.

\section{Energy geometry and configurational structure}\label{sec:geometry}

The previous section used two sets of variables for the same model, and the numerical method will use a third one. This is not a complication for its own sake. The variables that are best suited to \emph{derive} a model are not necessarily those in which its equations are best \emph{solved}, and neither need be those in which its thermodynamic structure is most \emph{transparent}:
\begin{itemize}
	\item the HTC variables \(\qq^\varphi\), with the relative velocities \(\ww_a\), are dictated by the variational derivation \cite{FerrariPeshkovDumbser};
	\item the variables \(\qq\) of Section~\ref{subsec:why_alpha} are those in which the scheme is built: they are conservative wherever the physics is, their admissible set is simple, and the fluxes, the non-conservative products and the energy identities of Section~\ref{sec:framework} are written in them;
	\item the phase variables \(Z\) introduced below are those in which the energy is convex and homogeneous of degree one.
\end{itemize}
The variables \(\qq\) are the backbone of the method rather than an exclusive choice, and they act as a bridge between the other two representations. The numerical dissipation and the direction of the thermodynamic correction are constructed in \(Z\) and transferred to \(\qq\); by the discrete product rule of Section~\ref{sec:equivalence}, the same scheme also advances the original variables \(\qq^\varphi\), with identical semi-discrete solutions. This section describes the structure of the energy in \(Z\) and explains why it is needed.

\subsection{Intrinsic and extensive phase variables}\label{subsec:phase_variables}

For each phase we collect the densities per unit volume \emph{of the phase} and per unit volume \emph{of the mixture}:
\begin{equation}\label{eq:YXZ}
	Y_a=(\rho_a,\bm m_a,\sigma_a),
	\qquad
	X_a=\alpha_aY_a=(\crho_a,\bm M_a,\eta_a),
	\qquad
	Z_a=(\alpha_a,X_a),
\end{equation}
with \(\bm m_a=\rho_a\vvv_a\), \(\sigma_a=\rho_as_a\), and \(\bm M_a=\alpha_a\bm m_a=\crho_a\vvv_a\) the momentum of phase \(a\) per unit mixture volume. We call \(Y_a\) the \emph{intrinsic state} of the phase: it describes the fluid as it would be seen from inside the phase, and it is the argument of a standard single-phase equation of state. The \emph{extensive} variables \(X_a\) also account for how much of the phase is present. The vector \(Z_a\) adds the volume fraction, and \(Z=(Z_1,\ldots,Z_\Nph)\), subject to \(\sum_a\alpha_a=1\).

The state \(\qq\) is a function of \(Z\):
\begin{equation}\label{eq:q_of_Z}
	\UU=\sum_a\bm M_a,
	\qquad
	\crho_a,
	\qquad
	\ww_a=\frac{\bm M_a}{\crho_a}-\frac{\bm M_\Nph}{\crho_\Nph},
	\qquad
	\eta_a,
	\qquad
	\rho=\sum_a\crho_a,
	\qquad
	\alpha_a .
\end{equation}
The two vectors have the same number of independent components: \(Z\) contains all the volume fractions, tied by \(\sum_a\alpha_a=1\), whereas \(\qq\) contains \(\Nph-1\) of them. The constraint \(\sum_a\alpha_a=1\) gives the same number of independent variables in the two representations. We call a variation \(\delta Z\) admissible if \(\sum_a\delta\alpha_a=0\).
The Jacobian \(\mathbf J=\partial\qq/\partial Z\) is a rectangular matrix, and for positive partial densities it maps the admissible variations one-to-one onto the variations of \(\qq\). The map \(Z\mapsto\qq\) is linear in every component except the relative velocities, which are differences of ratios. This single non-linearity is responsible for all the differences between the two representations discussed in Section~\ref{subsec:three_representations}.

\subsection{The perspective energy}\label{subsec:perspective}
The energy density of a single phase \(a\), per unit volume of that phase, is the energy density of a single compressible fluid,
\begin{equation}\label{eq:intrinsic_energy}
	E_a(Y_a)=\frac{|\bm m_a|^2}{2\rho_a}+\varepsilon_a(\rho_a,\sigma_a),
	\qquad
	\varepsilon_a(\rho_a,\sigma_a)=\rho_a\,e_a\!\left(\rho_a,\frac{\sigma_a}{\rho_a}\right).
\end{equation}
Multiplication by the volume fraction \(\alpha_aE_a(Y_a)\) gives its contribution to the mixture energy.
Using \(Y_a=X_a/\alpha_a\), we obtain
\begin{equation}\label{eq:perspective}
	\rhoE_a(Z_a)=\alpha_a\,E_a\!\left(\frac{X_a}{\alpha_a}\right),
	\qquad
	\rhoE=\sum_{a=1}^{\Nph}\rhoE_a(Z_a),
\end{equation}
which is the phase decomposition \eqref{eq:energy_phases}.
The function \(\alpha E(X/\alpha)\) is known as the \emph{perspective} of \(E\) in convex analysis \cite{BoydVandenberghe2004}. This transformation also relates thermodynamic functions written in specific variables to those written in densities; see \cite{LukacovaTheinWarneckeYuan2026} and \cite{Warnecke2025}.
In \(\rhoE_a(Z_a)\), the index \(a\) denotes a phase contribution, whereas a subscript identifying a state variable continues to denote a derivative.

\paragraph{Convexity and homogeneity}
We assume thermodynamically stable phase equations of state, with \(e_a\) convex in the specific volume and specific entropy.
This is equivalent to the convexity of \(\varepsilon_a\) in \((\rho_a,\sigma_a)\), and, together with the convexity of the kinetic term \(|\bm m|^2/(2\rho)\), to the convexity of the single-fluid energy \(E_a\) in \(Y_a\) \cite{God1961,MullerRuggeri1998,LukacovaTheinWarneckeYuan2026}. The perspective transfers this property to the extensive variables.

\begin{proposition}[Perspective energy]\label{prop:perspective}
Let \(E_a\) be convex on its domain. Then \(\rhoE_a\) is convex in \(Z_a\) for \(\alpha_a>0\), and it is positively homogeneous of degree one:
\[
	\rhoE_a(\lambda Z_a)=\lambda\,\rhoE_a(Z_a)\qquad\text{for all }\lambda>0 .
\]
Consequently \(\rhoE=\sum_a\rhoE_a\) is convex in \(Z\).
\end{proposition}

\begin{proof}
Homogeneity is immediate from \eqref{eq:perspective}. To verify convexity, consider two states \(Z'=(\alpha',X')\), \(Z''=(\alpha'',X'')\), a weight \(t\in[0,1]\), and set \(\alpha=t\alpha'+(1-t)\alpha''\), \(X=tX'+(1-t)X''\). The intrinsic state of the combination is itself a convex combination of the two intrinsic states,
\[
	\frac{X}{\alpha}
	=
	\frac{t\alpha'}{\alpha}\,\frac{X'}{\alpha'}+\frac{(1-t)\alpha''}{\alpha}\,\frac{X''}{\alpha''},
	\qquad
	\frac{t\alpha'}{\alpha}+\frac{(1-t)\alpha''}{\alpha}=1,
\]
with weights proportional to the volumes. Convexity of \(E_a\) then gives \(\alpha E_a(X/\alpha)\le t\alpha'E_a(X'/\alpha')+(1-t)\alpha''E_a(X''/\alpha'')\). The sum over the phases is convex, also when restricted to the constraint \(\sum_a\alpha_a=1\).
\end{proof}

Homogeneity of degree one has three consequences.
\begin{enumerate}
	\item \emph{The forces depend only on the intrinsic state.} Differentiating the phase energy \eqref{eq:perspective} gives its thermodynamic conjugates,
	\begin{equation}\label{eq:dual_Z}
		\pp^Z_a:=\frac{\partial\rhoE_a}{\partial Z_a}=\left(-p_a,\ \bm\theta_a\right),
		\qquad
		\bm\theta_a:=\nabla_YE_a(Y_a)=\left(g_a-\tfrac12|\vvv_a|^2,\ \vvv_a,\ T_a\right).
	\end{equation}
	The force conjugate to \(\alpha_a\) is therefore minus the phase pressure. The remaining entries are the usual single-fluid conjugates. They depend on \(Y_a\), but not separately on \(\alpha_a\). None of them depends on how much of the phase is present.
	\item \emph{Euler's relation.} \(\rhoE_a(Z_a)=\pp^Z_a\cdot Z_a\): the energy of a phase is the pairing of its forces with its own state. This identity will be used to express the viscous energy defects as Bregman divergences in Section~\ref{subsec:energy_production}.
	\item \emph{Convex, but not strictly.} Along the ray \(\lambda\mapsto\lambda Z_a\) the energy is linear, so the Hessian satisfies \(\nabla^2\rhoE_a(Z_a)\,Z_a=\mathbf 0\). Changing the amount of a phase at fixed intrinsic state is a direction of zero curvature. Strict convexity would require an additional, strictly convex interfacial contribution, which is not part of the present model.
\end{enumerate}
In particular, convexity already gives the inequality on which every dissipative flux with a sign is built: the forces are monotone,
\begin{equation}\label{eq:monotone_Z}
	\bigl(\pp^Z_a(Z_a'')-\pp^Z_a(Z_a')\bigr)\cdot\bigl(Z_a''-Z_a'\bigr)\ \ge\ 0
	\qquad\text{for any two states } Z_a',Z_a'' .
\end{equation}
Unlike strict convexity, this inequality allows distinct states to have the same thermodynamic forces.

\subsection{The configurational direction}\label{subsec:configurational_direction}

The flat direction has a simple physical meaning. A general variation of \(X_a=\alpha_aY_a\) satisfies
\[
    \delta X_a=Y_a\,\delta\alpha_a+\alpha_a\,\delta Y_a.
\]
We now consider variations at fixed intrinsic phase state, \(\delta Y_a=0\): the intrinsic density, phase velocity and specific entropy are held fixed, and only the amount of phase \(a\) changes. Then \(\delta X_a=Y_a\delta\alpha_a\), and the variation of the full phase vector is
\begin{equation}\label{eq:configurational_direction}
    \delta Z_a=(\delta\alpha_a,\ \delta X_a)=(1,\ Y_a)\,\delta\alpha_a.
\end{equation}
This variation changes how much of phase \(a\) is present and leaves what the phase is unchanged: \(Y_a=X_a/\alpha_a\) is constant along it. We call \eqref{eq:configurational_direction} the \emph{configurational direction} of phase \(a\). Along it the energy changes linearly, \(\rhoE_a\bigl(Z_a+\delta\alpha_a(1,Y_a)\bigr)=(\alpha_a+\delta\alpha_a)E_a(Y_a)\), and the forces \eqref{eq:dual_Z} do not change at all.

This observation explains both a difficulty and its remedy. The volume fraction is a geometric quantity, transported with the material velocity; the remaining components of \(Z_a\) also carry the acoustic waves. A numerical dissipation that treats the geometry and the acoustics with different coefficients is therefore desirable, but it cannot act on \(\alpha_a\) and on \(X_a\) independently. Indeed, differentiation of \(Y_a=X_a/\alpha_a\) gives
\begin{equation}\label{eq:decoupling}
    \delta Y_a=\frac{\delta X_a-Y_a\,\delta\alpha_a}{\alpha_a}.
\end{equation}
If the two variations are chosen independently, their mismatch is divided by the volume fraction. Where a phase is present only in traces, a small mismatch between the two dissipations would produce a large change of its density, momentum or entropy. Equation \eqref{eq:decoupling} also shows the remedy: a change of \(\alpha_a\) accompanied by the change \(\delta X_a=Y_a\delta\alpha_a\) of the extensive variables leaves the intrinsic state untouched.

We therefore distinguish a configurational contribution, which changes the phase amount, from an intrinsic contribution, which changes \(Y_a\). This allows different dissipative coefficients for the two contributions without discarding the product relation \(X_a=\alpha_aY_a\). Section~\ref{subsec:jump_decomposition} gives the corresponding decomposition for a finite jump across a cell face; the conditions on the two coefficients are discussed in Section~\ref{subsec:two_scales}.

\subsection{One energy, three representations}\label{subsec:three_representations}

The energy evaluated at a physical state is independent of the chosen representation, but its convexity is not. Table~\ref{tab:representations} summarizes the properties used here; further details are given in \ref{app:convexity}.

\begin{table}[t]
\centering
\small
\caption{The same energy in three sets of variables. ``Internal'' refers to \(\sum_a\crho_ae_a\), ``total'' to \(\rhoE\).}
\label{tab:representations}
\renewcommand{\arraystretch}{1.25}
\begin{tabularx}{\textwidth}{@{}>{\raggedright}p{2.9cm}>{\raggedright}X>{\raggedright}X>{\raggedright\arraybackslash}X@{}}
\toprule
Variables & Role & Internal energy & Total energy \\
\midrule
\addlinespace[2pt]
\(\qq^\varphi\) (with \(\varphi_a\), \(\ww_a\)) & derivation of the model; discrete equivalence (Section~\ref{sec:equivalence}) & not convex in general & not convex in general; indefinite for \(\ww_a\neq\mathbf 0\) at pressure equilibrium \\
\addlinespace[6pt]
\(\qq\) (with \(\alpha_a\), \(\ww_a\)) & evolved variables  & convex (sum of perspectives) & not convex for \(\ww_a\neq\mathbf 0\) \\
\addlinespace[6pt]
\(Z\) (phase variables) & dissipation and correction direction & convex & convex, homogeneous of degree one \\
\addlinespace[2pt]
\bottomrule
\end{tabularx}
\end{table}

\paragraph{From \(\varphi_a\) to \(\alpha_a\)}
In the \(\alpha\)-formulation, the internal energy reads \(\sum_a\alpha_a\varepsilon_a(\crho_a/\alpha_a,\eta_a/\alpha_a)\). It is a sum of convex perspectives in \((\alpha_a,\crho_a,\eta_a)\), also after inserting the affine relations \(\alpha_\Nph=1-\sum_{b<\Nph}\alpha_b\) and \(\crho_\Nph=\rho-\sum_{b<\Nph}\crho_b\), hence convex.
In the original formulation, \(\qq^\varphi\), the substitution \(\alpha_a=\varphi_a/\rho\) is nonlinear. It introduces an additional Hessian term involving the pressure differences. At pressure equilibrium the internal-energy Hessian remains positive semi-definite; for two phases it is indefinite whenever \((p_1-p_2)(\rho_1-\rho_2)<0\), as shown in \ref{app:phi_alpha}.

\paragraph{The relative velocities}
The use of \(\alpha_a\) thus recovers convexity of the internal energy, but not of the total energy when the relative velocities are retained as independent variables. The following result makes this distinction explicit.

\begin{proposition}[Non-convexity in the evolved variables]\label{prop:nonconvex_q}
Consider two phases with positive partial densities and \(\ww_1\neq\mathbf 0\). Then the Hessian of \(\rhoE\) with respect to \(\qq\) is indefinite.
\end{proposition}

\begin{proof}
Let \(\bm n\) be the variation obtained by increasing \(\alpha_1\) and decreasing \(\alpha_2\) at fixed intrinsic states and phase velocities. Per unit change of \(\alpha_1\), its entries in \(\qq\) are
\[
\begin{gathered}
    \delta\crho_1=\rho_1,\qquad
    \delta\crho_2=-\rho_2,\qquad
    \delta\UU=\rho_1\vvv_1-\rho_2\vvv_2,\qquad
    \delta\ww_1=\mathbf 0,\\
    \delta\eta_1=\rho_1s_1,\qquad
    \delta\eta_2=-\rho_2s_2,\qquad
    \delta\rho=\rho_1-\rho_2.
\end{gathered}
\]
Here \(\rho_1\) and \(\rho_2\) are intrinsic densities. Along this line,
\[
    \rhoE(\qq+t\bm n)
    =(\alpha_1+t)E_1(Y_1)+(\alpha_2-t)E_2(Y_2),
\]
so that \(\bm n^{\mathsmaller T}\nabla^2_{\qq}\rhoE\,\bm n=0\). Now let \(\bm e\) vary only \(\UU\), with unit length in the direction of \(\ww_1\). Since \(\rhoE_{\UU}=\VV\),
\[
    \bm e^{\mathsmaller T}\nabla^2_{\qq}\rhoE\,\bm e=\frac1\rho,
    \qquad
    \bm n^{\mathsmaller T}\nabla^2_{\qq}\rhoE\,\bm e
    =\partial_{\bm n}\VV\cdot\bm e
    =\frac{\rho_1\rho_2}{\rho^2}|\ww_1|=:b.
\]
The quadratic form restricted to the span of \(\bm n\) and \(\bm e\) is represented by
\[
	\begin{pmatrix}0&b\\ b&1/\rho\end{pmatrix},
	\qquad
	\det=-b^2<0,
\]
which is indefinite as soon as \(\ww_1\neq\mathbf 0\).
\end{proof}

The result holds for any nonzero relative velocity: an arbitrarily small \(\ww_1\) suffices. The proof also shows another point. The configurational direction has zero curvature in both representations, but only in \(Z\) does it belong to the kernel of the Hessian. In \(\qq\), for \(\ww_1\neq\mathbf 0\), it is coupled to the mixture momentum.

\subsection{Implication for numerical dissipation}\label{subsec:implications_dissipation}
For the numerical method, we must relate the thermodynamic conjugates in the two representations.
Because \(\qq\) is a function of \(Z\) and the energy is the same function of the physical state, the forces of the two representations are related at every state, through the Jacobian \(\mathbf J=\partial\qq/\partial Z\), by the chain rule on admissible variations,
\begin{equation}\label{eq:chain_rule_admissible}
	\pp\cdot\mathbf J\,\delta Z=\sum_a\pp^Z_a\cdot\delta Z_a
	\qquad\text{for every }\delta Z\text{ with }\textstyle\sum_a\delta\alpha_a=0 .
\end{equation}
Componentwise, \(\mathbf J^{\mathsmaller T}\pp=\pp^Z+p_\Nph\,\mathbf 1_\alpha\), where \(\mathbf 1_\alpha\) has unit entries on the volume fractions and zero entries elsewhere. The additional term is normal to the saturation constraint: its contraction with an admissible variation is \(p_\Nph\sum_a\delta\alpha_a=0\). Hence both sets of conjugates give the same variation of energy.

We must now distinguish these infinitesimal variations from the finite
jumps used by a numerical scheme. For two states \(L\) and \(R\), \(\Delta\qq\) is generally neither \(\mathbf J_L\Delta Z\) nor \(\mathbf J_R\Delta Z\).  Applying \eqref{eq:chain_rule_admissible} at each endpoint gives
\[
    \sum_a\Delta\pp^Z_a\cdot\Delta Z_a
    =\pp_R\cdot\mathbf J_R\Delta Z
     -\pp_L\cdot\mathbf J_L\Delta Z.
\]
Subtracting this expression from \(\Delta\pp\cdot\Delta\qq\) yields
\begin{equation}\label{eq:pairing_split}
	\Delta\pp\cdot\Delta\qq
	=
	\sum_a\Delta\pp^Z_a\cdot\Delta Z_a+T_w,
	\qquad
	\Delta(\cdot)=(\cdot)_R-(\cdot)_L ,
\end{equation}
where
\[
    T_w=\pp_R\cdot(\Delta\qq-\mathbf J_R\Delta Z)
       -\pp_L\cdot(\Delta\qq-\mathbf J_L\Delta Z).
\]
The linear rows of \eqref{eq:q_of_Z} make no contribution to this remainder. It comes only from the relative velocities, vanishes when all \(\ww_a\) vanish at both endpoints, and otherwise has no prescribed sign. Its explicit expression is given in \ref{app:pairing}. It is quadratic in small jumps, just as the pairing itself, and therefore cannot in general be neglected on the basis of the jump size alone.

\paragraph{The pairings as integrals of the Hessian}
Both pairings in \eqref{eq:pairing_split} can be read along the straight segment between the two states, as in the HTC literature \cite{HTC2022,ABGRALL2023127629}.
Consider first the evolved variables,
\[
    \qq(s)=\qq_L+s\,\Delta\qq,\qquad 0\le s\le1.
\]
Since \(\pp=\nabla_{\qq}\rhoE\), the variation of the forces along this
segment is
\[
    \frac{d}{ds}\pp(\qq(s))
    =\nabla^2_{\qq}\rhoE(\qq(s))\,\Delta\qq,
\]
and integration gives
\begin{equation}\label{eq:roe_hessian}
	\Delta\pp=\widetilde{\mathbf H}_{\qq}\,\Delta\qq,
	\qquad
	\widetilde{\mathbf H}_{\qq}=\int_0^1\nabla^2_{\qq}\rhoE\bigl(\qq(s)\bigr)\,ds,
	\qquad
	\Delta\pp\cdot\Delta\qq=\int_0^1\Delta\qq\cdot\nabla^2_{\qq}\rhoE\bigl(\qq(s)\bigr)\,\Delta\qq\,ds .
\end{equation}
The matrix \(\widetilde{\mathbf H}_{\qq}\) is the \emph{Roe matrix of the Hessian} of the energy: it maps the jump of the state exactly to the jump of its forces. The second identity gives a direct
interpretation of the pairing, namely the curvature of the energy in the direction of the jump, accumulated along the segment. For a convex energy this curvature is
non-negative. In the evolved variables of our model, however,
Proposition~\ref{prop:nonconvex_q} shows that this sign is not guaranteed
when the phases slip.

The same construction in the phase variables, along \(Z_a(s)=Z_{a,L}+s\,\Delta Z_a\), gives \[
    \Delta\pp_a^Z\cdot\Delta Z_a
    =
    \int_0^1
    \Delta Z_a\cdot\nabla^2\rhoE_a(Z_a(s))
    \,\Delta Z_a\,ds.
\]
Here the sign follows from the perspective structure
\(\rhoE_a(Z_a)=\alpha_a E_a(X_a/\alpha_a)\).
Indeed, differentiating \(X_a=\alpha_aY_a\) gives
\[
    \delta X_a-Y_a\,\delta\alpha_a=\alpha_a\,\delta Y_a,
\]
and the second variation of the perspective energy is
\begin{equation}\label{eq:perspective_hessian}
    \begin{aligned}
    \delta Z_a\cdot\nabla^2\rhoE_a(Z_a)\,\delta Z_a
    &=
    \frac{1}{\alpha_a}
    \bigl(\delta X_a-Y_a\,\delta\alpha_a\bigr)
    \cdot\nabla^2E_a(Y_a)
    \bigl(\delta X_a-Y_a\,\delta\alpha_a\bigr)
    \\
    &=
    \alpha_a\,\delta Y_a\cdot\nabla^2E_a(Y_a)\,\delta Y_a,
    \qquad \alpha_a>0.
    \end{aligned}
\end{equation}
Three properties can be read directly from \eqref{eq:perspective_hessian}. The integrand is non-negative. It vanishes along the configurational direction \(\delta X_a=Y_a\,\delta\alpha_a\) of Section~\ref{subsec:configurational_direction}. And it is the curvature of the single-fluid energy in the direction of the intrinsic change \(\delta Y_a\), weighted by the amount of phase \(\alpha_a\); equivalently, a given mismatch \(\delta X_a-Y_a\delta\alpha_a\) of the extensive variables is amplified by \(1/\alpha_a\), which is the decoupling issue of Section~\ref{subsec:configurational_direction}. The segment between two admissible states stays admissible, since \(\alpha_a(s)>0\) and the extensive variables are convex combinations of positive quantities.

The first term on the right of \eqref{eq:pairing_split} is non-negative by convexity in \(Z\). The full pairing in \(\qq\) need not be. This is the difficulty with a scalar viscosity proportional to \(\Delta\qq\), whose entropy production depends on \(\Delta\pp\cdot\Delta\qq\), and \eqref{eq:pairing_split} shows that this quantity can be negative. We will instead construct the viscous exchange in \(Z\) and transfer it through \(\mathbf J_L\) and \(\mathbf J_R\), so that the energy pairing is the one in the phase variables. The conservative rows receive the same face contribution on both sides, whereas the relative velocity rows need not. This construction is developed in Section~\ref{sec:viscosity}; the same side-wise chain rule is used for the correction direction in Section~\ref{subsec:abgrall}.

The conclusion of this section is the guiding principle of the method: \emph{build the scheme in \(\qq\), dissipate in \(Z\)}. The fluxes and the non-conservative products stay in \(\qq\) (Section~\ref{sec:framework}). The two ingredients that need a sign, the denominator of the thermodynamic correction and the numerical dissipation, are constructed in the phase variables, where the energy is convex and the configurational direction is explicit, and are then transferred to \(\qq\) (Sections~\ref{subsec:abgrall} and~\ref{sec:viscosity}); from \(\qq\), the same transfer reaches the original variables \(\qq^\varphi\) (Section~\ref{sec:equivalence}).

\section{Semi-discrete HTC framework}\label{sec:framework}
This section presents the semi-discrete finite volume scheme. We follow the HTC framework \cite{ABGRALL2018640,ABGRALL2023127629,HTC2022}, which imposes three discrete conditions: one on the reversible fluxes and non-conservative products, one on the relaxation sources and one on the numerical dissipation.
The first two are presented in this section. For the first, we change the direction of the standard thermodynamic correction: its denominator is non-negative when the energy is convex in the evolved variables, which is not the case here (Section~\ref{subsec:three_representations}), and the new direction restores its sign. The third condition is formulated here as a list of requirements and is the subject of Section~\ref{sec:viscosity}. As in the continuous model, where the reversible part conserves the energy and the irreversible part produces entropy, the two discrete ingredients have distinct roles: the thermodynamic correction restores the energy identity of the reversible part, which the uncorrected discretization violates, and the entropy production accounts for the energy that the dissipation converts into heat.

\subsection{Cell and face notation}\label{subsec:notation}

The domain is covered by orthogonal control volumes \(\Omega_i\) with centres \(\x_i\). For a neighbour \(j\in\mathcal N_i\), \(\partial\Omega_{ij}\) is the common face, \(\mathbf n_{ij}=-\mathbf n_{ji}\) its unit normal pointing from \(\Omega_i\) to \(\Omega_j\), and \(|\bm\delta_{ij}|=|\x_j-\x_i|\) the distance of the centres. Face quantities are projected on \(\mathbf n_{ij}\): we write \(\f_{\mathbf n}=\f_kn_{ij,k}\), \(\B_{\mathbf n}=\B_kn_{ij,k}\) and \(\mathrm F^{\rhoE}_{\mathbf n}\) for the normal component of the energy flux in \eqref{eq:energy_conservation}, and \(\Delta(\cdot)_{ij}=(\cdot)_j-(\cdot)_i\). The thermodynamic forces of cell \(i\) are \(\pp_i=\partial_{\qq}\rhoE(\qq_i)\).

With the system of Section~\ref{subsec:why_alpha} written as \(\partial_t\qq+\partial_k\f_k(\qq)+\B_k(\qq)\,\partial_k\qq=\S(\qq)\), the semi-discrete scheme reads
\begin{equation}\label{eq:semidiscrete}
	\frac{d\qq_i}{dt}
	=
	-\frac{1}{|\Omega_i|}\sum_{j\in\mathcal N_i}|\partial\Omega_{ij}|
	\Bigl(
		\bm\Phi_{ij}-\mathcal G^{(i)}_{ij}-\mathcal P^{(i)}_{ij}
	\Bigr)
	+\S(\qq_i).
\end{equation}
Each face contributes three terms to each of its two cells:
\begin{itemize}
	\item the \emph{reversible fluctuation} \(\bm\Phi_{ij}=\mathcal F_{ij}-\f_{\mathbf n}(\qq_i)+\mathcal B_{ij}\), built from a numerical flux \(\mathcal F_{ij}\) and a non-conservative fluctuation \(\mathcal B_{ij}\) (Sections~\ref{subsec:fluxes} and~\ref{subsec:abgrall});
	\item the \emph{dissipative contribution} \(\mathcal G^{(i)}_{ij}\), the numerical viscosity seen by cell \(i\);
	\item the \emph{entropy production} \(\mathcal P^{(i)}_{ij}\), whose only non-zero components are in the phase entropy equations.
\end{itemize}
The superscript \((i)\) emphasizes that the dissipative contribution and the production are assigned to a cell. For the rows of conserved quantities the dissipative contribution is a flux, \(\mathcal G^{(j)}_{ji}=-\mathcal G^{(i)}_{ij}\); for the non-conservative rows (relative velocities and volume fractions) the two cells may receive different contributions. This freedom will be essential in Sections~\ref{subsec:abgrall}, \ref{sec:viscosity} and~\ref{sec:equivalence}.

\subsection{Conservative fluxes and non-conservative contributions}\label{subsec:fluxes}

\paragraph{Face-wise energy condition}
Contracting the reversible part of \eqref{eq:semidiscrete} with \(\pp_i\), we obtain the evolution of the energy of cell \(i\),
\[
	\frac{d}{dt}\rhoE(\qq_i)=-\frac{1}{|\Omega_i|}\sum_{j\in\mathcal N_i}|\partial\Omega_{ij}|\,\pp_i\cdot\bm\Phi_{ij}.
\]
For each face we define the energy flux seen from cell \(i\),
\begin{equation}\label{eq:reversible_energy_flux}
	\mathcal F^{\rhoE}_{ij}:=\mathrm F^{\rhoE}_{\mathbf n}(\qq_i)+\pp_i\cdot\bm\Phi_{ij},
\end{equation}
namely the physical energy flux of the cell plus the pairing of its forces with the fluctuation. Since \(\sum_j|\partial\Omega_{ij}|\,\mathbf n_{ij}=\mathbf 0\) for a closed cell, the terms \(\mathrm F^{\rhoE}_{\mathbf n}(\qq_i)\) cancel after summation over the faces, and
\[
	\frac{d}{dt}\rhoE(\qq_i)+\frac{1}{|\Omega_i|}\sum_{j\in\mathcal N_i}|\partial\Omega_{ij}|\,\mathcal F^{\rhoE}_{ij}=0 .
\]
This is a discrete conservation law provided the flux is single-valued, \(\mathcal F^{\rhoE}_{ji}=-\mathcal F^{\rhoE}_{ij}\): the energy that leaves cell \(i\) through the face is the energy that enters cell \(j\). Since \(\mathbf n_{ji}=-\mathbf n_{ij}\), \(\mathcal F^{\rhoE}_{ji}=-\mathrm F^{\rhoE}_{\mathbf n}(\qq_j)+\pp_j\cdot\bm\Phi_{ji}\), and the condition is equivalent to
\begin{equation}\label{eq:face_energy_condition}
	\pp_i\cdot\bm\Phi_{ij}+\pp_j\cdot\bm\Phi_{ji}
	=
	\mathrm F^{\rhoE}_{\mathbf n}(\qq_j)-\mathrm F^{\rhoE}_{\mathbf n}(\qq_i).
\end{equation}
Thus, the energy fluctuations generated on the two sides of a face must combine into the jump of the physical energy flux. This is the discrete counterpart of the identity \(\pp\cdot(\partial_k\f_k+\B_k\partial_k\qq)=\partial_k\mathrm F^{\rhoE}_k\) that yields \eqref{eq:energy_conservation}.
When \eqref{eq:face_energy_condition} is satisfied, the numerical energy flux \eqref{eq:reversible_energy_flux} is single-valued and the reversible scheme satisfies a local semi-discrete energy balance.
Moreover, this flux need not be computed: it is implied by the fluctuations, and the total energy, which is not evolved, is conserved as a consequence.

\paragraph{Central flux and path-conservative products}
We now specify the preliminary numerical flux and the fluctuations
associated with the non-conservative products. We use a central flux
and approximate the path integral along the straight segment between
the two states by the midpoint rule \cite{Pares2006,Castro2006}
\begin{equation}\label{eq:central_ncp}
	\widetilde{\mathcal F}_{ij}=\tfrac12\bigl(\f_{\mathbf n}(\qq_i)+\f_{\mathbf n}(\qq_j)\bigr),
	\qquad
	\mathcal B_{ij}=\tfrac12\,\B_{\mathbf n}(\widetilde\qq_{ij})\,\Delta\qq_{ij},
	\qquad
	\widetilde\qq_{ij}=\tfrac12(\qq_i+\qq_j).
\end{equation}
These choices give the same uncorrected fluctuation for both the cells
\[
    \widetilde{\bm\Phi}_{ij}
    =
    \widetilde{\bm\Phi}_{ji}
    =
    \tfrac12\Delta\f_{\mathbf n,ij}+\mathcal B_{ij},
\]
namely half of the jump of the physical flux and half of the non-conservative product. In the volume fraction formulation, the non-conservative products account for the advection \(V_k\partial_k\alpha_a\) and the rotational term of the relative velocity equations.

The central flux has a thermodynamic interpretation \cite{Tadmor1987,ABGRALL2018640}. For a conservative system, \eqref{eq:face_energy_condition} states that the numerical flux is a discrete gradient, with respect to the forces, of the potential \(\mathcal K_{\mathbf n}(\pp)=\pp\cdot\f_{\mathbf n}-\mathrm F^{\rhoE}_{\mathbf n}\), which satisfies \(d\mathcal K_{\mathbf n}=\f_{\mathbf n}\cdot d\pp\). The exact flux would be the average of \(\f_{\mathbf n}\) along the straight segment in the space of forces; the central flux is its trapezoidal approximation. Neither this quadrature nor the midpoint rule for \(\mathcal B_{ij}\) satisfies \eqref{eq:face_energy_condition} exactly, and a correction is needed.

\subsection{Thermodynamic correction of the reversible part}\label{subsec:abgrall}
With the central flux and the midpoint rule \eqref{eq:central_ncp}, both cells receive the same uncorrected fluctuation, \(\widetilde{\bm\Phi}_{ji}=\widetilde{\bm\Phi}_{ij}\), and the face-wise energy condition \eqref{eq:face_energy_condition} is in general violated by the residual
\begin{equation}\label{eq:abgrall_residual}
	\mathcal R_{ij}
	=
	\mathrm F^{\rhoE}_{\mathbf n}(\qq_j)-\mathrm F^{\rhoE}_{\mathbf n}(\qq_i)
	-(\pp_i+\pp_j)\cdot\Bigl(\tfrac12\Delta\f_{\mathbf n,ij}+\mathcal B_{ij}\Bigr).
\end{equation}
Following \cite{ABGRALL2018640,ABGRALL2023127629}, we introduce
one scalar coefficient \(\theta_{ij}=\theta_{ji}\) per face to
compensate this residual. Correction terms proportional to
\(\theta_{ij}\) are added to the two fluctuations, and the
coefficient is determined by imposing \eqref{eq:face_energy_condition}.

\paragraph{The correction in the evolved variables}
The original construction uses the jump of the state \(\Delta\qq_{ij}\) as direction. The fluctuations become
\[
    \bm\Phi_{ij}
    =\widetilde{\bm\Phi}_{ij}
     -\tfrac12\theta_{ij}\Delta\qq_{ij},
    \qquad
    \bm\Phi_{ji}
    =\widetilde{\bm\Phi}_{ij}
     +\tfrac12\theta_{ij}\Delta\qq_{ij},
\]
i.e. the numerical flux becomes \(\widetilde{\mathcal F}_{ij}-\frac12\theta_{ij}\Delta\qq_{ij}\). Substitution into \eqref{eq:face_energy_condition} then requires
\[
    \tfrac12\theta_{ij}\,
    \Delta\pp_{ij}\cdot\Delta\qq_{ij}
    =\mathcal R_{ij},
    \qquad\text{i.e.}\qquad
    \theta_{ij}=\frac{2\,\mathcal R_{ij}}{\Delta\pp_{ij}\cdot\Delta\qq_{ij}} .
\]
For a strictly convex energy, the denominator of \(\theta_{ij}\) is positive whenever the two states differ. However, this property is not available in the evolved variables of our multiphase model, where the energy is not convex as soon as the phases slip (Proposition~\ref{prop:nonconvex_q}). By \eqref{eq:pairing_split}, the denominator is the sum of the non-negative configurational pairing \(\sum_a\Delta\pp^Z_{a,ij}\cdot\Delta Z_{a,ij}\) and of the remainder \(T_w\), which has no prescribed sign. The two contributions may therefore cancel even when the intrinsic phase states differ and the energy residual is nonzero. Near such a cancellation, \(\theta_{ij}\) can become large. Since \(|\theta_{ij}|\) enters the dissipative speeds (Section~\ref{subsec:two_scales}), this may severely restrict the time step, as observed in Section~\ref{subsec:test_shockcolumn}.

A negative denominator does not by itself prevent compatibility; indeed, the coefficient is signed and adjusts to the residual. The difficulty is the possible vanishing of the denominator without a corresponding vanishing of \(\mathcal R_{ij}\).

\paragraph{A direction built in the phase variables}
We now change the direction of the correction, while keeping the
uncorrected fluctuations and their residual \(\mathcal R_{ij}\)
unchanged. Instead of using \(\Delta\qq_{ij}\), we start from the
jump of the phase variables,
\[
    \Delta Z_{ij}=Z_j-Z_i.
\]
This jump, which is an admissible variation, can be transferred to each cell by the Jacobian
\(\mathbf J_s=\partial\qq/\partial Z\) evaluated at that cell. We define
\begin{equation}\label{eq:abgrall_flux}
	\bm\Phi_{ij}=\widetilde{\bm\Phi}_{ij}-\tfrac12\,\theta_{ij}\,\mathbf J_i\,\Delta Z_{ij},
	\qquad
	\bm\Phi_{ji}=\widetilde{\bm\Phi}_{ij}+\tfrac12\,\theta_{ij}\,\mathbf J_j\,\Delta Z_{ij}.
\end{equation}
There is still one scalar coefficient per face. The two directions,
however, need not coincide: they describe the change induced by the
same phase-space direction at two different states.

On the rows of \(\UU\), \(\crho_a\), \(\eta_a\), \(\rho\) and \(\alpha_a\), the map \eqref{eq:q_of_Z} is linear and both directions coincide with the original one
\[\mathbf J_i\Delta Z_{ij}=\mathbf J_j\Delta Z_{ij}=\Delta\qq_{ij}.\] The correction therefore retains the flux form
\(\widetilde{\mathcal F}_{ij}
-\frac12\theta_{ij}\Delta\qq_{ij}\)
on the conserved rows. Their conservation is unchanged, although
the value of the coefficient will in general differ from that of
the correction along \(\Delta\qq_{ij}\).

Only the relative velocity rows receive different directions.
For a phase velocity \(\vvv_a=\bm M_a/\crho_a\), differentiation
at the state of cell \(s\) gives
\[
    D_s\vvv_a[\Delta Z_{ij}]
    =
    \frac{\Delta\bm M_{a,ij}
          -\vvv_{a,s}\Delta\crho_{a,ij}}{\crho_{a,s}}.
\]
Using \(\bm M_a=\crho_a\vvv_a\), this becomes
\begin{equation}\label{eq:abgrall_w}
    D_i\vvv_a[\Delta Z_{ij}]
    =\frac{\crho_{a,j}}{\crho_{a,i}}\Delta\vvv_{a,ij},
    \qquad
    D_j\vvv_a[\Delta Z_{ij}]
    =\frac{\crho_{a,i}}{\crho_{a,j}}\Delta\vvv_{a,ij}.
\end{equation}
The direction on the row of
\(\ww_a=\vvv_a-\vvv_\Nph\) is then
\(D_s\vvv_a[\Delta Z_{ij}]
-D_s\vvv_\Nph[\Delta Z_{ij}]\).
For small jumps around a state with positive partial densities,
these directions differ from the corresponding finite jumps by
second-order terms. Their different values on the two sides are
consistent with the non-conservative form of the relative velocity
equations. The same transfer will be used for the numerical
dissipation in Section~\ref{subsec:transfer}.

The coefficient \(\theta_{ij}\) is again fixed by the face-wise energy condition,
not prescribed independently of the direction.
Substituting \eqref{eq:abgrall_flux} into
\eqref{eq:face_energy_condition}, the correction contributes
\[
    \tfrac12\theta_{ij}
    \left(
        \pp_j\cdot\mathbf J_j\Delta Z_{ij}
        -\pp_i\cdot\mathbf J_i\Delta Z_{ij}
    \right).
\]
At each endpoint, the chain rule
\eqref{eq:chain_rule_admissible} gives
\[
    \pp_s\cdot\mathbf J_s\Delta Z_{ij}
    =
    \sum_a\pp^Z_{a,s}\cdot\Delta Z_{a,ij}.
\]
Consequently, the energy condition reduces to
\begin{equation}\label{eq:abgrall_theta}
    \tfrac12\theta_{ij}\delta_{ij}=\mathcal R_{ij},
    \qquad
    \delta_{ij}
    :=
    \sum_a\Delta\pp^Z_{a,ij}\cdot\Delta Z_{a,ij},
    \qquad
    \theta_{ij}=\frac{2\mathcal R_{ij}}{\delta_{ij}}
    \quad\text{if }\delta_{ij}>0.
\end{equation}
Thus, the pairing in the phase variables appears as a consequence
of transferring the direction through the two Jacobians.
By the integral representation of
Section~\ref{subsec:implications_dissipation},
\[
    \delta_{ij}
    =
    \sum_a\int_0^1
    \Delta Z_{a,ij}\cdot
    \nabla^2\rhoE_a\bigl(Z_a(s)\bigr)
    \,\Delta Z_{a,ij}\,ds
    \ge0,
    \qquad
    Z_a(s)=Z_{a,i}+s\,\Delta Z_{a,ij}.
\]
Each phase contributes a non-negative curvature integral.
The denominator therefore cannot vanish through the cancellation
between the pairing in \(Z\) and \(T_w\) that occurs in
\(\Delta\pp_{ij}\cdot\Delta\qq_{ij}\).
It vanishes only along the flat directions of the perspective energy, that is, when every phase has the same intrinsic state on the two sides of the face (see the third remark below).

Equivalently, using the Bregman divergences introduced in Section~\ref{subsec:viscous_requirements},
\[
    \delta_{ij}
    =
    \sum_a\left[
        \alpha_{a,j}
        D_{E_a}(Y_{a,j}\mid Y_{a,i})
        +
        \alpha_{a,i}
        D_{E_a}(Y_{a,i}\mid Y_{a,j})
    \right].
\]
This expression is evaluated directly using the formulas of
\ref{app:bregman}; no numerical quadrature of the Hessian
is required.

Three remarks clarify the role of the correction.
\begin{enumerate}
	\item \emph{The correction remains reversible.} The coefficient \(\theta_{ij}\) has the sign that the residual requires; the correction may add or remove numerical diffusion. It restores an identity of the reversible part and produces no entropy.
	\item \emph{The coefficient is a signed speed.} Row by row, \(\mathbf J_s\Delta Z_{ij}\) has the dimensions of \(\Delta\qq_{ij}\), so \(\theta_{ij}\) has the dimension of a velocity, and on the conserved rows it enters the flux exactly as the coefficient of a Rusanov term. The magnitude \(|\theta_{ij}|\) will therefore appear in the scales of the numerical dissipation (Section~\ref{subsec:two_scales}) and in the positivity conditions (Section~\ref{subsec:positivity}).
	\item \emph{The denominator \(\delta_{ij}\) vanishes only with the intrinsic jumps.} For strictly convex single-fluid energies, \(\delta_{ij}=0\) if and only if \(Y_{a,i}=Y_{a,j}\) for every phase. On a contact, where the velocity and the pressure are also common, the residual vanishes as well and no correction is needed (Theorem~\ref{thm:contact}).
	With interphase slip, however, identical intrinsic states on
the two sides do not in general make the central discretization
energy-compatible: \(\delta_{ij}=0\) may then coexist with
\(\mathcal R_{ij}\ne0\); the pairing in \(\qq\) vanishes there as well. The phase-variable construction removes the indefinite
contribution \(T_w\), but not this common degeneracy.
\end{enumerate}

\subsection{Relaxation sources}\label{subsec:sources}

The relaxation sources are evaluated pointwise, \(\S_i=\S(\qq_i)\). Since in the continuous model the entropy productions \eqref{eq:entropy_production} exactly compensate the energy removed by friction and pressure relaxation, the same identity holds cell by cell,
\begin{equation}\label{eq:source_compatibility}
	\pp_i\cdot\S(\qq_i)
	=
	\rhoE_{\ww_a}\cdot\bm\Lambda_a+\rhoE_{\alpha_a}\Psi_a+\sum_aT_a\Pi_a
	=
	-\mathcal D_i+\mathcal D_i
	=
	0 .
\end{equation}
Nothing special is needed at the discrete level. In the computations the sources are treated implicitly; the identity \eqref{eq:source_compatibility} holds at whatever state they are evaluated.

\subsection{Energy compatibility of the dissipative part}\label{subsec:viscous_requirements}

The dissipative contributions of a face must not create or destroy total energy; they must move it between the two cells through a single-valued flux, and convert the rest into entropy.
As for the reversible part, we therefore require a single-valued
numerical energy flux \(\mathcal G^{\rhoE}_{ij}=-\mathcal G^{\rhoE}_{ji}\) and define the \emph{energy defects} of the two cells,
\begin{equation}\label{eq:energy_defects}
	\Gamma^{(i)}_{ij}:=\mathcal G^{\rhoE}_{ij}-\pp_i\cdot\mathcal G^{(i)}_{ij},
	\qquad
	\Gamma^{(j)}_{ji}:=\mathcal G^{\rhoE}_{ji}-\pp_j\cdot\mathcal G^{(j)}_{ji}.
\end{equation}
The defect is the energy carried by the viscous energy flux that the work of the viscous state contribution does not account for. It is assigned to the entropy of the cell,
\begin{equation}\label{eq:production_requirement}
	\pp_i\cdot\mathcal P^{(i)}_{ij}=\sum_{a=1}^{\Nph}T_{a,i}\,\Pi^{(i)}_{a,ij}=\Gamma^{(i)}_{ij},
\end{equation}
where \(\Pi^{(i)}_{a,ij}\) is the production of phase entropy \(\eta_a\) in cell \(i\). Since the temperatures are positive, non-negative productions exist if and only if \(\Gamma^{(i)}_{ij}\ge0\).

\begin{theorem}[Discrete energy conservation and entropy production]\label{thm:discrete_energy}
Assume periodic or impermeable boundaries. Let the reversible fluctuations satisfy \eqref{eq:face_energy_condition}, the sources \eqref{eq:source_compatibility}, and the dissipative contributions be conservative on the rows of \(\UU\), \(\crho_a\), \(\eta_a\), \(\rho\), with productions satisfying \eqref{eq:production_requirement}. Then every cell satisfies
\begin{equation}\label{eq:discrete_energy}
	\frac{d}{dt}\rhoE(\qq_i)
	+
	\frac{1}{|\Omega_i|}\sum_{j\in\mathcal N_i}|\partial\Omega_{ij}|\left(\mathcal F^{\rhoE}_{ij}-\mathcal G^{\rhoE}_{ij}\right)=0,
\end{equation}
and the total energy \(\sum_i|\Omega_i|\rhoE(\qq_i)\) is conserved. If, moreover, all energy defects are non-negative and the productions \(\Pi^{(i)}_{a,ij}\) are chosen non-negative, the mixture entropy \(\eta_i=\sum_a\eta_{a,i}\) satisfies a cell balance with a conservative numerical entropy flux and a non-negative production, and \(\sum_i|\Omega_i|\eta_i\) is non-decreasing.
\end{theorem}

\begin{proof}
Contract \eqref{eq:semidiscrete} with \(\pp_i\). The reversible part gives \(\mathcal F^{\rhoE}_{ij}-\mathrm F^{\rhoE}_{\mathbf n}(\qq_i)\) on each face; the terms \(\mathrm F^{\rhoE}_{\mathbf n}(\qq_i)\) cancel after summation over the faces of the cell, because \(\sum_j|\partial\Omega_{ij}|\,\mathbf n_{ij}=\mathbf 0\). The dissipative part gives \(\pp_i\cdot(\mathcal G^{(i)}_{ij}+\mathcal P^{(i)}_{ij})=\mathcal G^{\rhoE}_{ij}\) by \eqref{eq:energy_defects}--\eqref{eq:production_requirement}, and the sources give zero. This is \eqref{eq:discrete_energy}; summing over the cells, the single-valued fluxes cancel in pairs. For the entropy, sum the phase entropy rows of \eqref{eq:semidiscrete}: the reversible and dissipative contributions to these rows are conservative and form the numerical entropy flux, while the remaining terms are the productions \(\sum_a\Pi_{a,i}\ge0\) of the relaxation sources and \(\sum_a\Pi^{(i)}_{a,ij}\ge0\) of the faces.
\end{proof}

Note that \eqref{eq:energy_defects} does not require the two cells to receive
opposite dissipative contributions on the non-conservative
rows. Only the numerical energy flux must be single-valued.
This distinction will be used when transferring the
viscosity from the phase variables to \(\qq\).

\paragraph{The standard construction}
We recall how the HTC schemes of \cite{HTC2022,ABGRALL2023127629,BustoDumbser2024} handle the numerical dissipation, since the construction of Section~\ref{sec:viscosity} follows the same steps in different variables. The dissipative contribution is a scalar viscosity on the evolved variables,
\begin{equation}\label{eq:standard_viscosity}
	\mathcal G^{(i)}_{ij}=-\mathcal G^{(j)}_{ji}=\frac{\epsilon_{ij}}{|\bm\delta_{ij}|}\,\Delta\qq_{ij},
	\qquad \epsilon_{ij}\ge0 .
\end{equation}
It conserves every component of \(\qq\), but not, by itself, the energy.
Indeed, its work on the energies of the two cells is
\[
    \pp_i\cdot\mathcal G^{(i)}_{ij}
    +
    \pp_j\cdot\mathcal G^{(j)}_{ji}
    =
    -\frac{\epsilon_{ij}}{|\bm\delta_{ij}|}
    \,\Delta\pp_{ij}\cdot\Delta\qq_{ij},
\]
and the entropy productions are chosen to compensate this change.
A common numerical energy flux is defined using the
arithmetic average of the forces,
\[
    \mathcal G^{\rhoE}_{ij}
    =
    \frac{\epsilon_{ij}}{|\bm\delta_{ij}|}
    \,\frac{\pp_i+\pp_j}{2}\cdot\Delta\qq_{ij}.
\]
This choice can be interpreted as the trapezoidal approximation of
\(\int_0^1\pp(\qq(s))\cdot\Delta\qq_{ij}\,ds
=\rhoE(\qq_j)-\rhoE(\qq_i)\),
where \(\qq(s)=\qq_i+s\Delta\qq_{ij}\).
In this case, the two defects
\eqref{eq:energy_defects} are equal:
\begin{equation}\label{eq:standard_defects}
    \begin{aligned}
        \Gamma^{(i)}_{ij}
        =
        \Gamma^{(j)}_{ji}
        &=
        \frac{\epsilon_{ij}}{2|\bm\delta_{ij}|}
        \,\Delta\pp_{ij}\cdot\Delta\qq_{ij}
        \\
        &=
        \frac{\epsilon_{ij}}{2|\bm\delta_{ij}|}
        \int_0^1
        \Delta\qq_{ij}\cdot
        \nabla^2_{\qq}\rhoE\bigl(\qq(s)\bigr)
        \,\Delta\qq_{ij}\,ds ,
    \end{aligned}
\end{equation}
by \eqref{eq:roe_hessian}.
They are non-negative if the Hessian is positive
semi-definite along the segment, as it is when the energy
is convex in the evolved variables.
The Roe matrix is needed only to prove this sign;
the pairing
\(\Delta\pp_{ij}\cdot\Delta\qq_{ij}\)
can be evaluated directly in the computations
\cite{BustoDumbser2024}.
The trapezoidal rule is used here to choose the numerical
energy flux, not to approximate the discrete energy balance.
The latter is exact because the entropy productions compensate
the defects of that same flux.

\paragraph{Bregman divergences}
The same total defect in \eqref{eq:standard_defects} can be shared between the two cells in another way.
For a convex function \(f\), we introduce
\begin{equation}\label{eq:bregman}
	D_f(x'\,|\,x)=f(x')-f(x)-\nabla f(x)\cdot(x'-x),
\end{equation}
which is the non-negative Bregman divergence \cite{Bregman1967}. If \(f\) is twice continuously differentiable along
\(x(s)=x+s\Delta x\), with \(\Delta x=x'-x\),
Taylor's formula with integral remainder gives
\begin{equation}\label{eq:bregman_integral}
    \begin{aligned}
        D_f(x'\,|\,x)
        &=
        \int_0^1(1-s)\,
        \Delta x\cdot\nabla^2f\bigl(x(s)\bigr)
        \,\Delta x\,ds,
        \qquad
        D_f(x\,|\,x')
        &=
        \int_0^1s\,
        \Delta x\cdot\nabla^2f\bigl(x(s)\bigr)
        \,\Delta x\,ds.
    \end{aligned}
\end{equation}
The two weights add up to one, and hence
\[
    D_f(x'\,|\,x)+D_f(x\,|\,x')
    =
    \bigl(\nabla f(x')-\nabla f(x)\bigr)
    \cdot(x'-x).
\]
The standard defects \eqref{eq:standard_defects} are
therefore half of the symmetrized divergence of the energy,
\[
    \Gamma^{(i)}_{ij}
    =
    \Gamma^{(j)}_{ji}
    =
    \frac{\epsilon_{ij}}{2|\bm\delta_{ij}|}
    \left[
        D_{\rhoE}(\qq_j\,|\,\qq_i)
        +
        D_{\rhoE}(\qq_i\,|\,\qq_j)
    \right].
\]
Now replace the arithmetic average of the forces by their
exact average along the segment,
\(
    \overline\pp_{ij}
    =
    \int_0^1\pp(\qq(s))\,ds.
\)
The corresponding energy flux is
\[
    \mathcal G^{\rhoE}_{ij}
    =
    \frac{\epsilon_{ij}}{|\bm\delta_{ij}|}
    \,\overline\pp_{ij}\cdot\Delta\qq_{ij}
    =
    \frac{\epsilon_{ij}}{|\bm\delta_{ij}|}
    \left[\rhoE(\qq_j)-\rhoE(\qq_i)\right].
\]
Subtracting the work of the viscous contribution in cell
\(i\), we obtain
\[
    \begin{aligned}
        \Gamma^{(i)}_{ij}
        =
        \frac{\epsilon_{ij}}{|\bm\delta_{ij}|}
        \left[
            \rhoE(\qq_j)-\rhoE(\qq_i)
            -\pp_i\cdot(\qq_j-\qq_i)
        \right]
        =
        \frac{\epsilon_{ij}}{|\bm\delta_{ij}|}
        \,D_{\rhoE}(\qq_j\,|\,\qq_i),
        	\qquad
	    \Gamma^{(j)}_{ji}
    =
    \frac{\epsilon_{ij}}{|\bm\delta_{ij}|}
    \,D_{\rhoE}(\qq_i\,|\,\qq_j).
    \end{aligned}
\]
Thus the Bregman divergences are not additional dissipative
terms: they are the two energy defects associated with this
choice of numerical energy flux.
The viscous state contribution is unchanged, and so is the
sum of the defects; only their partition between the cells
changes. The configurational viscosity of Section~\ref{sec:viscosity} applies this energy-defect
construction in the phase variables.

\paragraph{Convexity in the evolved variables}
Convexity along the segment in \(\qq\) is a sufficient
condition for both partitions to have non-negative defects.
For the model considered here this hypothesis fails: by Proposition~\ref{prop:nonconvex_q} the Hessian \(\nabla^2_{\qq}\rhoE\) is indefinite wherever the phases slip, and the pairing \(\Delta\pp\cdot\Delta\qq\) can be negative (Section~\ref{subsec:implications_dissipation} and \ref{app:pairing}).

The reference viscosity used in Section~\ref{sec:results},
labelled Rusanov, is \eqref{eq:standard_viscosity} with
\(\epsilon_{ij}/|\bm\delta_{ij}|=\frac12s_{ac}\),
where \(s_{ac}\) is defined in
Section~\ref{subsec:two_scales}.
Where its production is negative, the solver sets it to zero.
This preserves the sign of the production, but no longer
compensates the corresponding energy defect.
Unless the dissipative contribution or the energy flux is
also modified, the compatibility identity
\eqref{eq:production_requirement} is then lost.

\paragraph{Further requirements of the multiphase problem}
Non-negative energy defects are necessary for the above
construction, but do not ensure the other properties needed
by a multiphase scheme. Consider a dissipation written in terms of the thermodynamic forces, \(\mathcal G=\mathbf K\,\Delta\pp\) with a symmetric positive semi-definite \(\mathbf K\). Its pairing \(\Delta\pp\cdot\mathbf K\Delta\pp\) is non-negative by construction. Two observations explain why this property alone is
not sufficient here.

First, such a dissipation cannot act on a contact with
uniform intrinsic phase states, a common velocity and
a common pressure.
All the forces
\eqref{eq:forces_phi}--\eqref{eq:forces_alpha}
are then continuous, \(\Delta\pp=\mathbf0\), whereas
the volume fractions may jump. The dissipation vanishes, leaving the volume fractions
to the central discretization alone, which does not
preserve their positivity. More generally, since \(\mathbf K\) is positive semi-definite, \(\Delta\pp\cdot\mathbf K\Delta\pp=0\) implies \(\mathbf K\Delta\pp=\mathbf 0\): a dissipation of this form cannot produce an exchange that is entropy-neutral, such as the transport of a phase at fixed intrinsic state.

Second, a non-negative pairing does not determine how the
phase masses are exchanged between cells, or how the
dissipation acts on the products
\(\crho_a=\alpha_a\rho_a\) and
\(\eta_a=\crho_as_a\) that define the state. For the multiphase model we therefore ask, in addition:
\begin{enumerate}
	\item[(R1)] \emph{Mass exchange.} The partial densities and the volume fractions remain positive, under an explicit time-step condition, and each phase receives mass only from the same phase in the neighbouring cells.
	\item[(R2)] \emph{Multiplicative structure.} A change of the amount of a phase must not alter its intrinsic state: the dissipation must respect the configurational direction of Section~\ref{subsec:configurational_direction}.
	\item[(R3)] \emph{Two scales.} The volume fractions, which carry the geometry of the phases, may receive a dissipation smaller than the acoustic one, without violating (R1)--(R2).
	\item[(R4)] \emph{Contact equilibria.} Phases with uniform intrinsic states, a common velocity and a common pressure remain in this equilibrium, and the interface between them is only transported.
\end{enumerate}
Section~\ref{sec:viscosity} constructs a dissipation that satisfies the energy requirement with non-negative defects, without any convexity assumption in \(\qq\), together with (R1)--(R4).

\section{A thermodynamically compatible configurational viscosity}\label{sec:viscosity}

This section constructs the numerical dissipation. Following the principle of Section~\ref{sec:geometry}, it is built phase by phase in the variables \(Z_a\) and transferred to \(\qq\) afterwards. Each step answers one of the requirements of Section~\ref{subsec:viscous_requirements}. Throughout the section we consider one face, oriented from cell \(i\) to cell \(j\); to lighten the notation we write \(L=i\), \(R=j\), \(\Delta(\cdot)=(\cdot)_R-(\cdot)_L\), and we drop the phase index when no confusion arises.

\subsection{An exact decomposition of the phase jump}\label{subsec:jump_decomposition}

The jump of the extensive variables of a phase across the face has two origins: the amount of phase changes, and its intrinsic state changes. The product relation \(X=\alpha Y\) separates them exactly
\[
    \delta X=Y\,\delta\alpha+\alpha\,\delta Y.
\]
We now seek the corresponding decomposition for the finite
jump across a cell face.
\begin{lemma}[Decomposition of the jump]
\label{lem:jump_decomposition}
For every phase and every pair of states with
\(\alpha_L,\alpha_R>0\),
\begin{equation}\label{eq:jump_decomposition}
    \Delta X
    =
    \Delta(\alpha Y)
    =
    Y^*\,\Delta\alpha+h\,\Delta Y,
    \qquad
    Y^*=
    \frac{\alpha_LY_L+\alpha_RY_R}{\alpha_L+\alpha_R},
    \qquad
    h=
    \frac{2\alpha_L\alpha_R}{\alpha_L+\alpha_R}.
\end{equation}
\end{lemma}

\begin{proof}
Subtracting \(Y^*\Delta\alpha\) from \(\Delta X\) gives
\[
    \begin{aligned}
        \Delta X-Y^*\Delta\alpha
        =
        \alpha_RY_R-\alpha_LY_L
        -
        \frac{\alpha_LY_L+\alpha_RY_R}
             {\alpha_L+\alpha_R}
        (\alpha_R-\alpha_L)
        =
        \frac{2\alpha_L\alpha_R}{\alpha_L+\alpha_R}
        (Y_R-Y_L).
    \end{aligned}
\]
\end{proof}
Both averages have a simple meaning. The state \(Y^*=(X_L+X_R)/(\alpha_L+\alpha_R)\) is the intrinsic state obtained by mixing the two portions of the phase, a volume-weighted average of the intrinsic states, and \(h\) is the harmonic mean of the two volume fractions. The arithmetic averages would also give an exact identity, \(\Delta X=\overline Y\Delta\alpha+\overline\alpha\Delta Y\); the volume-weighted one is the choice for which the dissipation of Section~\ref{subsec:exchange_form} takes an exchange form with non-negative weights, and for which the intrinsic jump is weighted by the trace side.
The full phase jump therefore splits into
\begin{equation}\label{eq:jump_split_Z}
    \Delta Z
    =
    \underbrace{\Delta\alpha\,(1,Y^*)}_{\text{configurational}}
    +
    \underbrace{h\,(0,\Delta Y)}_{\text{intrinsic}}.
\end{equation}
The first term is the configurational direction of
Section~\ref{subsec:configurational_direction}, evaluated
at \(Y^*\): a change of volume fraction is accompanied by
the corresponding change of the extensive variables.
The second term has no volume fraction component and
accounts for the difference between the intrinsic states.

Two limiting cases help interpret the decomposition.
If \(Y_L=Y_R\), the intrinsic contribution vanishes and
the entire jump is configurational.
If \(\alpha_L=\alpha_R=\alpha\), the configurational
contribution vanishes and \(h=\alpha\), so that
\(\Delta X=\alpha\Delta Y\).
When the volume fractions differ strongly, the harmonic
coefficient is controlled by the smaller one,
\[
    h\le 2\min(\alpha_L,\alpha_R).
\]
The intrinsic jump is therefore weighted by the amount
of phase available on the trace side, rather than by
the arithmetic mean of the two volume fractions.

\subsection{Two dissipative scales}\label{subsec:two_scales}

The decomposition \eqref{eq:jump_split_Z} allows us to assign
different dissipative coefficients to changes of the phase
amount and changes of its intrinsic state (requirements (R2)--(R3)). The volume fraction is not dissipated independently of the
extensive variables: its contribution is accompanied by
the corresponding density, momentum and entropy, in the
ratio \(Y_a^*\). For each phase we define the dissipative contribution as
\begin{equation}\label{eq:two_scales}
	\mathcal G_{Z,a}
	=
	\tfrac12\Bigl[\,s_f\,\Delta\alpha_a\,(1,\ Y_a^*)+s_{ac}\,h_a\,(0,\ \Delta Y_a)\Bigr],
\end{equation}
where \(s_f\) and \(s_{ac}\) have the dimensions of a speed
and satisfy \(0\le s_f\le s_{ac}\). The coefficient \(s_f\) controls the configurational
contribution, including the volume fraction exchange
\(\frac12\,s_f\Delta\alpha_a\).
The coefficient \(s_{ac}\) controls the intrinsic
contribution, which includes changes of density and
velocity associated with acoustic waves.
This separation does not require a characteristic
decomposition of the system.
If the coefficients coincide, the exact jump identity gives
\[
    \mathcal G_{Z,a}
    =
    \tfrac12\,s_{ac}\Delta Z_a.
\]
We recover the scalar Rusanov viscosity in the phase
variables, not in the evolved variables \(\qq\).
When \(s_f<s_{ac}\), the configurational contribution is
reduced while the intrinsic contribution retains the
acoustic coefficient.
In particular, for identical intrinsic states on the two
sides, \(\Delta Y_a=0\), the exchange depends only on
\(s_f\).
This is how the construction reduces interface diffusion
without discarding the product relation
\(X_a=\alpha_aY_a\).

\paragraph{Choice of the speeds}
Let \(\lambda_{ij}\) be an upper bound, over the two cells, of the characteristic speeds of the system in the direction \(\mathbf n_{ij}\) and of the normal phase velocities \(|\vvv_{a,s}\cdot\mathbf n_{ij}|\). Let \(\theta_{ij}\) be the thermodynamic correction
of \eqref{eq:abgrall_theta}.
We choose
\begin{equation}\label{eq:speeds}
	s_{ac}=\lambda_{ij}+|\theta_{ij}|,
	\qquad
	s_f=\max_{a,\ s\in\{L,R\}}|\vvv_{a,s}\cdot\mathbf n_{ij}|+|\theta_{ij}| .
\end{equation}
The acoustic speed is the coefficient of the Rusanov flux, augmented by the magnitude of the correction, which enters the reversible flux exactly as a signed viscosity (Section~\ref{subsec:abgrall}). The {phase speed} \(s_f\) is the speed at which material moves across the face, augmented in the same way; it is the value required by the positivity arguments of Section~\ref{subsec:positivity}. Since \(\lambda_{ij}\) bounds the phase velocities,
\(0\le s_f\le s_{ac}\).
The same \(s_f\) is used for every phase.
Consequently,
\[
    \sum_a(\mathcal G_{Z,a})_\alpha
    =
    \tfrac12\,s_f\sum_a\Delta\alpha_a
    =
    0.
\]
The exchange preserves the saturation constraint and is
an admissible variation for the chain rule of
Section~\ref{subsec:implications_dissipation}.
\subsection{Exchange form with non-negative weights}\label{subsec:exchange_form}

The definition \eqref{eq:two_scales}, based on the decomposition \eqref{eq:jump_split_Z} into configurational and intrinsic contributions, explains the two scales.
To see how the same operator moves phase mass between
cells, it is useful to write it directly in terms of
the two states.

\begin{proposition}[Exchange form]\label{prop:exchange_form}
Assume \(s_{ac}>0\), and let
\[
    r=\frac{s_f}{s_{ac}}\in[0,1],
    \qquad
    d_a=
    \frac{\alpha_{a,R}-\alpha_{a,L}}
         {\alpha_{a,R}+\alpha_{a,L}}.
\]
Then
\begin{equation}\label{eq:exchange_form}
    \mathcal G_{Z,a}
    =
    \tfrac12\,s_{ac}
    \bigl(b_{a,R}\,Z_{a,R}-b_{a,L}\,Z_{a,L}\bigr),
    \qquad
    b_{a,L}=1+(1-r)d_a,
    \qquad
    b_{a,R}=1-(1-r)d_a,
\end{equation}
with the weights satisfying
\[
    r\le b_{a,L},\,b_{a,R}\le2-r.
\]
\end{proposition}

\begin{proof}
By the jump decomposition,
\[
    \mathcal G_{Z,a}
    =
    \tfrac12\,s_{ac}
    \left[
        \Delta Z_a
        -(1-r)\Delta\alpha_a(1,Y_a^*)
    \right].
\]
The definition of \(Y_a^*\) gives
\[
    \Delta\alpha_a(1,Y_a^*)
    =
    d_a(Z_{a,L}+Z_{a,R}).
\]
Substitution and collection of the two endpoint states
yield \eqref{eq:exchange_form}.
The bounds follow from \(|d_a|\le1\) and \(0\le r\le1\).
\end{proof}

The exchange form says that each phase sends to the neighbouring cell a fraction of what it has, \(\frac12\,s_{ac}\,b_{a,L}\,Z_{a,L}\), and receives a fraction of what the neighbour has, \(\frac12\,s_{ac}\,b_{a,R}\,Z_{a,R}\). Mass, momentum, entropy and volume of a phase are exchanged with the same weights; no coefficient involves the ratio of the masses on the two sides, and phase \(a\) receives only phase \(a\) (requirement (R1)). The dimensionless weights \(b\) are bounded below by \(r\): the outflow of a cell is proportional to its own content. The following proposition makes the consequences explicit for the exchange alone.

\begin{proposition}[Exchange step]\label{prop:exchange_step}
Consider one forward Euler step of the isolated exchange
in \(Z\), with all coefficients evaluated at the
beginning of the step.
For a face \(ij\), denote by \(b^{(i)}_{a,ij}\) and
\(b^{(j)}_{a,ij}\) the weights multiplying the states
of cells \(i\) and \(j\), respectively, and set
\[
    \nu_{ij}
    =
    \Delta t\,
    \frac{|\partial\Omega_{ij}|}{2|\Omega_i|}
    s_{ac,ij}.
\]
Then
\[
    Z_{a,i}^{+}
    =
    c_{ii}\,Z_{a,i}
    +
    \sum_{j\in\mathcal N_i}c_{ij}\,Z_{a,j},
    \qquad
    c_{ii}
    =
    1-\sum_{j\in\mathcal N_i}
          \nu_{ij}\,b^{(i)}_{a,ij},
    \qquad
    c_{ij}
    =
    \nu_{ij}\,b^{(j)}_{a,ij},
\]
where the dependence of the coefficients \(c\) on the
phase is left implicit.
All coefficients are non-negative provided
\begin{equation}\label{eq:exchange_step_condition}
    \Delta t
    \sum_{j\in\mathcal N_i}
    \frac{|\partial\Omega_{ij}|}{2|\Omega_i|}
    s_{ac,ij}\,b^{(i)}_{a,ij}
    \le1.
\end{equation}
Under this condition, non-negative partial densities
and volume fractions remain non-negative, and the
saturation constraint is preserved.
Where the updated phase volume and mass are positive,
the new intrinsic density is a volume-weighted average
of the intrinsic densities in the stencil, while
the new velocity and specific entropy are mass-weighted
averages.
In particular, a uniform intrinsic phase state is
preserved.
\end{proposition}

\begin{proof}
The coefficients follow by collecting the incoming
and outgoing terms of \eqref{eq:exchange_form}.
The incoming coefficients are non-negative because
\(b\ge0\); condition \eqref{eq:exchange_step_condition}
makes the remaining coefficient \(c_{ii}\) non-negative.
The same coefficients act on all components of \(Z_a\).
Thus, writing \(k\) for the cell itself and its neighbours,
\[
    \alpha_{a,i}^{+}
    =
    \sum_k c_{ik}\,\alpha_{a,k},
    \qquad
    X_{a,i}^{+}
    =
    \sum_k c_{ik}\,X_{a,k}.
\]
The positivity statements follow.
The saturation constraint is preserved because the
volume fraction exchanges sum to zero over the phases
on every face.

Since \(X_{a,k}=\alpha_{a,k}Y_{a,k}\), division by
\(\alpha_{a,i}^{+}>0\) gives
\[
    Y_{a,i}^{+}
    =
    \sum_k
    \frac{c_{ik}\,\alpha_{a,k}}
         {\sum_l c_{il}\,\alpha_{a,l}}
    \,Y_{a,k}.
\]
These normalized weights are non-negative and sum to one.
The density component therefore gives the stated
volume-weighted average.
For the phase velocity and specific entropy, division
by the updated partial density instead gives
\[
    \vvv_{a,i}^{+}
    =
    \sum_k
    \frac{c_{ik}\,\crho_{a,k}}
         {\sum_l c_{il}\,\crho_{a,l}}
    \,\vvv_{a,k},
    \qquad
    s_{a,i}^{+}
    =
    \sum_k
    \frac{c_{ik}\,\crho_{a,k}}
         {\sum_l c_{il}\,\crho_{a,l}}
    \,s_{a,k}.
\]
These are the mass-weighted averages.
\end{proof}

The coefficients \(c_{ik}\) need not sum to one:
the amount of a phase in the cell changes during the
exchange.
It is after division by the updated phase volume or
mass that the weights sum to one.
This is the role of the multiplicative structure:
incoming mass carries the momentum and entropy of
the same phase, rather than changing those quantities
independently.

The proposition concerns the isolated dissipative exchange advanced in \(Z\). It is not a maximum-principle statement for the complete
time step in \(\qq\), since the map from phase momenta
to relative velocities is nonlinear.
Nor does it bound the specific entropy after the
production required by energy compatibility is added.
The mass and volume fraction rows, however, are linear
in \(Z\); their positivity in the presence of the
reversible contributions is established separately
in Theorem~\ref{thm:positivity}.

\subsection{Side-wise transfer to the model variables}\label{subsec:transfer}

The exchange has been constructed in the phase variables
\(Z\), whereas the numerical scheme evolves \(\qq\).
Since \(\qq\) is a function of \(Z\) (Section~\ref{subsec:phase_variables}), we transfer a variation \(\delta Z\)
using the same chain rule that relates variations of the two representations.
At the state of cell \(s\),
\[
    \delta\qq=\mathbf J_s\,\delta Z,
    \qquad
    \mathbf J_s
    =
    \left.\frac{\partial\qq}{\partial Z}\right|_s .
\]
On the face oriented from \(i=L\) to \(j=R\), the left
cell receives \(\mathcal G_Z\) and the right cell receives
\(-\mathcal G_Z\).
Their contributions in the evolved variables are therefore
\begin{equation}\label{eq:transfer}
    \mathcal G^{(L)}:=\mathcal G^{(i)}_{ij}
    =
    \mathbf J_L\,\mathcal G_Z,
    \qquad
    \mathcal G^{(R)}:=\mathcal G^{(j)}_{ji}
    =
    -\mathbf J_R\,\mathcal G_Z,
    \qquad
    \mathcal G_Z
    =
    (\mathcal G_{Z,1},\ldots,\mathcal G_{Z,\Nph}).
\end{equation}
The exchange in \(Z\) is computed once.
What differs between the two cells is the Jacobian through
which that exchange acts.

\paragraph{Linear rows}
The rows of \(\UU\), \(\crho_a\), \(\eta_a\), \(\rho\)
and \(\alpha_a\) in \eqref{eq:q_of_Z} are linear in \(Z\).
Their Jacobians are independent of the state, so that
\eqref{eq:transfer} gives equal and opposite contributions
to the two cells.

For example, the contributions to the left cell are
\[
    \begin{gathered}
        (\mathcal G^{(L)})_{\UU}=\sum_a\mathcal G_{\bm M,a},
        \qquad
        (\mathcal G^{(L)})_{\rho}=\sum_a\mathcal G_{\crho,a},
        \qquad
        (\mathcal G^{(L)})_{\alpha_a}=\tfrac12\,s_f\Delta\alpha_a,
        \\
        (\mathcal G^{(L)})_{\crho_a}=\mathcal G_{\crho,a},
        \qquad
        (\mathcal G^{(L)})_{\eta_a}=\mathcal G_{\eta,a}.
    \end{gathered}
\]
The right cell receives their opposites.
The dissipative contributions on these rows are thus
single-valued fluxes, and the conservation of mixture momentum, phase masses and phase entropies is retained.

\paragraph{Relative velocity rows}
The relative velocities are differences of ratios.
For a phase velocity \(\vvv_a=\bm M_a/\crho_a\), the chain rule gives
\[
    \delta\vvv_a
    =
    \frac{\delta\bm M_a-\vvv_a\,\delta\crho_a}
         {\crho_a}.
\]
The velocity changes because momentum enters or leaves
the cell, but also because the mass with which that
momentum is divided changes.
Denote by
\(\mathcal G^{(L)}_{\vvv_a}\) and
\(\mathcal G^{(R)}_{\vvv_a}\)
the phase-velocity contributions induced by the exchange
in the two cells.
Using \eqref{eq:exchange_form}, the numerator on the
left side is
\[
    \mathcal G_{\bm M,a}
    -
    \vvv_{a,L}\,\mathcal G_{\crho,a}
    =
    \tfrac12\,s_{ac}\,b_{a,R}\,\crho_{a,R}\,
    (\vvv_{a,R}-\vvv_{a,L}).
\]
The outgoing mass carries the velocity of the left cell,
so its contribution cancels in this difference.
The same calculation on the right, including the opposite
sign of the exchange, yields
\begin{equation}\label{eq:transfer_w}
    \begin{aligned}
        \mathcal G^{(L)}_{\vvv_a}
        &=
        \tfrac12\,s_{ac}\,b_{a,R}\,
        \frac{\crho_{a,R}}{\crho_{a,L}}\,
        (\vvv_{a,R}-\vvv_{a,L}),
        \qquad
        \mathcal G^{(R)}_{\vvv_a}
        &=
        \tfrac12\,s_{ac}\,b_{a,L}\,
        \frac{\crho_{a,L}}{\crho_{a,R}}\,
        (\vvv_{a,L}-\vvv_{a,R}).
    \end{aligned}
\end{equation}
Each contribution in \eqref{eq:transfer_w} drives the phase velocity towards
the velocity of the same phase in the neighbouring cell, at a rate proportional to the mass that enters relative to the mass present.
The relative velocity rows receive the difference of
the corresponding phase-velocity contributions:
\[
    \begin{aligned}
        (\mathcal G^{(L)})_{\ww_a}
        &=
        \mathcal G^{(L)}_{\vvv_a}
        -
        \mathcal G^{(L)}_{\vvv_\Nph},
        \qquad
        (\mathcal G^{(R)})_{\ww_a}
        &=
        \mathcal G^{(R)}_{\vvv_a}
        -
        \mathcal G^{(R)}_{\vvv_\Nph}.
    \end{aligned}
\]
These contributions need not be equal and opposite.
Conservation is unaffected, because the relative velocities are not conserved quantities; this is the freedom of Section~\ref{subsec:notation}.

\paragraph{Energy pairing}
The same transfer also preserves the work of the
thermodynamic forces, separately in each cell.
The chain rule \eqref{eq:chain_rule_admissible} applies because the exchange is admissible (Section~\ref{subsec:two_scales}). Consequently,
\begin{equation}\label{eq:transfer_duality}
    \pp_L\cdot\mathcal G^{(L)}
    =
    \sum_a\pp^Z_{a,L}\cdot\mathcal G_{Z,a},
    \qquad
    \pp_R\cdot\mathcal G^{(R)}
    =
    -\sum_a\pp^Z_{a,R}\cdot\mathcal G_{Z,a}.
\end{equation}
The total change of energy, \(\pp_L\cdot\mathcal G^{(L)}+\pp_R\cdot\mathcal G^{(R)}=-\sum_a\Delta\pp^Z_a\cdot\mathcal G_{Z,a}\), is the one computed in the phase variables, without any convexity assumption in \(\qq\). The energy flux and the entropy productions that complete the balance are constructed in Section~\ref{subsec:energy_production}.

\begin{figure}[t]
\centering
\begin{tikzpicture}[>=stealth,font=\small]
	\draw[thick] (0,0) rectangle (4,2.4);
	\draw[thick] (4,0) rectangle (8,2.4);
	\node[anchor=south west] at (0.1,0.05) {cell \(i=L\)};
	\node[anchor=south east] at (7.9,0.05) {cell \(j=R\)};
	\node[draw,rounded corners,fill=gray!10,align=center] (GZ) at (4,3.6) {exchange in \(Z\)\\\(\mathcal G_Z=\tfrac12\,s_{ac}(b_RZ_R-b_LZ_L)\)};
	\node[align=center] (GL) at (2,1.35) {\(\mathcal G^{(L)}=\mathbf J_L\mathcal G_Z\)\\[1mm] defect \(\Gamma^{(L)}\ge0\)};
	\node[align=center] (GR) at (6,1.35) {\(\mathcal G^{(R)}=-\mathbf J_R\mathcal G_Z\)\\[1mm] defect \(\Gamma^{(R)}\ge0\)};
	\draw[->,thick] (GZ.south) ++(-0.4,0) -- (2,2.4);
	\draw[->,thick] (GZ.south) ++(0.4,0) -- (6,2.4);
	\node[font=\footnotesize,anchor=east] at (2.5,2.9) {\(\mathbf J_L\)};
	\node[font=\footnotesize,anchor=west] at (5.5,2.9) {\(\mathbf J_R\)};
	\draw[<->,thick,dashed] (3.2,-0.35) -- (4.8,-0.35);
	\node at (4,-0.7) {single-valued energy flux \(\mathcal G^{\rhoE}_{ij}\)};
\end{tikzpicture}
\caption{The configurational viscosity on one face. The exchange is built once, in the phase variables; each cell receives it through its own chain rule. The conservative rows coincide on the two sides, the rows of the relative velocities do not. The energy flux is single-valued, and each cell converts its own non-negative energy defect into the entropy of its phases.}
\label{fig:face}
\end{figure}
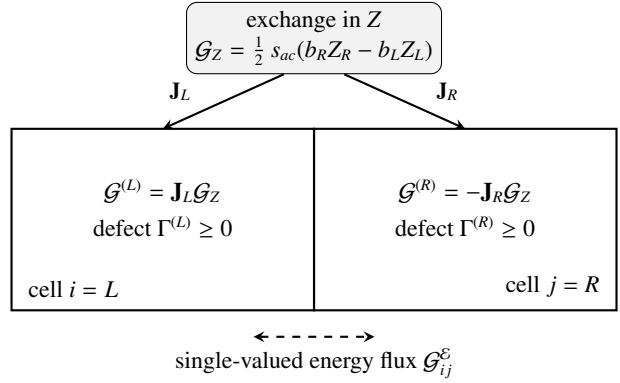

This is why the transfer is side-wise. A single-valued flux in \(\qq\) would require \(\mathbf J_L\,\mathcal G_Z=\mathbf J_R\,\mathcal G_Z\), which holds on the linear rows but not, by \eqref{eq:transfer_w}, on the relative velocity rows. Its work on at least one cell would then differ in general from that of the phase forces, and convexity in \(Z\) would no longer suffice to establish the sign of the energy defects. With the side-wise transfer, each cell receives the contribution prescribed by its own chain rule, which fixes both the state variation and its work, while the linear rows keep a single-valued flux. Figure~\ref{fig:face} summarizes the construction.

\subsection{Energy flux and entropy production}\label{subsec:energy_production}

The exchange and its transfer to the evolved variables have now been specified.
To complete the dissipative part, we close the energy balance of the exchange following the steps of the standard construction of Section~\ref{subsec:viscous_requirements}.
Namely, a common energy flux is chosen, the defects of the two cells are identified, and they are converted into entropy. Two things change: the exchange is built in the phase variables, where each phase energy is convex, and the energy flux uses the exact energies of the exchanged states, which for \(s_f=s_{ac}\) corresponds to the exact average of the forces along the segment rather than the trapezoidal one. The defects are then the Taylor remainders \eqref{eq:bregman_integral} of the phase energies.

\paragraph{Bregman divergences of the phase energies}
For the single-fluid energy \eqref{eq:intrinsic_energy}, with \(\bm\theta=\nabla_YE\) from \eqref{eq:dual_Z}, the divergence \eqref{eq:bregman} reads
\begin{equation}\label{eq:bregman_fluid}
	D_{E}(Y'\,|\,Y)
	=
	\tfrac12\rho'|\vvv'-\vvv|^2
	+
	\Bigl[\varepsilon(\rho',\sigma')-\varepsilon(\rho,\sigma)-g\,(\rho'-\rho)-T\,(\sigma'-\sigma)\Bigr],
\end{equation}
where \(g=\partial_\rho\varepsilon\) and
\(T=\partial_\sigma\varepsilon\) are evaluated at the
unprimed state.
The first term is the kinetic energy of the primed
state measured relative to the velocity of the unprimed
state.
The remaining terms form the Bregman divergence of the
internal energy density.
Both parts are non-negative under the convexity
assumptions of Section~\ref{subsec:perspective}.

For the perspective energy, Euler's identity gives
\[
    \rhoE_a(Z)
    =
    \pp^Z_a(Z)\cdot Z.
\]
The terms involving the energy and its tangent at the
reference state therefore cancel:
\[
    \begin{aligned}
        D_{\rhoE_a}(Z'\,|\,Z)
        =
        \rhoE_a(Z')-\rhoE_a(Z)
        -\pp^Z_a(Z)\cdot(Z'-Z)
        =
        \rhoE_a(Z')-\pp^Z_a(Z)\cdot Z'.
    \end{aligned}
\]
Using \(\pp^Z_a=(-p_a,\bm\theta_a)\),
\(Z'=\alpha'(1,Y')\), and
\(E_a(Y)=\bm\theta_a\cdot Y-p_a\), we obtain
\begin{equation}\label{eq:bregman_perspective}
    D_{\rhoE_a}(Z'\,|\,Z)
    =
    \rhoE_a(Z')-\pp^Z_a(Z)\cdot Z'
    =
    \alpha'\,D_{E_a}(Y'\,|\,Y).
\end{equation}
The divergence depends only on the two intrinsic states, weighted by the volume fraction of the first argument.
It vanishes when \(Y'=Y\), even if
\(\alpha'\ne\alpha\).
A change of phase amount at fixed intrinsic state
changes the energy linearly and leaves no Taylor
remainder. Unlike the Roe matrix of the Hessian, \eqref{eq:bregman_fluid} has a closed form for the equations of state used here, which is evaluated without quadrature and without cancellations (\ref{app:bregman}).

\paragraph{Energy flux and defects}
In \eqref{eq:exchange_form}, each cell sends a contribution
proportional to its own phase state.
We assign to this contribution the energy of the same
state, with the same coefficient:
\begin{equation}\label{eq:energy_flux}
	\mathcal G^{\rhoE}_{ij}=\tfrac12\,s_{ac}\sum_a\bigl(b_{a,R}\,\rhoE_a(Z_{a,R})-b_{a,L}\,\rhoE_a(Z_{a,L})\bigr),
	\qquad
	\mathcal G^{\rhoE}_{ji}=-\mathcal G^{\rhoE}_{ij}.
\end{equation}
Once this flux is chosen, the defects
\eqref{eq:energy_defects} are fixed.

\begin{theorem}[Energy defects]\label{thm:defects}
Assume positive phase volume fractions and densities,
convex single-fluid energies, and
\(0\le s_f\le s_{ac}\).
With the transfer \eqref{eq:transfer} and the energy
flux \eqref{eq:energy_flux}, the defects of the two
cells,
\(\Gamma^{(L)}:=\Gamma^{(i)}_{ij}\) and
\(\Gamma^{(R)}:=\Gamma^{(j)}_{ji}\), are
\begin{equation}\label{eq:defects}
    \begin{aligned}
        \Gamma^{(L)}
        =
        \tfrac12\,s_{ac}
        \sum_a
        b_{a,R}\,\alpha_{a,R}\,
        D_{E_a}(Y_{a,R}\,|\,Y_{a,L}),
        \qquad
        \Gamma^{(R)}
        =
        \tfrac12\,s_{ac}
        \sum_a
        b_{a,L}\,\alpha_{a,L}\,
        D_{E_a}(Y_{a,L}\,|\,Y_{a,R}).
    \end{aligned}
\end{equation}
Both defects are non-negative, without any convexity
assumption on the energy in \(\qq\).
\end{theorem}

\begin{proof}
The side-wise identity \eqref{eq:transfer_duality}
allows the work of the state contribution to be
computed in the phase variables:
\[
    \Gamma^{(L)}
    =
    \mathcal G^{\rhoE}_{ij}
    -
    \sum_a
    \pp^Z_{a,L}\cdot\mathcal G_{Z,a}.
\]
Substituting the state and energy exchanges gives
\[
    \begin{aligned}
        \Gamma^{(L)}
        =
        \tfrac12\,s_{ac}\sum_a
        \Bigl\{
        &b_{a,R}
        \left[
            \rhoE_a(Z_{a,R})
            -
            \pp^Z_{a,L}\cdot Z_{a,R}
        \right]
        -
        b_{a,L}
        \left[
            \rhoE_a(Z_{a,L})
            -
            \pp^Z_{a,L}\cdot Z_{a,L}
        \right]
        \Bigr\}.
    \end{aligned}
\]
The second bracket is zero by Euler's identity.
It corresponds to the outgoing contribution, whose
energy is exactly accounted for by the forces of
the cell from which it leaves.
The first bracket is
\(D_{\rhoE_a}(Z_{a,R}\,|\,Z_{a,L})\):
the incoming state is compared with the tangent
energy at the receiving state.
Equation \eqref{eq:bregman_perspective} gives the
first formula in \eqref{eq:defects}.
For the right cell,
\[
    \Gamma^{(R)}
    =
    -\mathcal G^{\rhoE}_{ij}
    +
    \sum_a
    \pp^Z_{a,R}\cdot\mathcal G_{Z,a},
\]
and the same argument with the states reversed
gives the second formula.
Each term is non-negative because the exchange
weights and volume fractions are non-negative
and the single-fluid energies are convex.
\end{proof}

The two defects are generally different, but their
sum is the energy removed by the viscous state
contributions:
\[
    \Gamma^{(L)}+\Gamma^{(R)}
    =
    \sum_a
    \Delta\pp^Z_a\cdot\mathcal G_{Z,a}.
\]
Their separate signs are established by
\eqref{eq:defects}, not merely by the sign of this sum.
For \(s_f=s_{ac}\), all weights are one.
The state exchange is
\(\frac12\,s_{ac}\Delta Z_a\), and the total energy
flux is
\(\frac12\,s_{ac}\sum_a\Delta\rhoE_a\).
The two defects then reduce to
\[
    \Gamma^{(L)}
    =
    \tfrac12\,s_{ac}
    \sum_a
    D_{\rhoE_a}(Z_{a,R}\,|\,Z_{a,L}),
    \qquad
    \Gamma^{(R)}
    =
    \tfrac12\,s_{ac}
    \sum_a
    D_{\rhoE_a}(Z_{a,L}\,|\,Z_{a,R}).
\]
This is the exact-average partition of
Section~\ref{subsec:viscous_requirements}, written
in the phase variables.
For \(s_f<s_{ac}\), the exchange is weighted rather
than proportional to \(\Delta Z_a\).
Homogeneity allows the same construction to hold:
the outgoing contribution still has zero defect,
and the incoming contribution retains its
non-negative Bregman remainder.

\paragraph{Entropy production}
Each cell defect is already expressed as a sum
of non-negative phase contributions. We assign each contribution to the entropy of the phase it comes from:
\begin{equation}\label{eq:production}
    \begin{aligned}
        \Pi^{(L)}_a
        :=
        \Pi^{(i)}_{a,ij}
        &=
        \frac{s_{ac}}{2T_{a,L}}\,
        b_{a,R}\,\alpha_{a,R}\,
        D_{E_a}(Y_{a,R}\,|\,Y_{a,L}),
        \\
        \Pi^{(R)}_a
        :=
        \Pi^{(j)}_{a,ji}
        &=
        \frac{s_{ac}}{2T_{a,R}}\,
        b_{a,L}\,\alpha_{a,L}\,
        D_{E_a}(Y_{a,L}\,|\,Y_{a,R}).
    \end{aligned}
\end{equation}
The temperatures are those of the receiving cells.
They are positive, so both productions are
non-negative, and
\[
    \sum_aT_{a,L}\Pi^{(L)}_a=\Gamma^{(L)},
    \qquad
    \sum_aT_{a,R}\Pi^{(R)}_a=\Gamma^{(R)}.
\]
Thus \eqref{eq:production_requirement} holds on
each side of the face. Three features distinguish this production from that of the standard construction. It is non-negative without any clipping. It is computed phase by phase, so that no rule is needed to share a mixture defect among the phases. And it vanishes whenever the two intrinsic states of a phase coincide, whatever the volume fractions: a purely configurational exchange produces no entropy. This is consistent with the observation of Section~\ref{subsec:viscous_requirements} that a contact is entropy neutral, and it is precisely what a symmetric dissipation in the forces cannot reproduce.

\section{Properties of the scheme}\label{sec:properties}\label{subsec:properties}

We now collect the properties of the scheme obtained
by combining the reversible discretization of
Section~\ref{sec:framework} with the configurational
viscosity of Section~\ref{sec:viscosity}.
We first establish the energy and entropy balances,
then the positivity of phase masses and volume fractions
under suitable time-step conditions, and the preservation
of contact equilibria.
Finally, we show that the side-wise chain rule also
preserves the discrete relation between the volume fraction
formulation and the original density-weighted formulation.
The energy balance and the equivalence of the two
formulations are semi-discrete statements.

\subsection{Energy conservation and entropy production}\label{subsec:energy_entropy}

The construction of Section~\ref{subsec:energy_production}
provides the non-negative energy defects required by
the general compatibility result of
Section~\ref{subsec:viscous_requirements}.
Combining it with the reversible part gives the
following property of the complete semi-discrete scheme.

\begin{theorem}[Energy conservation and entropy production]\label{thm:energy_entropy}
Assume the hypotheses of Theorem~\ref{thm:defects},
the face-wise compatibility condition
\eqref{eq:face_energy_condition} for the reversible
fluctuations, and the relaxation sources of
Section~\ref{subsec:sources}.
Then the scheme satisfies the local energy balance
\eqref{eq:discrete_energy}, with the single-valued
numerical flux
\(\mathcal F^{\rhoE}_{ij}-\mathcal G^{\rhoE}_{ij}\).
The entropy of each phase, and hence the mixture entropy, satisfies a cell balance with a conservative numerical flux and a non-negative production.
\end{theorem}

\begin{proof}
The viscous contributions are conservative on the
rows of \(\UU\), \(\crho_a\), \(\eta_a\) and \(\rho\),
by Section~\ref{subsec:transfer}.
Their energy flux is single-valued.
Theorem~\ref{thm:defects} gives non-negative defects,
and the productions \eqref{eq:production} satisfy
\[
    \sum_a T_{a,L}\Pi_a^{(L)}=\Gamma^{(L)},
    \qquad
    \sum_a T_{a,R}\Pi_a^{(R)}=\Gamma^{(R)}.
\]
The hypotheses of Theorem~\ref{thm:discrete_energy}
are therefore satisfied.
\end{proof}

\subsection{Positivity of phase masses and volume fractions}\label{subsec:positivity}

The exchange step of
Proposition~\ref{prop:exchange_step} has non-negative weights.
We now show that the choice of \(s_f\) also controls
the central transport and the signed thermodynamic
correction, so that their combination preserves non-negative phase masses and the bounds of the volume fractions.
For simplicity, consider a uniform one-dimensional
grid with cell size \(\Delta x\).
Each face \(i+\frac12\) is oriented from \(i\) to
\(i+1\), and \(L,R\) denote its left and right states.
On Cartesian grids the argument applies face by face, with \(d_{a,i}\) and \(d_{\alpha,i}\) collecting the outgoing rates of all faces of the cell.

The volume fraction equation uses the mixture
velocity of the midpoint state,
\[
    \widetilde V_{i+\frac12}
    =
    \frac{U_i+U_{i+1}}{\rho_i+\rho_{i+1}}
    =
    \frac{
        \sum_a
        \left(
            \crho_{a,i}v_{a,i}
            +
            \crho_{a,i+1}v_{a,i+1}
        \right)
    }{\rho_i+\rho_{i+1}}.
\]
Since the partial densities are positive, this is
a mass-weighted average of the endpoint phase
velocities. In particular,
\[
    |\widetilde V_{i+\frac12}|
    \le
    \max_{\substack{a\\s\in\{i,i+1\}}}|v_{a,s}|.
\]

\begin{theorem}[Positivity of phase masses and volume fractions]\label{thm:positivity}
Consider the scheme without relaxation sources,
evaluated at admissible states, with the speeds
\eqref{eq:speeds}.
For each phase, the partial densities satisfy
\begin{equation}\label{eq:positivity_mass}
	\frac{d\crho_{a,i}}{dt}
	=
	\frac{1}{2\Delta x}\Bigl[
		\bigl(s_{ac}\,b_{a,R}-v_{a,i+1}+\theta\bigr)_{i+\frac12}\crho_{a,i+1}
		+\bigl(s_{ac}\,b_{a,L}+v_{a,i-1}+\theta\bigr)_{i-\frac12}\crho_{a,i-1}
	\Bigr]
	-d_{a,i}\,\crho_{a,i},
\end{equation}
where
\[
    d_{a,i}
    =
    \frac{1}{2\Delta x}
    \left[
        (s_{ac}\,b_{a,L}+\theta)_{i+\frac12}
        +
        (s_{ac}\,b_{a,R}+\theta)_{i-\frac12}
    \right]
    \ge0.
\]
The volume fractions satisfy
\begin{equation}\label{eq:positivity_alpha}
    \begin{aligned}
        \frac{d\alpha_{a,i}}{dt}
        =
        \frac{1}{2\Delta x}\Bigl[
        &(s_f-\widetilde V+\theta)_{i+\frac12}
        (\alpha_{a,i+1}-\alpha_{a,i})
        \\
        &+
        (s_f+\widetilde V+\theta)_{i-\frac12}
        (\alpha_{a,i-1}-\alpha_{a,i})
        \Bigr].
    \end{aligned}
\end{equation}
All coefficients multiplying neighbouring values
are non-negative.
Consequently, a forward Euler step preserves \(\crho_{a,i}\ge0\) if \(\Delta t\,d_{a,i}\le1\). It satisfies a local maximum principle for the
volume fractions if
\[
    \Delta t\,d_{\alpha,i}\le1,
    \qquad
    d_{\alpha,i}
    :=
    \frac{1}{2\Delta x}
    \left[
        (s_f-\widetilde V+\theta)_{i+\frac12}
        +
        (s_f+\widetilde V+\theta)_{i-\frac12}
    \right].
\]
Under this condition, \(\alpha_{a,i}^{+}\) lies
between the minimum and maximum of
\(\alpha_{a,i-1}\), \(\alpha_{a,i}\) and
\(\alpha_{a,i+1}\).
Thus volume fractions initially in \([0,1]\)
remain in that interval.
The same properties hold for SSP Runge--Kutta
methods, provided every forward Euler substep
satisfies these restrictions, with its effective
step size and the coefficients evaluated at the
corresponding stage.
\end{theorem}

\begin{proof}
The phase mass equation has physical flux
\(\crho_av_a\) and no non-conservative product.
Its total numerical flux on one face is
\[
    \widehat F_a
    =
    \tfrac12
    \left(
        \crho_{a,L}v_{a,L}
        +
        \crho_{a,R}v_{a,R}
    \right)
    -
    \tfrac12\theta
    (\crho_{a,R}-\crho_{a,L})
    -
    \tfrac12s_{ac}
    \left(
        b_{a,R}\crho_{a,R}
        -
        b_{a,L}\crho_{a,L}
    \right).
\]
Taking the difference of the fluxes at
\(i+\frac12\) and \(i-\frac12\) gives
\eqref{eq:positivity_mass}.
The central-flux terms containing the cell's own
quantity \(\crho_{a,i}v_{a,i}\) cancel between
the two faces.

The lower bound on the exchange weights gives
\(s_{ac}\,b\ge s_f\).
By the choice of the phase speed,
\[
    s_{ac}\,b+\theta
    \ge
    s_f+\theta
    =
    \max_{a,s}|v_{a,s}|+|\theta|+\theta
    \ge
    \max_{a,s}|v_{a,s}|.
\]
The incoming coefficients in
\eqref{eq:positivity_mass} are therefore
non-negative, as is the outgoing rate \(d_{a,i}\).
In a forward Euler step, the coefficient of
\(\crho_{a,i}\) is \(1-\Delta t\,d_{a,i}\);
the stated restriction makes it non-negative.

For the volume fractions, the midpoint
non-conservative fluctuation, the correction
and the viscous exchange give
\eqref{eq:positivity_alpha}.
Its two coefficients are non-negative because
\[
    s_f+\theta
    \ge
    \max_{a,s}|v_{a,s}|
    \ge
    |\widetilde V|.
\]
The Euler update has coefficient
\(1-\Delta t\,d_{\alpha,i}\) on the cell's
own value, and the other two coefficients sum
to \(\Delta t\,d_{\alpha,i}\).
It is therefore a convex combination under
the stated restriction.
Finally, phase masses and volume fractions are
linear functions of the evolved variables.
Convex combinations of states satisfying these
bounds satisfy them again, which gives the
SSP extension whenever the stage evaluations
are well-defined.
\end{proof}

The phase speed \eqref{eq:speeds} is exactly what these two arguments need: the material speed for the transport, and \(|\theta|\) for the correction, which by \eqref{eq:abgrall_flux} acts on the partial densities and on the volume fractions as a signed viscosity.
\subsection{Preservation of contact equilibria}\label{subsec:contacts}
We now consider a contact at mechanical equilibrium.
The volume fractions may vary from cell to cell, while
the intrinsic state of each phase is uniform.
The scheme should then transport the amount
of each phase without changing its density, specific
entropy, velocity or pressure.

\begin{theorem}[Contact equilibria]\label{thm:contact}
Suppose that, for every phase, the intrinsic density
and specific entropy are constant across the grid,
and that all phases share a common velocity \(\VV\)
and a common pressure.
The volume fractions may vary arbitrarily, subject to
\(\alpha_{a,i}>0\) and
\(\sum_a\alpha_{a,i}=1\).
Then, without relaxation sources, the semi-discrete
scheme preserves these intrinsic states and the common
velocity and pressure, and produces no entropy.

On a one-dimensional uniform grid, the volume fractions
evolve according to the first-order upwind scheme
\begin{equation}\label{eq:contact_upwind}
    \frac{d\alpha_{a,i}}{dt}
    =
    \frac{1}{2\Delta x}
    \Bigl[
        \bigl(|V|-V\bigr)
        (\alpha_{a,i+1}-\alpha_{a,i})
        +
        \bigl(|V|+V\bigr)
        (\alpha_{a,i-1}-\alpha_{a,i})
    \Bigr],
\end{equation}
where \(V\) is the common velocity in the grid direction.
\end{theorem}

\begin{proof}
The intrinsic states \(Y_a\) are the same in every cell,
although they may differ between phases.
Consequently,
\[
    \Delta Z_a=\Delta\alpha_a(1,Y_a).
\]
Since the relative velocities vanish, the corresponding
states \(\qq\) form an affine family parametrized by the
volume fractions.
The thermodynamic forces are constant on this family, since they depend only on the intrinsic states, the common velocity and the relative velocities (Section~\ref{subsec:balances});
denoting their common value by \(\pp_0\), we have
\[
    \Delta\pp=\mathbf 0,
    \qquad
    \Delta\rhoE=\pp_0\cdot\Delta\qq,
    \qquad
    \delta_{ij}=0.
\]

We first check that the reversible residual also vanishes.
On a face with normal \(\mathbf n\), let
\(V=\VV\cdot\mathbf n\).
The mass and entropy fluxes are advected with this
common velocity, the pressure contribution has no jump,
and the relative-velocity contributions vanish.
The midpoint rule for the volume fraction equation
therefore gives
\[
    \tfrac12\Delta\f_{\mathbf n,ij}
    +\mathcal B_{ij}
    =
    \tfrac12 V\,\Delta\qq.
\]
Moreover, the physical energy flux reduces to
\(V(\rhoE+p)\), with constant pressure \(p\), so that
\[
    \mathcal R_{ij}
    =
    V\,\Delta\rhoE
    -
    (\pp_i+\pp_j)\cdot
    \bigl(\tfrac12 V\,\Delta\qq\bigr)
    =0.
\]
Thus no thermodynamic correction is needed, and
\(\theta_{ij}=0\).

For the viscous part, \(\Delta Y_a=\mathbf 0\), and
\eqref{eq:two_scales} reduces to
\[
    \mathcal G_{Z,a}
    =
    \tfrac12 s_f\,\Delta\alpha_a(1,Y_a).
\]
This exchange changes the volume fraction and the
extensive phase variables in the same ratio \(Y_a\).
Because \(\Delta Z_a\) is tangent to the contact family, on which \(\qq\) is an affine function of \(Z\), both Jacobians map it to \(\Delta\qq\):
\[
    \mathcal G^{(L)}
    =\tfrac12 s_f\,\Delta\qq,
    \qquad
    \mathcal G^{(R)}
    =-\tfrac12 s_f\,\Delta\qq.
\]
In particular, it does not generate relative velocity.
The intrinsic Bregman divergences also vanish,
\(D_{E_a}(Y_{a,R}\mid Y_{a,L})=0\), so that all
entropy productions in \eqref{eq:production} are zero.

Both the reversible and the viscous contributions
therefore remain tangent to the same contact family:
\[
    \frac{d\crho_{a,i}}{dt}
    =
    \rho_a\frac{d\alpha_{a,i}}{dt},
    \qquad
    \frac{d\eta_{a,i}}{dt}
    =
    \rho_a s_a\frac{d\alpha_{a,i}}{dt},
    \qquad
    \frac{d\UU_i}{dt}
    =
    \VV\frac{d\rho_i}{dt},
    \qquad
    \frac{d\ww_{a,i}}{dt}=\mathbf 0.
\]
Hence the intrinsic densities and specific entropies,
the common velocity and, through the equations of state,
the pressures remain unchanged.
Finally, \(\theta=0\), \(\widetilde V=V\) and
\(s_f=|V|\); substituting these values into
\eqref{eq:positivity_alpha} gives
\eqref{eq:contact_upwind}.
\end{proof}

A stationary contact is therefore preserved exactly.
An advected contact is transported with the numerical
diffusion of the first-order upwind scheme, determined
by the material velocity rather than the acoustic speed.
The equilibrium is preserved even though the
volume fraction profile may broaden.
This establishes requirement (R4).

\subsection{Discrete product rule and equivalence of the two formulations}\label{sec:equivalence}

We finally return to the original variables
\(\varphi_a=\rho\alpha_a\).
Proposition~\ref{prop:continuous_equivalence} states that the formulations in \(\varphi_a\) and in \(\alpha_a\) are equivalent for smooth solutions, through the product rule \(\mathcal R_{\varphi_a}=\alpha_a\mathcal R_\rho+\rho\,\mathcal R_{\alpha_a}\). However, at the discrete level, a numerical scheme written in one of the two sets of variables does not automatically satisfy a discrete version of this identity.
We now show that, in the present framework, the discrete product rule is a further application of the side-wise transfer of Sections~\ref{subsec:abgrall} and~\ref{subsec:transfer}: every contribution of a face to a cell is transferred to the other set of variables through the chain rule of that cell. The conservative flux of \(\varphi_a\) and the average velocity at the faces follow as consequences, and the thermodynamic correction can act on all the components.

\paragraph{The product rule as a side-wise chain rule}\label{subsec:product_rule}

We consider the volume fraction scheme of Sections~\ref{sec:framework}--\ref{sec:viscosity} and write its update of cell \(i\) as a sum of contributions,
\[
	\frac{d\qq_i}{dt}=\sum_{j\in\mathcal N_i}\bm\Psi^{(i)}_{ij}+\S(\qq_i),
	\qquad
	\bm\Psi^{(i)}_{ij}=-\frac{|\partial\Omega_{ij}|}{|\Omega_i|}\Bigl(\bm\Phi_{ij}-\mathcal G^{(i)}_{ij}-\mathcal P^{(i)}_{ij}\Bigr).
\]
The map from \(\qq\) to \(\qq^\varphi\) changes only the
volume fraction components.
Let
\(\mathbf J^\varphi_i
=\partial\qq^\varphi/\partial\qq\)
be its Jacobian at the state of cell \(i\).
Its non-trivial rows are
\begin{equation}\label{eq:J_phi}
    \delta\varphi_a
    =
    \alpha_{a,i}\,\delta\rho
    +\rho_i\,\delta\alpha_a.
\end{equation}

\begin{definition}[Transferred scheme]\label{def:transferred_scheme}
The scheme for \(\qq^\varphi\) associated with the volume fraction scheme is
\begin{equation}\label{eq:transferred_scheme}
	\frac{d\qq^\varphi_i}{dt}=\sum_{j\in\mathcal N_i}\mathbf J^\varphi_i\,\bm\Psi^{(i)}_{ij}+\mathbf J^\varphi_i\,\S(\qq_i),
\end{equation}
where the state on the right-hand side is recovered
from \(\qq^\varphi_i\) using
\(\alpha_{a,i}=\varphi_{a,i}/\rho_i\).
\end{definition}

All other components are updated as before. For \(\varphi_a\), the contribution of face \(ij\) is \(\alpha_{a,i}\Psi^{(i),\rho}_{ij}+\rho_i\Psi^{(i),\alpha_a}_{ij}\); applied to the sources, the same rule gives \(\rho\Psi_a=\Phi_a\), since the mixture density has no source, and recovers the pressure-relaxation source of the original formulation.

\begin{lemma}[Discrete product rule]\label{lem:discrete_product_rule}
The transferred scheme satisfies, in every cell and for every \(a<\Nph\),
\begin{equation}\label{eq:discrete_product_rule}
    \frac{d\varphi_{a,i}}{dt}
    =
    \rho_i\frac{d\alpha_{a,i}}{dt}
    +\alpha_{a,i}\frac{d\rho_i}{dt}.
\end{equation}
\end{lemma}
\begin{proof}
Apply \eqref{eq:J_phi} to each face contribution and
to the source in \eqref{eq:transferred_scheme}, and sum.
The resulting two sums are precisely the right-hand
sides of the equations for \(\rho_i\) and
\(\alpha_{a,i}\).
\end{proof}

The product rule is thus a statement about each contribution, on each side of each face, and not about the fluxes. The question is what the transferred scheme looks like as a discretization of the original equations \eqref{eq:system_phi}.

\paragraph{The uncorrected reversible part}\label{subsec:phi_reversible}
Consider a face oriented from cell \(i\) to cell \(j\),
and write
\(V_{\mathbf n}=\VV\cdot\mathbf n_{ij}\).
Before the thermodynamic correction, the density
fluctuation is
\(\tfrac12\Delta(\rho V_{\mathbf n})\),
and the volume fraction fluctuation is
\(\tfrac12\widetilde V_{\mathbf n}\Delta\alpha_a\)
on both sides.
The midpoint rule \eqref{eq:central_ncp} gives
\begin{equation}\label{eq:roe_velocity}
    \widetilde V_{\mathbf n}
    =
    \VV(\widetilde\qq_{ij})\cdot\mathbf n_{ij}
    =
    \frac{(\UU_i+\UU_j)\cdot\mathbf n_{ij}}
         {\rho_i+\rho_j}.
\end{equation}
In particular,
\(
    \widetilde V_{\mathbf n}(\rho_i+\rho_j)
    =
    (\rho V_{\mathbf n})_i
    +(\rho V_{\mathbf n})_j.
\)
The transferred fluctuation of \(\varphi_a\) in cell
\(i\) is therefore
\[
    \widetilde\Phi^{\varphi_a}_{ij}
    =
    \tfrac12\alpha_{a,i}
        \Delta(\rho V_{\mathbf n})
    +
    \tfrac12\rho_i\widetilde V_{\mathbf n}
        \Delta\alpha_a,
\]
and the fluctuation in cell \(j\) has the same expression
with the coefficients \(\alpha_{a,j}\) and \(\rho_j\).

\begin{proposition}[Conservative central flux for
\(\varphi_a\)]\label{prop:phi_flux}
These transferred fluctuations are generated by the
single-valued numerical flux
\begin{equation}\label{eq:phi_flux}
    \mathcal F^{\varphi_a}_{ij}
    =
    \tfrac12\widetilde V_{\mathbf n}
    \bigl(
        \alpha_{a,i}\rho_j
        +\rho_i\alpha_{a,j}
    \bigr).
\end{equation}
More precisely,
\[
    \widetilde\Phi^{\varphi_a}_{ij}
    =
    \mathcal F^{\varphi_a}_{ij}
    -(\varphi_a V_{\mathbf n})_i,
    \qquad
    \widetilde\Phi^{\varphi_a}_{ji}
    =
    (\varphi_a V_{\mathbf n})_j
    -\mathcal F^{\varphi_a}_{ij}.
\]
\end{proposition}

\begin{proof}
Add the physical flux of cell \(i\) to its transferred
fluctuation:
\[
\begin{aligned}
    (\varphi_a V_{\mathbf n})_i
    +\widetilde\Phi^{\varphi_a}_{ij}
    &=
    \tfrac12\alpha_{a,i}
    \bigl[
        (\rho V_{\mathbf n})_i
        +(\rho V_{\mathbf n})_j
    \bigr]
    +\tfrac12\rho_i\widetilde V_{\mathbf n}
        (\alpha_{a,j}-\alpha_{a,i})
    \\
    &=
    \tfrac12\widetilde V_{\mathbf n}
    \bigl(
        \alpha_{a,i}\rho_j
        +\rho_i\alpha_{a,j}
    \bigr).
\end{aligned}
\]
Subtracting the transferred fluctuation of cell \(j\)
from its physical flux gives the same expression.
Moreover, reversing the face orientation changes the
sign of \(\widetilde V_{\mathbf n}\), while leaving the
factor in parentheses unchanged.
Hence
\(\mathcal F^{\varphi_a}_{ji}
=-\mathcal F^{\varphi_a}_{ij}\).
\end{proof}

The density-weighted average \eqref{eq:roe_velocity} is essential to this result. It is not an additional choice made for the \(\varphi_a\)-equation, but the evaluation of \(\UU/\rho\) at the midpoint of the evolved states; an arithmetic average of the velocities would not, in general, give a conservative flux.

\paragraph{Thermodynamic correction}
The correction of Section~\ref{subsec:abgrall} acts on
all components.
Since \(\rho\) and \(\alpha_a\) are linear functions of
\(Z\), its contributions to their reversible
fluctuations are
\(-\tfrac12\theta_{ij}\Delta\rho\) and
\(-\tfrac12\theta_{ij}\Delta\alpha_a\) in cell \(i\),
with the opposite signs in cell \(j\).
By the side-wise product rule, and since \(\alpha_{a,i}\Delta\rho+\rho_i\Delta\alpha_a=\Delta\varphi_a-\Delta\rho\,\Delta\alpha_a\) and \(\alpha_{a,j}\Delta\rho+\rho_j\Delta\alpha_a=\Delta\varphi_a+\Delta\rho\,\Delta\alpha_a\), the transferred corrections of the \(\varphi_a\)-fluctuations are
\[
    \mathcal C^{\varphi_a}_{ij}
    =
    -\tfrac12\theta_{ij}\Delta\varphi_a
    +\tfrac12\theta_{ij}\Delta\rho\,\Delta\alpha_a,
    \qquad
    \mathcal C^{\varphi_a}_{ji}
    =
    +\tfrac12\theta_{ij}\Delta\varphi_a
    +\tfrac12\theta_{ij}\Delta\rho\,\Delta\alpha_a .
\]
The first term is a conservative flux correction; the second, equal in the two cells, is the difference between the finite jump of a product and its differential at either endpoint. The correction thus acts on all components, with an unchanged energy contribution (see below), but its \(\varphi_a\)-component is not in general a conservative flux.

\paragraph{Dissipation and entropy production}
For the configurational viscosity, we compose the
transfer from \(Z\) to \(\qq\) with the map from
\(\qq\) to \(\qq^\varphi\).
The Jacobian of the combined map is
\(\mathbf J^\varphi_s\mathbf J_s\), evaluated separately
at each side \(s\) of the face.
The resulting contributions are
\begin{equation}\label{eq:phi_dissipation}
\begin{aligned}
    \mathcal G^{(i)}_{\varphi_a}
    =
    \alpha_{a,i}\,\mathcal G_\rho
    +\rho_i\,\mathcal G_{\alpha_a},
    \qquad
    \mathcal G^{(j)}_{\varphi_a}
    =
    -\bigl(
        \alpha_{a,j}\,\mathcal G_\rho
        +\rho_j\,\mathcal G_{\alpha_a}
    \bigr),
\end{aligned}
\end{equation}
where the quantities on the right are oriented from
\(i\) to \(j\):
\[
    \mathcal G_\rho
    =
    \tfrac12\,s_{ac}
    \sum_b
    \bigl(
        b_{b,R}\,\crho_{b,R}
        -b_{b,L}\,\crho_{b,L}
    \bigr),
    \qquad
    \mathcal G_{\alpha_a}
    =
    \tfrac12\,s_f\Delta\alpha_a.
\]

Equation~\eqref{eq:phi_dissipation} is not an independent prescription for the viscosity of \(\varphi_a\), but the chain rule for the exchange already constructed in the phase variables; as for the relative velocity rows, each cell uses its own Jacobian. Since the map \(\qq\mapsto\qq^\varphi\) leaves the entropy densities unchanged and the production vector has no density or volume fraction components, the phase entropy productions are exactly those of Section~\ref{subsec:energy_production}.

\paragraph{Semi-discrete equivalence}
\label{subsec:discrete_equivalence}

\begin{theorem}[Semi-discrete equivalence]
\label{thm:discrete_equivalence}
Let the volume fraction scheme be that of
Sections~\ref{sec:framework}--\ref{sec:viscosity},
including the thermodynamic correction on all
components, and let the scheme for \(\qq^\varphi\)
be \eqref{eq:transferred_scheme}.
For corresponding initial data, the two semi-discrete
systems have corresponding solutions,
\[
    \varphi_{a,i}(t)=\rho_i(t)\alpha_{a,i}(t),
\]
with identical values of all other state variables,
for as long as the solutions remain admissible and
\(\rho_i>0\).

The transfer preserves the face-wise energy identities
and the phase entropy productions.
Consequently, under the assumptions of the preceding
results, the corresponding semi-discrete solutions have
the same energy balance, entropy production, phase
masses, volume fractions and contact equilibria.
\end{theorem}

\begin{proof}
Differentiate the change of variables along a solution
of the volume fraction scheme.
The chain rule gives
\eqref{eq:transferred_scheme}, and
Lemma~\ref{lem:discrete_product_rule} gives
\(\varphi_{a,i}=\rho_i\alpha_{a,i}\) at all times
if this relation holds initially.
Conversely, for \(\rho_i>0\), the map is invertible
through \(\alpha_{a,i}=\varphi_{a,i}/\rho_i\).
Applying the inverse change of variables to a solution
of the transferred scheme recovers the
volume fraction scheme.

For the energy, the forces satisfy
\[
    \pp_i
    =
    (\mathbf J^\varphi_i)^{\mathsmaller T}
    \pp^\varphi_i.
\]
Hence every face contribution obeys
\[
    \pp^\varphi_i\cdot
    \mathbf J^\varphi_i\bm\Psi^{(i)}_{ij}
    =
    \pp_i\cdot\bm\Psi^{(i)}_{ij},
\]
and the same identity holds for the source.
The work of each contribution is therefore identical
in the two representations.
The numerical energy fluxes and the energy defects
can be taken to be the same, while the entropy
productions are unchanged because their state
components are unchanged.
The remaining properties follow from the
correspondence of the physical states.
\end{proof}

\paragraph{Conservation and time integration}

Two distinctions are important.
First, semi-discrete equivalence does not mean that
every row of the transferred scheme is conservative.
The uncorrected reversible part has the common flux
\eqref{eq:phi_flux}, but the correction and the
dissipation satisfy
\[
    \mathcal C^{\varphi_a}_{ij}
    +\mathcal C^{\varphi_a}_{ji}
    =
    \theta_{ij}\Delta\rho\,\Delta\alpha_a,
\]
and
\[
    \mathcal G^{(i)}_{\varphi_a}
    +\mathcal G^{(j)}_{\varphi_a}
    =
    -\Delta\alpha_a\,\mathcal G_\rho
    -\Delta\rho\,\mathcal G_{\alpha_a}.
\]
These sums need not vanish, so the transferred scheme is not, in general, conservative in \(\varphi_a\), even without relaxation. This is the price of reproducing the volume fraction scheme through the cell-wise product rule; it does not affect the conservation of the phase masses, mixture momentum and total energy, and for small jumps the additional terms are at least quadratic in the jumps.

Second, the equivalence concerns the semi-discrete
systems.
A time integrator applied independently in the two
representations need not preserve the product rule.
For example, with semi-discrete rates \(L_\rho\) and \(L_{\alpha_a}\) in one cell, the product of two forward Euler updates contains the term \((\Delta t)^2L_\rho L_{\alpha_a}\), which a forward Euler update of the transferred \(\varphi_a\)-equation does not.
The same lack of invariance under a nonlinear change
of variables applies in general to Runge--Kutta
methods.
Fully discrete equivalence, and the transfer of
time-step-dependent positivity statements, therefore
require a compatible time update rather than two
independent applications of the same integrator.

The transfers from \(Z\) to \(\qq\) and from \(\qq\) to \(\qq^\varphi\) thus use the same cell-wise chain rule, which fixes each contribution and preserves its work against the thermodynamic forces: the energy compatibility and the non-negative production established in the phase variables are retained in the other representations.
\section{Numerical results}\label{sec:results}\label{subsec:time_integration}

All computations use the scheme of Sections~\ref{sec:framework}--\ref{sec:viscosity} on uniform Cartesian grids. When relaxation sources are present, the semi-discrete system is integrated with the implicit--explicit scheme IMEX-SSP3(3,3,2) of Pareschi and Russo \cite{PareschiRusso2005}. The fluxes, non-conservative products and numerical
viscosity are treated explicitly, while the algebraic
relaxation sources are treated implicitly in each cell.
The complete IMEX method is second-order accurate.
In the absence of sources, it reduces to the
third-order SSP Runge--Kutta method
\cite{ShuOsher1988,GottliebShuTadmor2001}.
The time step follows a CFL condition based on the speeds \(s_{ac}\) and on the exchange rates of Section~\ref{subsec:transfer}.
The restrictions are checked at the stages, and a
step is repeated with a smaller time step when needed.
For the phase masses, which have no relaxation source,
the positivity argument of Theorem~\ref{thm:positivity}
applies to the SSP stage combinations under the stated
forward Euler bounds, provided the stage states remain
admissible.
Pressure relaxation also changes the volume fractions;
their admissibility during the implicit updates is
monitored, but is not covered by that theorem.
The reference results labelled \emph{Rusanov} use the same reversible part, sources and time integrator, with the scalar dissipation \(\frac12s_{ac}\Delta\qq\) of Section~\ref{subsec:viscous_requirements} in place of the configurational viscosity.

\subsection{Convergence: a stationary two-phase vortex}\label{subsec:test_vortex}
We first measure the accuracy on a smooth, stationary
two-phase vortex without relaxation.
Following the construction of \cite{Thomann2023},
adapted to the present model without thermal impulses,
two ideal gases
(\(\gamma_1=1.4\), \(\gamma_2=2\),
\(c_{v,1}=c_{v,2}=2\))
rotate with azimuthal velocities
\[
    v_a(r)=r\,v_a^0 e^{\nu_a(1-r^k)},
    \qquad
    (v_1^0,\nu_1)=(1,1),
    \qquad
    (v_2^0,\nu_2)=(1.1,1.2),
\]
and volume fraction
\(\alpha_1(r)=0.1+0.33e^{1-r^k}\), with \(k=2\).
The specific entropies are uniform.
The intrinsic densities are obtained by integrating
the radial momentum balance and the stationary
relative-velocity equation \eqref{eq:system_phi_w}.
The domain is \([-5,5]^2\), with periodic boundaries.
We compute the solution up to \(t=0.1\) on
\(N\times N\) grids, with \(N=48,96,192,384\),
using no numerical viscosity, the configurational
viscosity, or the Rusanov viscosity.

\begin{table}[!htbp]
\centering
\small
\caption{Two-phase vortex, \(k=2\), \(t=0.1\): \(L_2\) errors and observed orders of accuracy.}
\label{tab:vortex}
\footnotesize
\setlength{\tabcolsep}{3.5pt}
\begin{tabular}{llcccccccc}
\toprule
& \(N\) & \(\alpha_1\) & order & \(\rho_1\) & order & \(v_{1,x}\) & order & \(v_{2,x}\) & order \\
\midrule
no viscosity & 48 & \(1.68\times10^{-3}\) & -- & \(2.21\times10^{-2}\) & -- & \(8.21\times10^{-3}\) & -- & \(1.22\times10^{-2}\) & -- \\
 & 96 & \(4.67\times10^{-4}\) & 1.85 & \(5.62\times10^{-3}\) & 1.97 & \(2.24\times10^{-3}\) & 1.87 & \(2.92\times10^{-3}\) & 2.06 \\
 & 192 & \(1.20\times10^{-4}\) & 1.97 & \(1.41\times10^{-3}\) & 1.99 & \(5.73\times10^{-4}\) & 1.97 & \(7.33\times10^{-4}\) & 1.99 \\
 & 384 & \(3.01\times10^{-5}\) & 1.99 & \(3.53\times10^{-4}\) & 2.00 & \(1.44\times10^{-4}\) & 1.99 & \(1.84\times10^{-4}\) & 2.00 \\
\midrule
configurational & 48 & \(2.32\times10^{-2}\) & -- & \(2.00\times10^{-1}\) & -- & \(1.72\times10^{-1}\) & -- & \(2.52\times10^{-1}\) & -- \\
 & 96 & \(1.16\times10^{-2}\) & 1.00 & \(1.07\times10^{-1}\) & 0.89 & \(8.93\times10^{-2}\) & 0.95 & \(1.43\times10^{-1}\) & 0.82 \\
 & 192 & \(5.83\times10^{-3}\) & 0.99 & \(5.57\times10^{-2}\) & 0.95 & \(4.54\times10^{-2}\) & 0.98 & \(7.92\times10^{-2}\) & 0.85 \\
 & 384 & \(2.93\times10^{-3}\) & 1.00 & \(2.83\times10^{-2}\) & 0.97 & \(2.29\times10^{-2}\) & 0.99 & \(4.21\times10^{-2}\) & 0.91 \\
\midrule
Rusanov & 48 & \(7.03\times10^{-2}\) & -- & \(2.20\times10^{-1}\) & -- & \(1.97\times10^{-1}\) & -- & \(2.46\times10^{-1}\) & -- \\
 & 96 & \(3.60\times10^{-2}\) & 0.96 & \(1.19\times10^{-1}\) & 0.89 & \(1.04\times10^{-1}\) & 0.92 & \(1.29\times10^{-1}\) & 0.93 \\
 & 192 & \(1.82\times10^{-2}\) & 0.98 & \(6.17\times10^{-2}\) & 0.95 & \(5.35\times10^{-2}\) & 0.96 & \(6.61\times10^{-2}\) & 0.96 \\
 & 384 & \(9.15\times10^{-3}\) & 0.99 & \(3.15\times10^{-2}\) & 0.97 & \(2.72\times10^{-2}\) & 0.98 & \(3.35\times10^{-2}\) & 0.98 \\
\bottomrule
\end{tabular}
\end{table}

\begin{figure}[!htbp]
\centering
\begin{subfigure}[b]{0.49\textwidth}
\centering
\includegraphics[scale=\figscale]{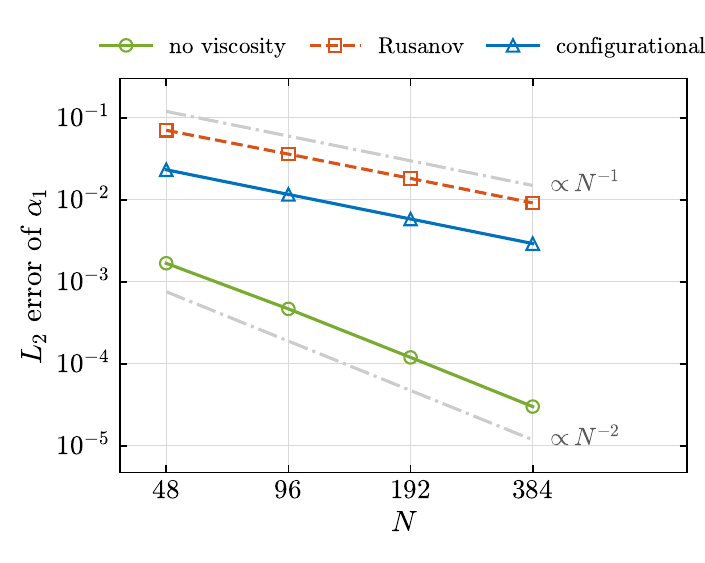}
\caption{\(L_2\) error of \(\alpha_1\)}
\end{subfigure}\hfill
\begin{subfigure}[b]{0.49\textwidth}
\centering
\includegraphics[scale=\figscale]{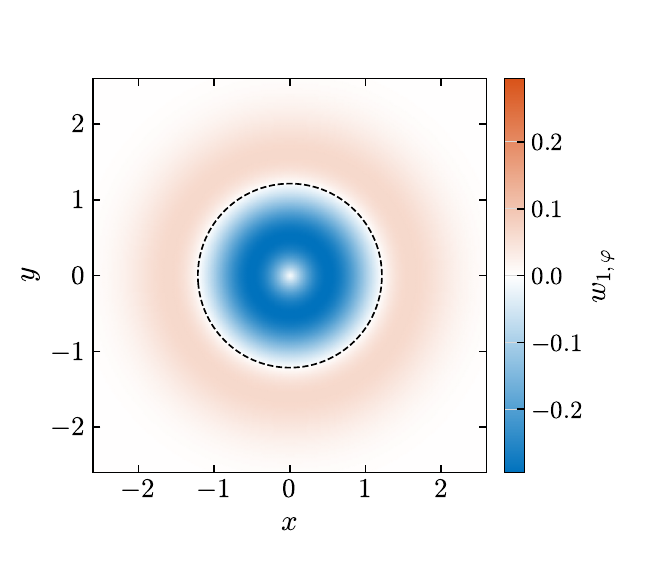}
\caption{relative velocity \(w_{1,\varphi}\)}
\end{subfigure}
\caption{Two-phase vortex at \(t=0.1\). Left: \(L_2\) error of the volume fraction versus the number of cells per direction. Right: azimuthal component of the relative velocity \(\ww_1=\vvv_1-\vvv_2\) on the \(384\times384\) grid without numerical dissipation; the dashed circle is the locus \(\ww_1=\mathbf 0\).}
\label{fig:vortex_alpha}
\end{figure}

Table~\ref{tab:vortex} shows the
\(L_2\) errors measured against the initial
profiles of the primitive variables and the orders of convergence.
Without viscosity, the accuracy of the
reversible discretization shows second-order convergence.
With either viscosity, the observed orders approach
one, as expected for a first-order dissipative scheme.
The main difference is in the volume fraction:
using \(s_f\) instead of the acoustic scale reduces
its error by a factor of about three on every grid
(Figure~\ref{fig:vortex_alpha}, left).
The errors in the other variables remain comparable.

In (Figure~\ref{fig:vortex_alpha}, on the right), it is interesting to note that the test is designed such that the phases are out of phase with one another, with the
azimuthal relative velocity changing sign in the vicinity of
\(r=1.22\)
Although the energy is non-convex in \(\qq\) at states
with slip, no negative Rusanov production is observed
in this test.
Non-convexity does not require the pairing of every
neighbouring pair of states to be negative.
The shock tube of Section~\ref{subsec:test_shocktube}
provides a case where negative pairings do occur.

\FloatBarrier
\subsection{Transported contact}\label{subsec:test_contact}
We next consider the contact equilibrium of
Theorem~\ref{thm:contact}.
Two ideal gases
(\(\gamma_1=1.4\), \(\gamma_2=1.6\))
have uniform intrinsic densities
\(\rho_1=1\), \(\rho_2=0.125\),
common velocity \(u=1\) and common pressure \(p=1\).
The volume fraction changes from
\(\alpha_1=1-\epsilon\) to \(\alpha_1=\epsilon\)
across \(x=-0.2\), with
\(\epsilon=3\times10^{-5}\).
The transition is smoothed on the scale of one cell, leaving the intrinsic
states unchanged.
The domain is \([-0.5,0.5]\), discretized with
\(288\) cells in \(x\) and three identical rows in \(y\).
Relaxation processes are not taken into account, and the two viscosities
are compared at \(t=0.05\).

\begin{figure}[!htbp]
\centering
\begin{minipage}[c]{0.47\textwidth}
\centering
\includegraphics[scale=\figscale]{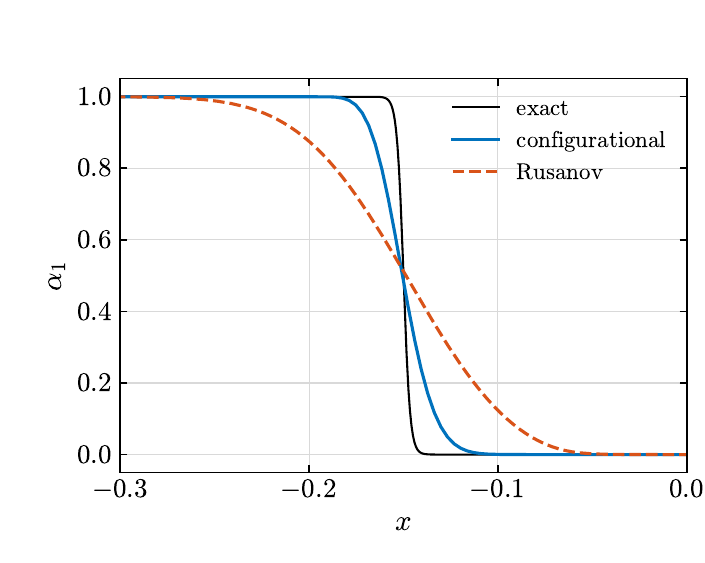}
\end{minipage}\hfill
\begin{minipage}[c]{0.51\textwidth}
\centering
\footnotesize
\begin{tabular}{lcc}
\toprule
& configurational & Rusanov \\
\midrule
\(L_1\) error of \(\alpha_1\) & \(8.49\times10^{-3}\) & \(3.19\times10^{-2}\) \\
\(\max|p-1|\) & \(4.8\times10^{-14}\) & \(1.9\times10^{-13}\) \\
\(\max|u-1|\) & \(4.7\times10^{-14}\) & \(3.1\times10^{-13}\) \\
\(\max|p_1-p_2|\) & \(1.3\times10^{-9}\) & \(5.8\times10^{-10}\) \\
\(\max|\ww_1|\) & \(9.0\times10^{-9}\) & \(1.5\times10^{-9}\) \\
\bottomrule
\end{tabular}
\end{minipage}
\caption{Transported contact at \(t=0.05\), \(288\) cells. Left: volume fraction near the interface against the translated initial profile. Right: \(L_1\) error of the volume fraction and largest deviations from the contact equilibrium.}
\label{fig:contact}
\end{figure}

Both computations preserve the contact equilibrium
with small numerical deviations
(Figure~\ref{fig:contact}, right).
The mixture pressure and velocity remain uniform
to about \(10^{-13}\); the phase pressure difference
and relative velocity remain of order \(10^{-9}\)
and \(10^{-8}\), respectively.
The Rusanov viscosity also preserves the equilibrium
because its scalar action keeps the states in the
same affine contact family.
The difference lies in the diffusion of the interface.
Both schemes transport it at the speed of the material,
but Rusanov's viscosity diffuses it on the acoustic timescale.
Configurational viscosity, on the other hand, reduces to
the upwind scheme at the speed of the material described in
Theorem~\ref{thm:contact}.
The partial densities remain positive in both runs.

\FloatBarrier
\subsection{Multiphase shock tube}\label{subsec:test_shocktube}
We consider a Sod-type shock tube for the two gases
of Section~\ref{subsec:test_contact}, without relaxation.
The left state is
\((\rho_1,\rho_2,p)=(1,1,1)\),
with \(\alpha_1=1-\epsilon\), and the right state is
\((0.125,0.125,0.1)\), with \(\alpha_1=\epsilon\),
where \(\epsilon=3\times10^{-5}\).
The domain \([-0.5,0.5]\) is discretized with
\(576\) cells in \(x\) and three identical rows in \(y\),
and the final time is \(t=0.2\).

\begin{figure}[!htbp]
\centering
\includegraphics[scale=\figscale]{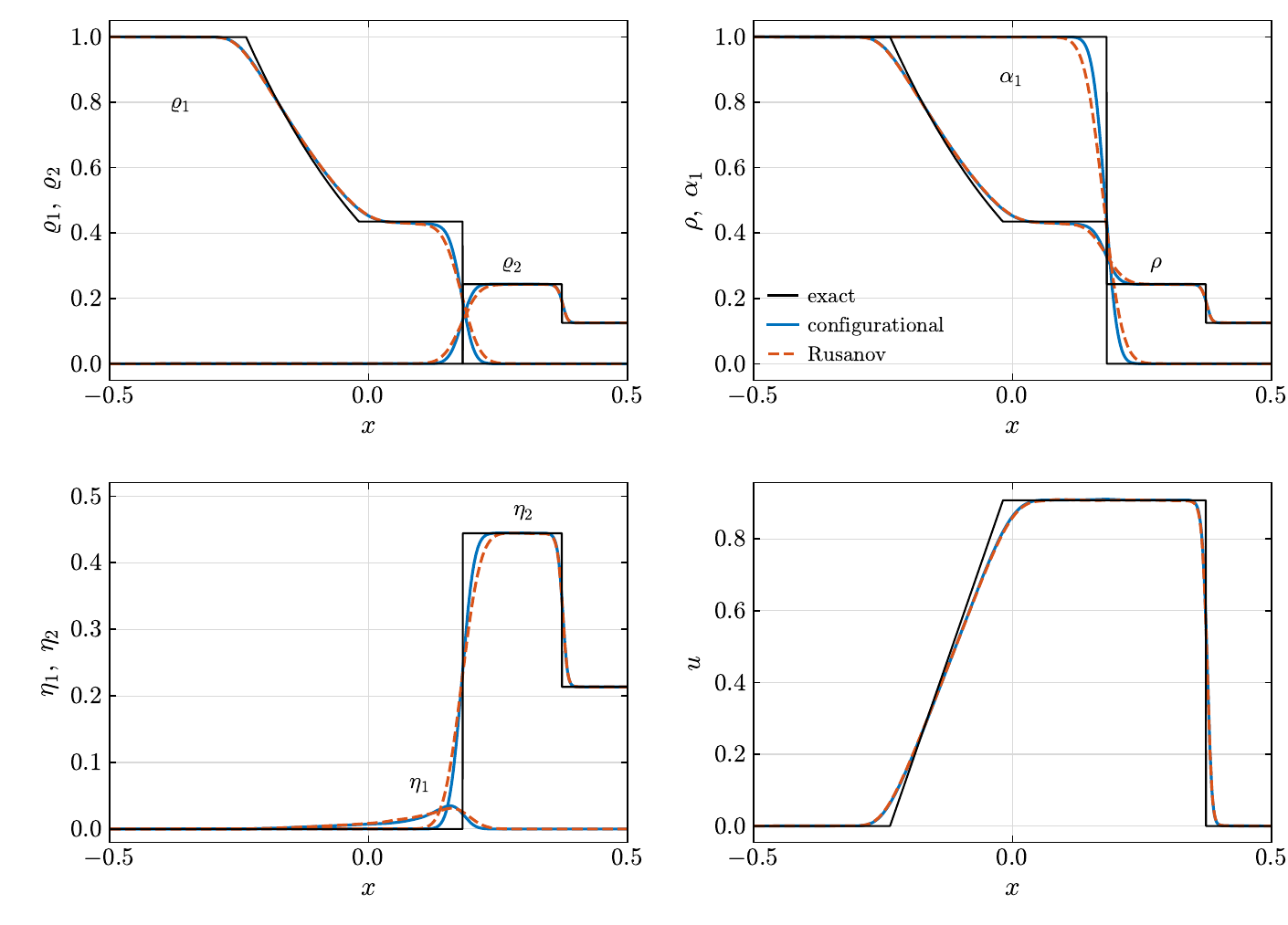}
\caption{Multiphase shock tube at \(t=0.2\),
\(576\) cells: partial densities \(\crho_1\) and
\(\crho_2\) (top left), mixture density and volume
fraction of phase~1 (top right), entropy densities
\(\eta_1\) and \(\eta_2\) (bottom left), and mixture
velocity (bottom right).
The reference curves are the Riemann solution for
two pure gases with a sharp initial discontinuity; the
left state has zero entropy, so that \(\eta_1\) vanishes
in the reference solution.}
\label{fig:rp13_profiles}
\end{figure}

Both viscosities reproduce the wave structure of the
exact solution for two pure gases
(Figure~\ref{fig:rp13_profiles}).
The configurational viscosity gives a narrower volume
fraction transition, and its \(L_1\) differences from
the reference in density, velocity and pressure are
also smaller.
The partial densities and the entropy densities of the
two phases follow the reference, with a sharper
transition at the contact for the configurational
viscosity.
The volume fractions and partial densities remain
positive in both computations.

\begin{figure}[!htbp]
\centering
\begin{subfigure}[b]{4.75cm}
\centering
\includegraphics[scale=\figscale]{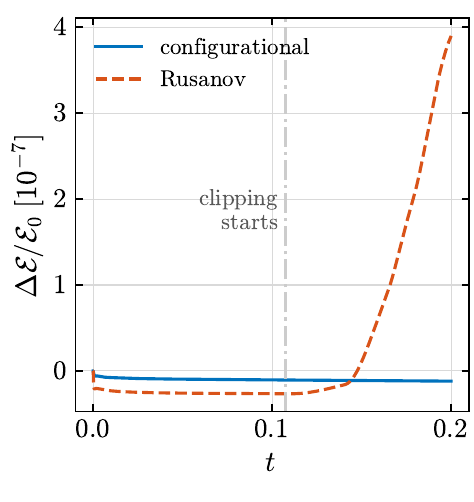}
\caption{total energy}\label{fig:rp13_energy}
\end{subfigure}\hfill
\begin{subfigure}[b]{5.1cm}
\centering
\includegraphics[scale=\figscale]{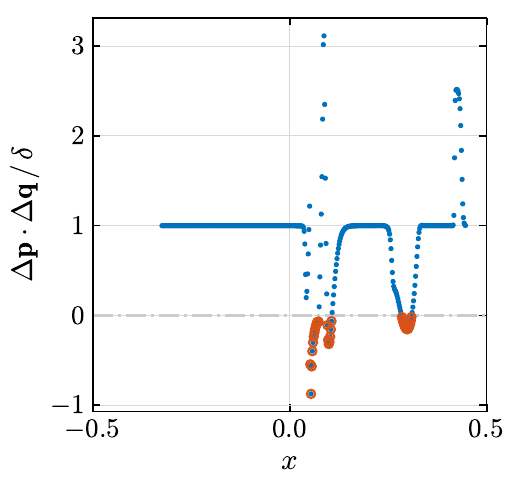}
\caption{\(\Delta\pp\cdot\Delta\qq/\delta_{ij}\) at \(t=0.2\)}\label{fig:rp13_ratio}
\end{subfigure}\hfill
\begin{subfigure}[b]{5.0cm}
\centering
\includegraphics[scale=\figscale]{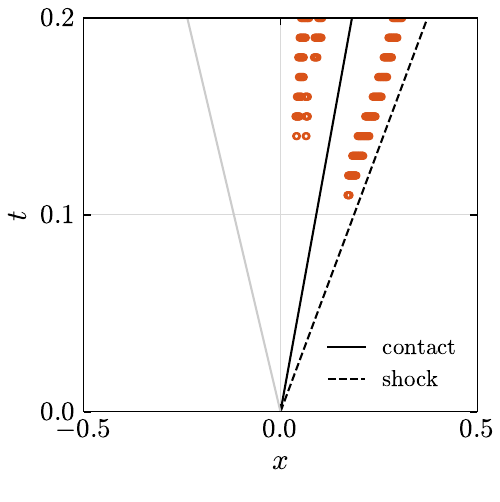}
\caption{negative pairing in \((x,t)\)}\label{fig:rp13_xt}
\end{subfigure}
\caption{Multiphase shock tube, \(576\) cells.
(a) Relative change of the total energy \(\int\rhoE\,dx\); the dash-dotted line marks the first evaluation in which the solver clips a negative production of the Rusanov viscosity.
(b) Rusanov viscosity: ratio of \(\Delta\pp\cdot\Delta\qq\) to the phase-variable pairing \(\delta_{ij}\) at \(t=0.2\).
(c) Faces with negative Rusanov pairing at the output times; the lines show the contact and shock trajectories of the two-gas Riemann reference.}
\label{fig:rp13_thermo}\label{fig:rp13_sign}
\end{figure}

Both viscosities reproduce the entropy increase across
the shock, with similar entropy-density profiles consistent
with the Riemann reference
(Figure~\ref{fig:rp13_profiles}, bottom left).
Compared with the reference, the numerical profiles also
show additional entropy in the rarefaction and near the
contact, where \(\eta_1\), which vanishes in the reference,
takes small positive values.

The local productions, however, differ in sign.
From \(t\approx0.11\), the Rusanov pairing
\(\Delta\pp\cdot\Delta\qq\) becomes negative in
regions where a trace phase slips relative to the
dominant phase; its ratio to the non-negative
phase-variable pairing \(\delta_{ij}\) reaches \(-0.87\)
(Figures~\ref{fig:rp13_ratio} and~\ref{fig:rp13_xt}).
These negative values arise from the remainder
\(T_w\) in \eqref{eq:pairing_split}.
To enforce non-negative entropy production, the solver
clips the resulting negative productions to zero.
This clipping breaks the viscous energy identity:
it removes the term that would compensate the energy
increase caused by the viscous state contributions.
The Rusanov computation consequently develops a positive
energy drift after the onset of clipping
(Figure~\ref{fig:rp13_energy}).
The configurational viscosity instead satisfies the
viscous energy identity with non-negative entropy
production, without clipping.
Its relative total-energy change remains about
\(1.2\times10^{-8}\); the energy compatibility is
semi-discrete and does not imply exact conservation
after time integration.

\subsection{Shock--water column interaction}\label{subsec:test_shockcolumn}
The last test considers a planar shock in air
interacting with a water column
\cite{SaurelAbgrall,SHTCSurfaceTension}.
Air is an ideal gas with \(\gamma=1.4\), and water
a stiffened gas with \(\gamma=4.7\) and
\(p_\infty=4.7\times10^8\,\mathrm{Pa}\).
The unshocked phases are at \(10^5\,\mathrm{Pa}\),
with densities \(1.18\) and
\(998.2\,\mathrm{kg\,m^{-3}}\).
A residual volume fraction of \(3\times10^{-5}\)
is retained for each trace phase.
The column has radius \(3.2\,\mathrm{mm}\) and is
centred at the origin; the incident shock has
Mach number \(1.3\) and starts at
\(x=-3.65\,\mathrm{mm}\).

The domain is \([-20,20]^2\,\mathrm{mm}^2\),
with periodic boundaries in \(y\) and outflow
conditions in \(x\).
Pressure and velocity relaxation are active and stiff, with
\(\tau^{\mathrm p}=\tau^{\mathrm w}=10^{-8}\,\mathrm s\)
(see \cite{ExpInt2} for a dedicated solver for stiff
relaxation in Baer--Nunziato models).
We first compare the two viscosities on a
\(192\times192\) grid up to \(t=8\,\mu\mathrm s\).
Both use the thermodynamic correction built in \(Z\),
unless stated otherwise.

\begin{figure}[!htbp]
\centering
\begin{subfigure}[b]{0.42\textwidth}
\centering
\includegraphics[scale=\figscale]{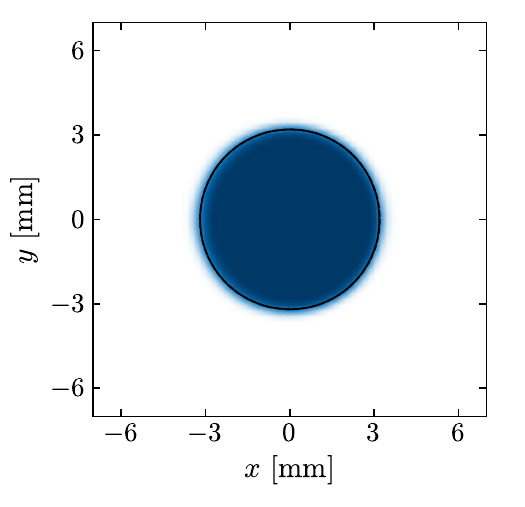}
\caption{configurational viscosity}
\end{subfigure}\hspace{4mm}
\begin{subfigure}[b]{0.50\textwidth}
\centering
\includegraphics[scale=\figscale]{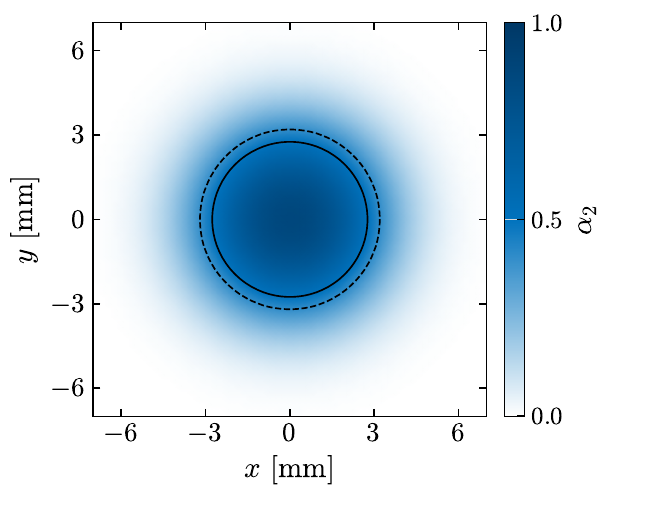}
\caption{Rusanov viscosity}
\end{subfigure}
\caption{Shock--water column interaction at \(t=8\,\mu\mathrm s\), \(192\times192\) cells: volume fraction of water. Dashed: initial interface \(\alpha_2=0.5\); solid: \(\alpha_2=0.5\) at \(t=8\,\mu\mathrm s\).}
\label{fig:ic11_maps}
\end{figure}

The large acoustic speed of water makes the
difference in interface diffusion particularly
visible (Figure~\ref{fig:ic11_maps}).
With the configurational viscosity, the equivalent radius
of the region \(\alpha_2>0.5\) remains at its
initial value on this grid, \(3.18\,\mathrm{mm}\).
With Rusanov, the diffuse is wider and the equivalent radius
decreases to \(2.75\,\mathrm{mm}\).
The partial densities and volume fractions remain
positive in both computations.

The energy computed from the state is compared with
an independently transported energy that accounts
for the boundary fluxes.
At \(t=8\,\mu\mathrm s\), their relative difference
is \(3.5\times10^{-8}\) with the configurational
viscosity, a difference that comes from the time
discretization, and \(1.4\times10^{-7}\) with Rusanov.
As in the shock tube, the Rusanov computation
requires clipping of negative entropy productions;
the configurational viscosity requires none.

\begin{figure}[!htbp]
\centering
\includegraphics[scale=\figscale]{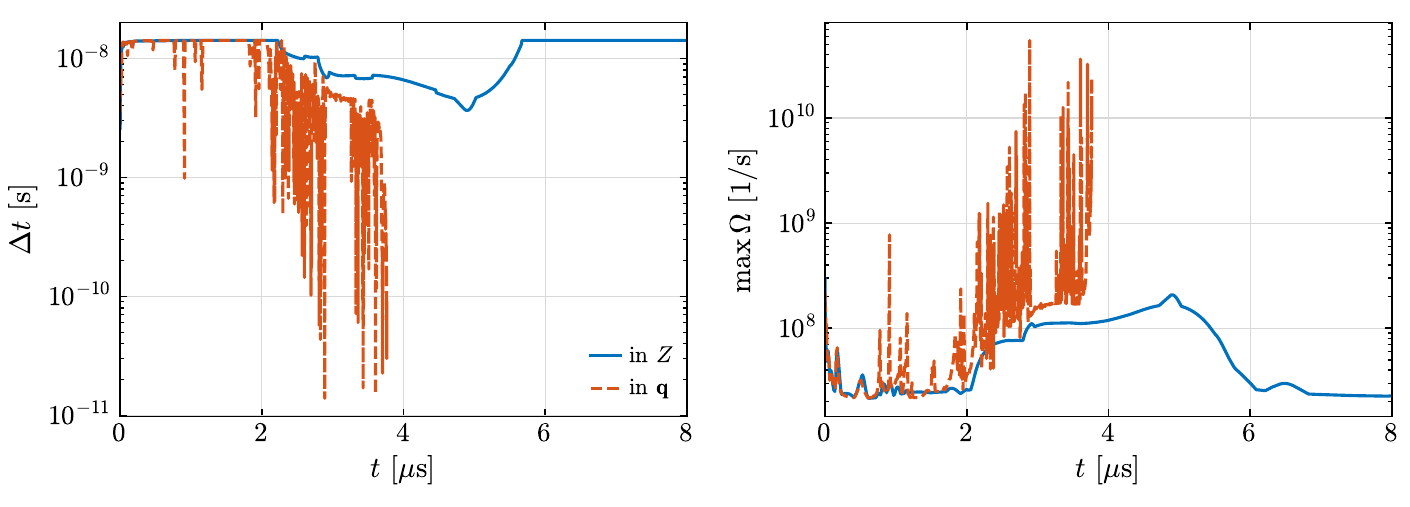}
\caption{Shock--water column interaction with the configurational viscosity: time step (left) and largest exchange rate of the viscosity (right) versus time, with the thermodynamic correction built in the phase variables (solid) and with the original correction in the evolved variables (dashed). Steps shortened to reach the output times are omitted.}
\label{fig:ic11_correction}
\end{figure}

We also repeat the configurational computation with
the original thermodynamic correction along
\(\Delta\qq\).
Near-cancellation in its denominator
\(\Delta\pp\cdot\Delta\qq\) produces large correction
coefficients and increases the dissipative speeds.
After \(t=2\,\mu\mathrm s\), the largest exchange
rate grows by more than two orders of magnitude,
and the time step falls to about
\(10^{-11}\,\mathrm s\)
(Figure~\ref{fig:ic11_correction}).
This run was stopped at \(t=3.8\,\mu\mathrm s\)
after \(982\) steps.
With the direction built in \(Z\), the computation
reaches \(t=8\,\mu\mathrm s\) in \(855\) steps.
Thus, in this test, the phase-variable direction
substantially reduces the time-step restriction
associated with the correction.

\begin{figure}[!htbp]
\centering
\includegraphics[scale=\figscale]{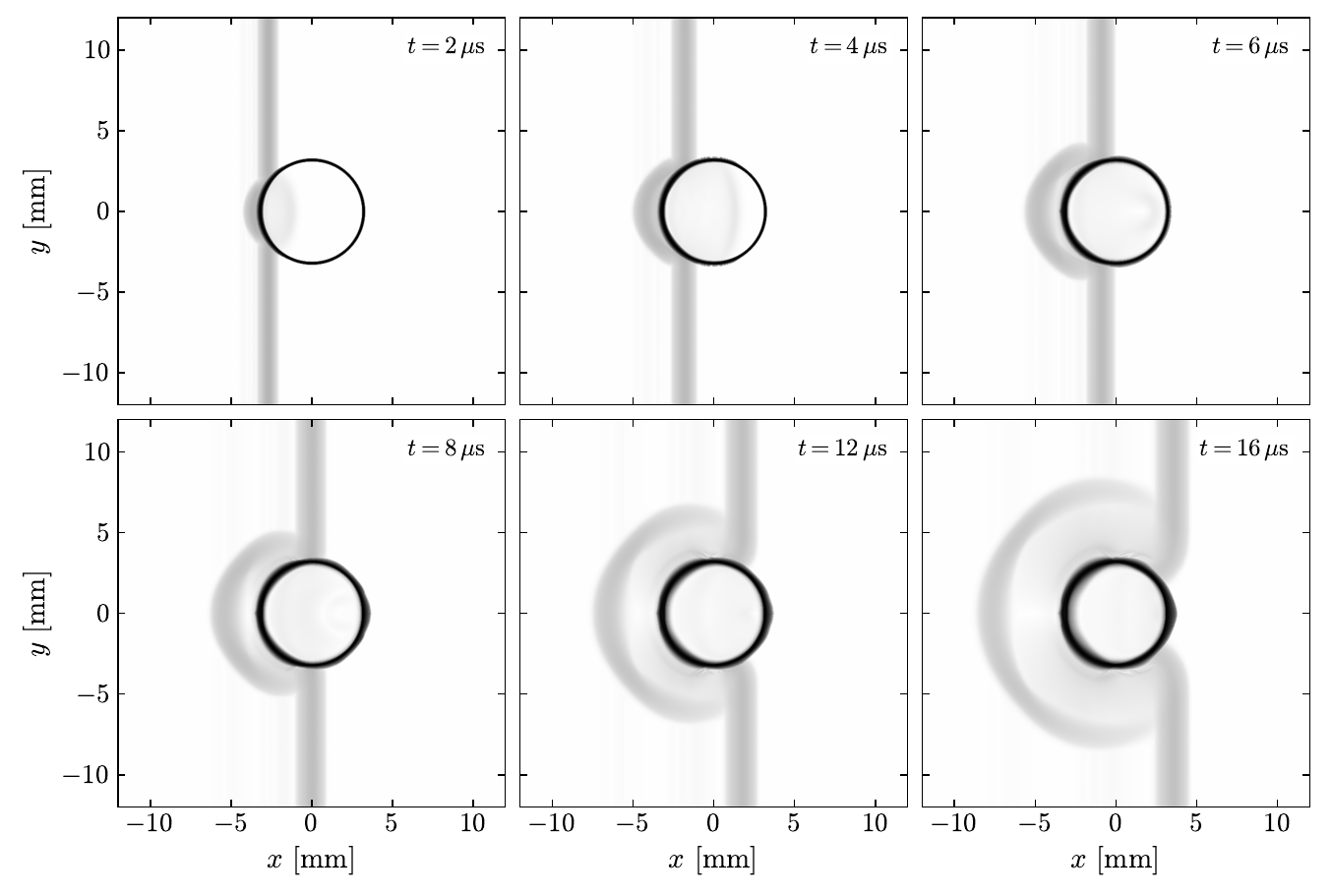}
\caption{Shock--water column interaction on \(1152\times1152\) cells with the configurational viscosity: numerical schlieren image of the mixture density, \(1-\bigl[\ln(1+|\nabla\rho|/g_0)/\ln(1+g_1/g_0)\bigr]^{0.7}\), where \(g_0\) and \(g_1\) are the \(90\)th percentile of \(|\nabla\rho|\) in the air and its \(99.9\)th percentile in the window at \(t=16\,\mu\mathrm s\).}
\label{fig:ic11_schlieren}
\end{figure}

Finally, Figure~\ref{fig:ic11_schlieren} shows an
additional computation on \(1152\times1152\) cells,
continued to \(t=16\,\mu\mathrm s\).
This run also uses a sensor to reduce the acoustic
dissipation away from shocks:
\[
    s_{ac}'=
    s_f+\psi_{ij}(s_{ac}-s_f),
    \qquad
    \psi_{ij}\in[0.15,1],
\]
where \(\psi_{ij}\) increases with
\(|\theta_{ij}|/\lambda_{ij}\).
The inequality \(s_f\le s_{ac}'\le s_{ac}\)
preserves the non-negative exchange weights and
Bregman productions of the viscous construction.
The image shows the incident shock, the wave
reflected by the column, the waves diffracted
around it and those transmitted into the water,
while the interface remains a thin ring.
The relative difference between the state energy
and the transported energy is \(3.1\times10^{-9}\).

\section{Conclusions}\label{sec:conclusions}
The variables used to derive a model, to discretize its
equations and to expose its thermodynamic structure need
not coincide.
We have used this distinction to construct a configurational
numerical viscosity for a multiphase HTC model whose energy
is not convex in the evolved variables.
The fluxes and non-conservative products remain in these
variables, while the viscous exchange and the direction of
the thermodynamic correction are constructed in the phase
variables, where the energy is convex and its configurational
direction is explicit.
Each cell receives these contributions through its own
chain rule.

The jump decomposition separates changes in the amount of
a phase from changes in its intrinsic state, allowing
distinct configurational and acoustic dissipation scales.
The exchange uses non-negative weights shared by the volume,
mass, momentum and entropy of each phase.
Under the stated explicit time-step conditions, it preserves
positive partial densities, admissible volume fractions
and contact equilibria.
The side-wise transfer preserves the work of the
thermodynamic forces, and a common viscous energy flux
gives non-negative Bregman defects in each cell and phase.
Together with the compatibility condition for the reversible
part, this yields the semi-discrete energy balance and
entropy inequalities without clipping.
The same chain rule gives the discrete product rule and
the semi-discrete equivalence with the density-weighted
volume-fraction formulation.

The smooth-vortex test shows second-order convergence for
the reversible discretization and first-order convergence
when viscosity is included.
The configurational viscosity reduces interface spreading
compared with the scalar Rusanov viscosity.
In the shock-tube and shock--water-column tests, the latter,
used within the same HTC framework, develops negative
entropy production whose clipping introduces an energy
defect; configurational viscosity restores full thermodynamic compatibility.
In the water-column test, the phase-variable correction
direction also avoids the severe time-step reduction
observed with the original direction.

The applicability of this construction to the complete
unified HTC model of mixtures \cite{FerrariPeshkovDumbser}
remains to be investigated.
The analysis also motivates further study of the model
itself, including possible formulations of relative motion
that address energy convexity at the modelling level.

\section*{Acknowledgments}
D.F.\ has received funding from the European Union's Horizon Europe research and innovation programme under the Marie Sk\l{}odowska-Curie grant agreement No.~101206752 (project OMNIA, \emph{Novel unified and extended mixture theories: theory and numerics}). Funded by the European Union. Views and opinions expressed are however those of the author(s) only and do not necessarily reflect those of the European Union or the European Research Executive Agency (REA). Neither the European Union nor the granting authority can be held responsible for them.
   M.D.\ was funded by the Fondazione Caritro via the project SOPHOS and by the European Research Council (ERC)
   under the European Union's Horizon Europe research and innovation programme via the project SOPHOS, grant
   agreement no.~ERC-ADG-2025-101265878-SOPHOS.  D.F., I.P.\ and M.D.\ are members of the Gruppo Nazionale per il
   Calcolo Scientifico dell'Istituto Nazionale di Alta Matematica (GNCS-INdAM).

\appendix

\section{Convexity in the different representations}\label{app:convexity}

This appendix collects the convexity statements used in Section~\ref{sec:geometry}. All of them follow from one elementary observation on how a change of variables acts on a Hessian, and the two non-linear changes of variables of the model have the same structure. We consider states with positive partial densities,
volume fractions in \((0,1)\), and positive temperatures,
and assume thermodynamically stable phase equations of
state. We give explicit formulas for two phases; the extension to \(\Nph\) phases is indicated where it is not immediate.

\subsection{How a change of variables acts on convexity}\label{app:change_of_variables}

\begin{lemma}[Hessian under a change of variables]\label{lem:change_of_variables}
Let \(\rhoE\) be a twice differentiable function of \(\bm x\), and let \(\bm x=\bm x(\bm y)\) be a smooth change of variables with Jacobian \(\mathbf J=\partial\bm x/\partial\bm y\). Then
\begin{equation}\label{eq:hessian_change}
	\nabla^2_{\bm y}\rhoE
	=
	\mathbf J^{\mathsmaller T}\,\nabla^2_{\bm x}\rhoE\,\mathbf J
	+
	\sum_k\rhoE_{x_k}\,\nabla^2_{\bm y}x_k .
\end{equation}
\end{lemma}

The proof is the chain rule applied twice. Formula \eqref{eq:hessian_change} separates two effects. If \(\mathbf J\) is invertible, the first is a congruence
and has the same numbers of positive, negative and zero
eigenvalues as \(\nabla^2_{\bm x}\rhoE\), so it transfers convexity, or its absence, unchanged. The second accounts for the curvature of the change of
variables, weighted by the thermodynamic forces.
It vanishes for affine transformations, which therefore
preserve convexity.
For a nonlinear transformation it can change the sign
of the quadratic form.

The two maps used here have a particularly simple
nonlinearity.
In the model, all changes of variables between \(\qq^\varphi\), \(\qq\) and \(Z\) are affine except two, and both are ratios: \(\alpha_a=\varphi_a/\rho\) between \(\qq^\varphi\) and \(\qq\), and \(\vvv_a=\bm M_a/\crho_a\) inside the relative velocities between \(Z\) and \(\qq\). The curvature of a ratio \(f(\zeta,\xi)=\xi/\zeta\) (for scalar \(\xi\), or for each component of a vector \(\xi\)) is
\begin{equation}\label{eq:ratio_hessian}
	\nabla^2f=
	\begin{pmatrix}
		2\xi/\zeta^3 & -1/\zeta^2\\
		-1/\zeta^2 & 0
	\end{pmatrix},
	\qquad
	\det\nabla^2f=-\frac{1}{\zeta^4}<0,
\end{equation}
so it is always indefinite. By \eqref{eq:hessian_change}, if the force multiplying this curvature vanishes, the two Hessians are congruent and a representation inherits the convexity, or the lack of it, of the other. The condition is sufficient, not necessary: away from it the curvature term may destroy convexity, as for the relative velocities, but it may also restore it (Remark~\ref{rem:curvature_restores}). For the redundant variables \(Z\), \eqref{eq:hessian_change} is an identity between quadratic forms on admissible variations, on which \(\mathbf J\) is invertible. We shall see that the relevant forces are the pressure difference for \(\varphi_a\leftrightarrow\alpha_a\) and the relative momentum for \(Z\leftrightarrow\qq\): both vanish at mechanical equilibrium, and both are generically non-zero away from it.

\subsection{A single phase and its perspective}\label{app:single_phase}

\begin{lemma}[Single-phase energy]\label{lem:single_phase}
Let
\[
    f_a(v,s)=e_a(1/v,s)
\]
be the specific internal energy expressed in specific
volume and specific entropy. Then \(f_a\) is convex if and only if \(\varepsilon_a(\rho,\sigma)=\rho\,e_a(\rho,\sigma/\rho)\) is convex. Under this assumption, the single-fluid energy \[
    E_a(Y)
    =\frac{|\bm m|^2}{2\rho}
     +\varepsilon_a(\rho,\sigma)
\] is convex in \(Y=(\rho,\bm m,\sigma)\).
\end{lemma}

\begin{proof}
The two internal-energy representations satisfy
\[
    \varepsilon_a(\rho,\sigma)
    =\rho\,f_a\left(\frac1\rho,\frac{\sigma}{\rho}\right),
    \qquad
    f_a(v,s)
    =v\,\varepsilon_a\left(\frac1v,\frac{s}{v}\right).
\]

One has \(\varepsilon_a(\rho,\sigma)=\rho\,f_a\bigl((1,\sigma)/\rho\bigr)\), i.e.\ \(\varepsilon_a\) is the perspective of \(f_a\) evaluated on the affine set where the first argument equals one; conversely \(f_a(v,s)=v\,\varepsilon_a\bigl((1,s)/v\bigr)\). The perspective of a convex function is convex \cite{BoydVandenberghe2004}, and so is its restriction to an affine set. The kinetic term \(|\bm m|^2/(2\rho)\) is the perspective of \(\frac12|\bm m|^2\).
\end{proof}

Proposition~\ref{prop:perspective} then gives the convexity of \(\rhoE\) in \(Z\). The flat directions of Section~\ref{subsec:perspective} are a consequence of extensivity and not of a particular equation of state: along \(\lambda Z_a\) the energy of the phase is linear, for internal and kinetic energy alike. The energy of the model is therefore convex but never strictly convex in \(Z\); strict convexity would require an additional, strictly convex interfacial contribution.

\subsection{Volume fractions versus density-weighted volume fractions}\label{app:phi_alpha}

\begin{proposition}[Internal energy in \(\qq\) and in \(\qq^\varphi\)]\label{prop:phi_alpha}
Consider two phases with thermodynamically stable equations of state, and let \(\rhoE_{\mathrm{int}}=\sum_a\crho_ae_a(\rho_a,s_a)\).
Then:
\begin{enumerate}
    \item[(i)]
    \(\rhoE_{\mathrm{int}}\) is convex in the variables
    \((\crho_1,\eta_1,\eta_2,\rho,\alpha_1)\)
    of \(\qq\).

    \item[(ii)]
    In the corresponding variables with
    \(\varphi_1=\rho\alpha_1\),
    \begin{equation}\label{eq:hessian_phi}
        \nabla^2_{\qq^\varphi}\rhoE_{\mathrm{int}}
        =
        \mathbf J^{\mathsmaller T}
        \nabla^2_{\qq}\rhoE_{\mathrm{int}}\,\mathbf J
        +(p_2-p_1)
        \nabla^2\left(\frac{\varphi_1}{\rho}\right),
    \end{equation}
    where \(\mathbf J\) is the Jacobian of
    \(\qq^\varphi\mapsto\qq\).
    At pressure equilibrium this Hessian is
    positive semi-definite.

    \item[(iii)]
    If \((p_1-p_2)(\rho_1-\rho_2)<0\), the
    internal-energy Hessian in \(\qq^\varphi\)
    is indefinite.
\end{enumerate}
\end{proposition}

\begin{proof}
\begin{enumerate}
    \item[(i)]
   We write
\[
    \rhoE_{\mathrm{int}}
    =
    \sum_a
    \alpha_a\varepsilon_a
    \left(\frac{\crho_a}{\alpha_a},
          \frac{\eta_a}{\alpha_a}\right).
\]
This is a sum of convex perspectives.
The substitutions
\(\alpha_2=1-\alpha_1\) and
\(\crho_2=\rho-\crho_1\) are affine, so they preserve
convexity.

    \item[(ii)]
    The map \(\qq^\varphi\mapsto\qq\) is the identity except for \(\alpha_1=\varphi_1/\rho\), which is the only nonlinear component of the map.
Its conjugate force, by \eqref{eq:forces_alpha}, is
\[
    \partial_{\alpha_1}\rhoE_{\mathrm{int}}=p_2-p_1.
\]
Equation~\eqref{eq:hessian_phi} follows from
\eqref{eq:hessian_change}.
At \(p_1=p_2\) only the congruence remains, and (i) gives semi-definiteness.

    \item[(iii)]
We consider a variation \(\bm k\) that changes
the volume fractions at fixed intrinsic densities and
specific entropies.
Writing
\(\delta\alpha=\delta\alpha_1=-\delta\alpha_2\),
its components satisfy
\[
\begin{gathered}
    \delta\crho_1=\rho_1\delta\alpha,
    \qquad
    \delta\crho_2=-\rho_2\delta\alpha,
    \qquad
    \delta\rho=(\rho_1-\rho_2)\delta\alpha,
    \\
    \delta\eta_1=\rho_1s_1\delta\alpha,
    \qquad
    \delta\eta_2=-\rho_2s_2\delta\alpha .
\end{gathered}
\]
The internal energy is linear along this direction,
hence
\(
    \bm k^{\mathsmaller T}
    \nabla^2_{\qq}\rhoE_{\mathrm{int}}\,\bm k=0.
\)
The corresponding tangent in the other representation
is \(\bm k^\varphi=\mathbf J^{-1}\bm k\), with
\(
    \delta\varphi_1
    =\alpha_1\delta\rho+\rho\,\delta\alpha.
\)
Using \eqref{eq:ratio_hessian}, the curvature of
\(\varphi_1/\rho\) along this tangent is
\[
    \frac{2}{\rho^2}\delta\rho
    \left(\alpha_1\delta\rho-\delta\varphi_1\right)
    =
    -\frac{2}{\rho}\delta\rho\,\delta\alpha.
\]
The congruence term contributes zero, and therefore
\[
    (\bm k^\varphi)^{\mathsmaller T}
    \nabla^2_{\qq^\varphi}\rhoE_{\mathrm{int}}\,
    \bm k^\varphi
    =
    \frac{2(\delta\alpha)^2}{\rho}
    (p_1-p_2)(\rho_1-\rho_2)<0.
\]

There are also directions of positive curvature.
A variation of \(\eta_1\) alone is unchanged by the
transformation and gives
\[
    \partial^2_{\eta_1}\rhoE_{\mathrm{int}}
    =
    \frac{\partial_s^2e_1}{\crho_1}
    =
    \frac{T_1}{c_{v,1}\crho_1}>0.
\]
Thus the Hessian has both positive and negative
directions.
\end{enumerate}

\end{proof}

For \(\Nph\) phases, the additional term in
\eqref{eq:hessian_phi} becomes
\[
    \sum_{a<\Nph}
    (p_\Nph-p_a)\,
    \nabla^2\left(\frac{\varphi_a}{\rho}\right).
\]
It vanishes when all phase pressures are equal.

\begin{remark}[The curvature term can also restore convexity]\label{rem:curvature_restores}
For \((p_1-p_2)(\rho_1-\rho_2)>0\) the curvature term of \eqref{eq:hessian_phi} is positive along \(\bm k^\varphi\), the direction that is flat in \(\qq\). This does not by itself establish convexity, since
other directions must also be considered.
Considering the total energy, this can compensate the coupling with the mixture momentum that makes the Hessian in \(\qq\) indefinite (Proposition~\ref{prop:nonconvex_q}).
For example, consider
\[
    e_a=\rho_a^{\gamma_a-1}\exp(s_a),
    \qquad
    \gamma_1=\frac75,\quad\gamma_2=\frac53,
\]
at the state
\[
    \alpha_1=\frac12,\quad
    \rho_1=2,\quad\rho_2=1,\quad
    s_1=\frac12,\quad s_2=0,\quad
    \VV=\mathbf0,\quad |\ww_1|=0.01.
\]
Here \(p_1\approx1.74\) and \(p_2\approx0.67\).
Direct evaluation gives an indefinite total-energy
Hessian in \(\qq\), but a positive definite one in
\(\qq^\varphi\).
At pressure equilibrium, by contrast, the additional
term vanishes.
The two total-energy Hessians are then congruent,
so the Hessian in \(\qq^\varphi\) is also indefinite
whenever \(\ww_1\ne\mathbf0\), by
Proposition~\ref{prop:nonconvex_q}.
\end{remark}

\subsection{Relative velocities}\label{app:relative_velocities}

We now compare the total energy in the phase variables
and in the evolved variables.
Recall from Section~\ref{subsec:configurational_direction} that the configurational direction changes the amount of the phases at fixed intrinsic states, \(\delta Z_a=(1,Y_a)\,\delta\alpha_a\). On admissible variations \(\sum_a\delta\alpha_a=0\), so that for two phases it is the single direction \(\delta Z=\bigl((1,Y_1),\,-(1,Y_2)\bigr)\,\delta\alpha_1\). Along it the energy in \(Z\) varies linearly and the forces do not change; the question is whether this flat direction survives the change of variables to \(\qq\).
The only nonlinear rows of \(Z\mapsto\qq\) are those
of the relative velocities.

\begin{proposition}[Energy in \(\qq\) and in \(Z\)]\label{prop:q_Z}
For two phases, let
\[
    \mathbf J=\frac{\partial\qq}{\partial Z},
    \qquad
    \mu=\frac{\crho_1\crho_2}{\rho}.
\]
On admissible variations,
\begin{equation}\label{eq:hessian_Z_q}
	\nabla^2_Z\rhoE
	=
	\mathbf J^{\mathsmaller T}\,\nabla^2_{\qq}\rhoE\,\mathbf J
	+
	\rhoE_{\ww_1}\cdot\nabla^2_Z\ww_1,
	\qquad
	\rhoE_{\ww_1}=\mu\,\ww_1,
	\qquad
	\mu=\frac{\crho_1\crho_2}{\rho}.
\end{equation}
Consequently
\begin{enumerate}
	\item[(i)] If \(\ww_1=\mathbf 0\), the Hessians in \(\qq\) and \(Z\) are congruent: \(\nabla^2_{\qq}\rhoE\) is positive semi-definite and the configurational direction, written in \(\qq\) as \(\mathbf J\,\delta Z\), belongs to its kernel.
	\item[(ii)] If \(\ww_1\neq\mathbf 0\), \(\nabla^2_{\qq}\rhoE\) is indefinite (Proposition~\ref{prop:nonconvex_q}), and the configurational direction is no longer in its kernel.
\end{enumerate}
\end{proposition}

\begin{proof}
All components of \(\qq(Z)\) in \eqref{eq:q_of_Z} are linear except
\[
    \ww_1=\frac{\bm M_1}{\crho_1}
          -\frac{\bm M_2}{\crho_2},
\]
whose force is \(\rhoE_{\ww_1}=\crho_1(\vvv_1-\VV)=\mu\ww_1\) by \eqref{eq:forces_phi_a}. Lemma~\ref{lem:change_of_variables} gives \eqref{eq:hessian_Z_q}.
If \(\ww_1=\mathbf0\), the additional term vanishes.
For every admissible variation \(\bm z\),
\[
    (\mathbf J\bm z)^{\mathsmaller T}
    \nabla^2_{\qq}\rhoE\,(\mathbf J\bm z)
    =
    \bm z^{\mathsmaller T}
    \nabla^2_Z\rhoE\,\bm z
    \ge0.
\]
The Jacobian maps admissible variations onto all
variations of \(\qq\), which proves positive
semi-definiteness.
For a configurational variation, the right-hand side
is zero.
A zero-curvature direction of a positive semi-definite
matrix belongs to its kernel, proving (i).

(ii)
For nonzero slip, Proposition~\ref{prop:nonconvex_q}
provides a configurational direction \(\bm n\) and
a mixture-momentum direction \(\bm e\) such that
\[
    \bm n^{\mathsmaller T}
    \nabla^2_{\qq}\rhoE\,\bm n=0,
    \qquad
    \bm e^{\mathsmaller T}
    \nabla^2_{\qq}\rhoE\,\bm n\ne0.
\]
The second identity implies
\(\nabla^2_{\qq}\rhoE\,\bm n\ne\mathbf0\).
Thus zero curvature does not imply membership in the
kernel when the Hessian is indefinite.
\end{proof}

For \(\Nph\) phases, the additional term becomes
\[
    \sum_{a<\Nph}
    \rhoE_{\ww_a}\cdot\nabla^2_Z\ww_a,
    \qquad
    \rhoE_{\ww_a}=\crho_a(\vvv_a-\VV).
\]
The argument in (i) therefore applies when all
relative velocities vanish.

\subsection{The pairing of two states}\label{app:pairing}

Convexity is a property of the Hessian at one state; a numerical dissipation compares two states. We now make explicit the remainder \(T_w\) in
\eqref{eq:pairing_split}.

\begin{proposition}[Pairing in \(\qq\) and in \(Z\)]\label{prop:Tw}
Let \(L\) and \(R\) be two states and let \(D_s\ww_a[\cdot]\) denote the differential of \(\ww_a(Z)\) at the state \(s\in\{L,R\}\). Then
\[
    \Delta\pp\cdot\Delta\qq
    =
    \sum_a\Delta\pp_a^Z\cdot\Delta Z_a+T_w,
\]
where
\begin{equation}\label{eq:Tw_general}
    T_w
    =
    \sum_{a<\Nph}
    \left[
        \rhoE_{\ww_a,R}\cdot
        \bigl(\Delta\ww_a-D_R\ww_a[\Delta Z]\bigr)
        -
        \rhoE_{\ww_a,L}\cdot
        \bigl(\Delta\ww_a-D_L\ww_a[\Delta Z]\bigr)
    \right].
\end{equation}
Equivalently, in terms of the drift velocities of the
phases relative to the mixture,
\(
    \bm d_{a,s}=\vvv_{a,s}-\VV_s,
\)
the remainder is
\begin{equation}\label{eq:Tw_explicit}
    T_w
    =
    \sum_{a=1}^{\Nph}
    \Delta\crho_a\,
    (\bm d_{a,L}+\bm d_{a,R})
    \cdot\Delta\vvv_a .
\end{equation}
In particular \(T_w=0\) if all relative velocities vanish in both states, \(T_w\) is quadratic in the jumps, like the pairing itself, and it takes both signs.
\end{proposition}

\begin{proof}
Since both states satisfy the saturation constraint,
\(\Delta Z\) is an admissible variation.
Applying \eqref{eq:chain_rule_admissible} at each
endpoint gives
\[
    \sum_a\Delta\pp_a^Z\cdot\Delta Z_a
    =
    \pp_R\cdot\mathbf J_R\Delta Z
    -
    \pp_L\cdot\mathbf J_L\Delta Z.
\]
Subtracting this identity from
\(\Delta\pp\cdot\Delta\qq\) gives
\[
    T_w
    =
    \pp_R\cdot(\Delta\qq-\mathbf J_R\Delta Z)
    -
    \pp_L\cdot(\Delta\qq-\mathbf J_L\Delta Z).
\]
The linear rows contribute zero.
Only the relative-velocity rows remain, yielding
\eqref{eq:Tw_general}.

For a phase velocity
\(\vvv_a=\bm M_a/\crho_a\), the endpoint
differentials satisfy
\[
    D_R\vvv_a[\Delta Z]
    =\frac{\crho_{a,L}}{\crho_{a,R}}\Delta\vvv_a,
    \qquad
    D_L\vvv_a[\Delta Z]
    =\frac{\crho_{a,R}}{\crho_{a,L}}\Delta\vvv_a.
\]
Hence
\[
\begin{aligned}
    \Delta\vvv_a-D_R\vvv_a[\Delta Z]
    =
    \frac{\Delta\crho_a}{\crho_{a,R}}
    \Delta\vvv_a,
    \qquad
    \Delta\vvv_a-D_L\vvv_a[\Delta Z]
    =
    -\frac{\Delta\crho_a}{\crho_{a,L}}
    \Delta\vvv_a .
\end{aligned}
\]
Insert these expressions into
\eqref{eq:Tw_general}, using
\[
    \rhoE_{\ww_a,s}=\crho_{a,s}\bm d_{a,s},
    \qquad
    \sum_{a=1}^{\Nph}\crho_{a,s}\bm d_{a,s}
    =\mathbf0.
\]
The result is \eqref{eq:Tw_explicit}, including
the contribution of the reference phase.
The formula contains products of mass and velocity
jumps, weighted by the endpoint drifts.
It therefore vanishes when all drifts vanish and
is generally quadratic in small jumps, that is, of the
same order as the pairing itself, so that it can change
its sign.
To see that either sign is possible, consider two
phases with
\[
    \vvv_{1,L}=\vvv_{2,L}=\vvv_{2,R}=\mathbf0,
    \qquad
    \vvv_{1,R}=\ww_{1,R}\ne\mathbf0,
    \qquad
    \Delta\crho_2=0.
\]
Then
\[
    T_w
    =
    \frac{\mu_R}{\crho_{1,R}}
    |\ww_{1,R}|^2\,\Delta\crho_1,
\]
whose sign is that of \(\Delta\crho_1\).
\end{proof}

Equation~\eqref{eq:Tw_explicit} identifies the
coupling responsible for the difference between
the two pairings.
It involves changes of phase mass and phase velocity
in the presence of slip.
These conditions allow a nonzero remainder, but do
not force one: individual terms may vanish or cancel.
In particular, zero slip at both endpoints is
sufficient for \(T_w=0\), independently of pressure
equilibrium.

For two phases, negative pairings exist arbitrarily
close to any admissible state with
\(\ww_1\ne\mathbf0\).
Indeed, Proposition~\ref{prop:nonconvex_q} provides
a direction \(\bm\xi\) with
\[
    \bm\xi^{\mathsmaller T}
    \nabla^2_{\qq}\rhoE(\qq_L)\,\bm\xi<0.
\]
Taking \(\qq_R=\qq_L+h\bm\xi\), which remains
admissible for sufficiently small \(h\), gives
\[
    \Delta\pp\cdot\Delta\qq
    =
    h^2\bm\xi^{\mathsmaller T}
    \nabla^2_{\qq}\rhoE(\qq_L)\,\bm\xi
    +o(h^2)<0.
\]
Thus reducing the jump size alone does not guarantee
a non-negative pairing in the evolved variables.

\section{Characteristic structure and asymptotic limits}\label{app:characteristics}
This appendix records two complementary properties of the fluid subsystem for two phases: the frozen characteristic structure at equal phase velocities, and the mechanically relaxed acoustic response obtained when interfacial friction and pressure relaxation are stiff, together with the dispersion spectrum between the two. The results are stated for the original formulation \eqref{eq:system_phi}; by Proposition~\ref{prop:continuous_equivalence}, they apply equally to the volume fraction formulation used by the numerical method.

\subsection{Frozen characteristic structure for two phases}\label{characteristic.str1}

We consider two phases without relaxation. Writing the original system as
\(
\partial_t\qq^\varphi+\partial_k\f_k(\qq^\varphi)+\B_k(\qq^\varphi)\,\partial_k\qq^\varphi=\mathbf 0,
\)
where \(\B_k\) collects the non-conservative terms of the relative velocity equation, the state vector reduces to
\begin{equation}\label{eq.state.var4}
	\qq^\varphi=
	\left(\UU,\,\crho_1,\,\ww_1,\,\eta_1,\,\eta_2,\,\rho,\,\varphi_1\right)^{\mathsmaller T}.
\end{equation}
The characteristic speeds are most easily read in the primitive phase variables
\begin{equation}\label{eq.state.var3}
	\rr=
	\left(\VV,\,\ww_1,\,\rho_1,\,\rho_2,\,s_1,\,s_2,\,\alpha_1\right)^{\mathsmaller T},
\end{equation}
in which the system takes the quasi-linear form
\begin{equation}\label{eq.QL}
	\partial_t\rr+\mathbf C_k(\rr)\,\partial_k\rr=\mathbf 0,
	\qquad
	\mathbf C_k=\frac{\partial\rr}{\partial\qq^\varphi}\left(\frac{\partial\f_k}{\partial\qq^\varphi}+\B_k\right)\frac{\partial\qq^\varphi}{\partial\rr}.
\end{equation}
The characteristic structure of this system is consistent with the SHTC two-phase models of \cite{Rom2001,Romenski2004}; related analyses are available for the barotropic case \cite{Thein2022} and for the single-entropy limit \cite{RomDrikToro2010}. Here we record the explicit spectrum of the fully thermal, non-equilibrium subsystem.

Let the two phases move with the same velocity, \(\ww_1=\mathbf 0\), and consider the direction \(x=x_1\). The transverse components \(V_2,V_3,w_{1,2},w_{1,3}\) are simply transported, \(\partial_t+V_1\partial_x\), and decouple. In the longitudinal variables \((V_1,w_{1,1},\rho_1,\rho_2,s_1,s_2,\alpha_1)\) the matrix \(\mathbf C_1\) reads
\begin{equation}\label{eq.matrix.c1}
	\mathbf C_1^{\parallel}=
	\begin{pmatrix}
		V_1 & 0 & \dfrac{\alpha_1C_1^2}{\rho} & \dfrac{\alpha_2C_2^2}{\rho} & \dfrac{\alpha_1\rho_1^2e^{(1,1)}_1}{\rho} & \dfrac{\alpha_2\rho_2^2e^{(1,1)}_2}{\rho} & \dfrac{p_1-p_2}{\rho}\\[2mm]
		0 & V_1 & \xi_1 & -\xi_2 & -\zeta_1 & \zeta_2 & 0\\[1mm]
		\rho_1 & \rho_1c_2 & V_1 & 0 & 0 & 0 & 0\\
		\rho_2 & -\rho_2c_1 & 0 & V_1 & 0 & 0 & 0\\
		0 & -s_1c_2 & 0 & 0 & V_1 & 0 & 0\\
		0 & s_2c_1 & 0 & 0 & 0 & V_1 & 0\\
		0 & 0 & 0 & 0 & 0 & 0 & V_1
	\end{pmatrix},
\end{equation}
where \(c_a=\crho_a/\rho\) are the mass fractions, \(e^{(i,j)}_a=\partial^i_{\rho_a}\partial^j_{s_a}e_a\) denote the derivatives of the specific internal energy, \(C_a\) is the adiabatic sound speed of phase \(a\),
\begin{equation}
	C_a^2=\left.\frac{\partial p_a}{\partial\rho_a}\right|_{s_a}=2\rho_ae^{(1,0)}_a+\rho_a^2e^{(2,0)}_a,
\end{equation}
and
\begin{equation}
	\xi_a=\frac{\partial g_a}{\partial\rho_a}=\frac{C_a^2}{\rho_a}-s_ae^{(1,1)}_a,
	\qquad
	\zeta_a=-\frac{\partial g_a}{\partial s_a}=s_ae^{(0,2)}_a-\rho_ae^{(1,1)}_a .
\end{equation}
The first row is the momentum equation, driven by the gradient of \(P=\alpha_1p_1+\alpha_2p_2\); the second is the relative velocity equation, driven by the gradient of \(g_1-g_2\) (Section~\ref{subsec:balances}); the next four express that each phase carries its mass with its own velocity while the phase entropies are carried with \(\VV\).

Together with the four transverse components, the system has seven eigenvalues equal to \(V_1\), associated with linearly degenerate waves moving with the fluid, and four thermo-acoustic eigenvalues
\begin{equation}\label{eq.eigenmix}
	\lambda_{8-11}=V_1\pm\frac{1}{\rho\sqrt2}\sqrt{K_1\pm\sqrt{K_1^2-K_2}},
\end{equation}
with
\begin{align}
	K_1&=\rho\Bigl(\rho\,(C_1^2+C_2^2)+\alpha_1\rho_1s_2^2e^{(0,2)}_2+\alpha_2\rho_2s_1^2e^{(0,2)}_1-2\alpha_2\rho_2\rho_1s_1e^{(1,1)}_1-2\alpha_1\rho_1\rho_2s_2e^{(1,1)}_2\Bigr),
	\nonumber\\
	K_2&=4\rho^2\Bigl[\rho^2C_1^2C_2^2
	+\alpha_1\rho_1s_2^2\Bigl(\alpha_1\rho_1C_1^2e^{(0,2)}_2+\alpha_2\rho_2\bigl(C_2^2e^{(0,2)}_2-\rho_2^2(e^{(1,1)}_2)^2\bigr)\Bigr)
	\nonumber\\
	&\hspace{11mm}
	+\alpha_2\rho_2s_1^2\Bigl(\alpha_2\rho_2C_2^2e^{(0,2)}_1+\alpha_1\rho_1\bigl(C_1^2e^{(0,2)}_1-\rho_1^2(e^{(1,1)}_1)^2\bigr)\Bigr)
	\nonumber\\
	&\hspace{11mm}
	-2\rho_1\rho_2\Bigl(\alpha_1\rho\,s_2C_1^2e^{(1,1)}_2+\alpha_2\rho\,s_1C_2^2e^{(1,1)}_1-\alpha_1\rho_1\alpha_2\rho_2\,s_1s_2\,e^{(1,1)}_1e^{(1,1)}_2\Bigr)\Bigr].
	\label{eq.K1K2}
\end{align}
These are four mixed waves whose speeds involve both phases. For \(s_1=s_2=0\), \eqref{eq.K1K2} gives \(K_1=\rho^2(C_1^2+C_2^2)\), \(K_2=4\rho^4C_1^2C_2^2\), and the four speeds reduce to \(V_1\pm C_1\) and \(V_1\pm C_2\) for every composition. In the pure-phase limits three statements must be distinguished.
\begin{itemize}
	\item \emph{The subsystem of the dominant phase.} At \(\alpha_1=1\), the equations of \((V_1,\rho_1,s_1)\) decouple from \((w_{1,1},\rho_2,s_2)\), whose coefficients in the first, third and fifth rows of \eqref{eq.matrix.c1} vanish (the coupling with \(\alpha_1\) does not affect the spectrum, since \(\alpha_1\) is only transported), and they form the Euler system with sound speed \(C_1\).
	\item \emph{The complete frozen spectrum.} Before the variables of the trace phase are eliminated, the spectrum at \(\alpha_1=1\) still contains a pair for the vanishing phase: \eqref{eq.eigenmix} gives \(V_1\pm C_1\) and \(V_1\pm\widetilde C_2\), with \(\widetilde C_2^{\,2}=C_2^2-2\rho_2s_2e^{(1,1)}_2+s_2^2e^{(0,2)}_2\); symmetrically for \(\alpha_1\to0\).
	\item \emph{The relaxed limit.} With stiff friction and pressure relaxation the response is that of a single mechanically equilibrated mixture, with Wood's sound speed (Section~\ref{characteristic.str2}).
\end{itemize}
The speeds also depend on the phase entropies through the products \(s_ae^{(0,2)}_a\) and \(s_ae^{(1,1)}_a\), a consequence of the fact that the phase entropies are carried with the mixture velocity while the phase masses are carried with the phase velocities.

\subsection{Relaxed acoustic limit and dispersion}
\label{characteristic.str2}

The frozen eigenvalues in \eqref{eq.eigenmix} describe the short-wave response of the reduced hyperbolic system. We now retain the two irreversible mechanisms introduced in Section~\ref{subsec:balances}, namely interfacial friction and pressure relaxation \eqref{eq:relaxation}, show that their joint stiff limit enforces mechanical equilibrium and recovers Wood's acoustic response, and describe the dispersion spectrum that connects the two limits.

For the binary subsystem, the two source terms reduce to
\begin{equation}
\label{eq:binary_relaxation_sources}
\Lambda_{1,k}=-\lambda_{11,k}\rhoE_{w_{1,k}},
\qquad
\Phi_1=-\phi_{11}\rhoE_{\varphi_1},
\qquad k=1,2,3,
\end{equation}
where \(\lambda_{11,k}>0\) and \(\phi_{11}>0\) are proportional to the inverse relative velocity and pressure-relaxation times, respectively. The associated entropy-production terms are quadratic in the thermodynamic forces and therefore do not contribute to the linearized source operator about equilibrium.

Let \(\rr_0\) be a spatially uniform reference state with positive phase densities. In the stiff-relaxation limit, or equivalently at leading order for long-wave acoustic disturbances, the source terms in \eqref{eq:binary_relaxation_sources} impose
\begin{equation}
\label{eq:mechanical_equilibrium_forces}
\rhoE_{w_{1,k}}=0,
\qquad
\rhoE_{\varphi_1}=0.
\end{equation}
For two phases, by \eqref{eq:forces_phi_a}, these two conditions can be written explicitly as
\begin{equation}
\label{eq:mechanical_equilibrium_conditions}
\rhoE_{w_{1,k}}
=
\frac{\crho_1\crho_2}{\rho}\,w_{1,k},
\qquad
\rhoE_{\varphi_1}
=
\frac{p_2-p_1}{\rho}.
\end{equation}
Hence, \eqref{eq:mechanical_equilibrium_forces} yields
\begin{equation}
\label{eq:mechanical_equilibrium_state}
\ww_1=\mathbf 0,
\qquad
p_1=p_2.
\end{equation}
The two sources therefore constrain the mixture to the velocity- and pressure-equilibrium manifold. No temperature equilibration is assumed here: the phase entropies may remain distinct, while acoustic perturbations are phase-wise adiabatic.

To make this connection explicit, we linearize the source-augmented two-phase system about \(\rr_0\) and restrict attention to one-dimensional perturbations. With \(\widetilde{\rr}=\rr-\rr_0\), the resulting system can be written as
\begin{equation}
\label{eq:linear_relaxation_system}
\partial_t\widetilde{\rr}
+
\mathbb C\,\partial_x\widetilde{\rr}
=
\mathbb S\,\widetilde{\rr},
\end{equation}
where \(\mathbb C=\mathbf C_1(\rr_0)\) is the frozen matrix of \eqref{eq.QL} and \(\mathbb S\) is the Jacobian of the relaxation operator in the variables \(\rr\). For normal modes of the form
\begin{equation}
\label{eq:appendix_plane_wave}
\widetilde{\rr}(t,x)
=
\widehat{\rr}\,e^{i(kx-\omega t)},
\end{equation}
the complex frequencies are obtained from the dispersion relation
\begin{equation}
\label{eq:appendix_dispersion_relation}
\det\left(
 k\mathbb C+i\mathbb S-\omega\mathbb I
\right)=0.
\end{equation}
For a prescribed real wavenumber \(k\), the roots are computed numerically as the eigenvalues of \(k\mathbb C+i\mathbb S\). We write them as
\begin{equation}
\label{eq:omega_split}
	\omega_j(k)=\nu_j(k)-i\,\gamma_j(k),
	\qquad
	V_{\mathrm{ph},j}(k)=\frac{\nu_j(k)}{k},
\end{equation}
with \(\nu_j\) the oscillation frequency, \(\gamma_j\) the temporal damping rate and \(V_{\mathrm{ph},j}\) the phase velocity; with the convention \(e^{-i\omega t}\), \(\gamma_j>0\) means decay.

In the short-wave limit, \(k\to\infty\), the hyperbolic transport operator dominates the finite relaxation operator. Writing the eigenvalue problem as \((\mathbb C+i\epsilon\,\mathbb S)\widehat{\rr}=(\omega/k)\,\widehat{\rr}\) with \(\epsilon=1/k\), first-order perturbation theory gives, for a simple eigenvalue \(\lambda_j\) of \(\mathbb C\) with right and left eigenvectors \(\bm r_j\), \(\bm l_j\),
\begin{equation}
\label{eq:short_wave_limit}
	\omega_j(k)=k\,\lambda_j+i\,\frac{\bm l_j^{\mathsmaller T}\mathbb S\,\bm r_j}{\bm l_j^{\mathsmaller T}\bm r_j}+\mathcal O(k^{-1}),
	\qquad
	\lim_{k\to\infty}\gamma_j(k)=-\operatorname{Re}\frac{\bm l_j^{\mathsmaller T}\mathbb S\,\bm r_j}{\bm l_j^{\mathsmaller T}\bm r_j}.
\end{equation}
The propagative branches approach the four frozen thermo-acoustic speeds \eqref{eq.eigenmix}, which is the characteristic structure relevant to the high-frequency response and to a CFL restriction based on the unreduced hyperbolic subsystem, and their damping rates saturate at finite values: relaxation damps the frozen waves at a bounded rate without changing their speeds.

Conversely, we consider the propagative acoustic branches for which the frequency admits the long-wave expansion
\begin{equation}
\omega_\pm(k)
=
\pm c_{\mathrm{eq}}\,k
+
\mathcal O(k^2),
\qquad
k\to0.
\end{equation}
Equivalently, these are the roots for which \(\omega=\mathcal O(k)\) and the phase velocity \(\omega/k\) remains finite in the long-wave limit. Substitution into \eqref{eq:appendix_dispersion_relation} then yields, at leading order,
\begin{equation}
\mathbb S\,\widehat{\rr}^{(0)}=\mathbf 0.
\end{equation}
Thus, the leading perturbation is tangent to the mechanical-equilibrium manifold \eqref{eq:mechanical_equilibrium_state}: it satisfies \(\delta\ww_1=\mathbf 0\) and \(\delta p_1=\delta p_2\). The corresponding equilibrium acoustic subsystem has the single propagative sound speed
\begin{equation}
\label{eq:wood_sound_speed}
\frac{1}{\rho c_{\mathrm W}^2}
=
\frac{\alpha_1}{\rho_1 C_1^2}
+
\frac{\alpha_2}{\rho_2 C_2^2},
\end{equation}
which is Wood's formula for a phase-wise adiabatic, mechanically equilibrated binary mixture \cite{wood}. Thus, the long-wave acoustic sector is consistent with the velocity- and pressure-equilibrium limit associated with the five-equation Kapila-type model \cite{Bdzil}. This statement concerns the mechanically relaxed acoustic response only; it does not require the phase entropies to coincide.

\paragraph{Numerical evaluation}
The spectrum is evaluated for the air--water mixture of Table~\ref{tab:wood_reference_state}. The specific entropies are measured from the equilibrium state, \(s_1=s_2=0\), so that by \ref{characteristic.str1} the frozen speeds are exactly \(V_1\pm C_1\) and \(V_1\pm C_2\). The relaxation coefficients are \(\lambda_{11}=1/(\rho\,\tau^{\mathrm w})\), so that friction alone gives \(\partial_t\ww_1=-c_1c_2\,\ww_1/\tau^{\mathrm w}\), and \(\phi_{11}=\rho/(p_{\mathrm{mix}}\,\tau^{\mathrm p})\), so that pressure relaxation alone gives \(\partial_t\alpha_1=(p_1-p_2)/(\rho\,p_{\mathrm{mix}}\,\tau^{\mathrm p})\); here \(\tau^{\mathrm w}=10^{-9}\) and \(\tau^{\mathrm p}=10^{-6}\), in SI units, are used as relaxation coefficients rather than as times.

\begin{table}[t]
\centering
\small
\caption{Reference equilibrium state used for the dispersion analysis. Phase~1 is modelled as an ideal gas and phase~2 as a stiffened gas. The common equilibrium pressure and temperature are \(p_{\mathrm{mix}}=10^5\,\mathrm{Pa}\) and \(T_{\mathrm{mix}}=283\,\mathrm{K}\), respectively; the specific entropies are measured from the equilibrium state.}
\label{tab:wood_reference_state}
\begin{tabular}{lcc}
\toprule
& Phase 1 (air) & Phase 2 (water) \\
\midrule
Equation of state & ideal gas & stiffened gas \\
Reference density \(\rho_{0a}\) & \(1.23\,\mathrm{kg\,m^{-3}}\) & \(1000\,\mathrm{kg\,m^{-3}}\) \\
Adiabatic exponent \(\gamma_a\) & \(1.4\) & \(1.88\) \\
Specific heat \(c_{v,a}\) & \(718\,\mathrm{J\,kg^{-1}\,K^{-1}}\) & \(4184\,\mathrm{J\,kg^{-1}\,K^{-1}}\) \\
Equilibrium density \(\rho_a\) & \(1.23\,\mathrm{kg\,m^{-3}}\) & \(1000.05\,\mathrm{kg\,m^{-3}}\) \\
Equilibrium specific entropy \(s_a\) & \(0\) & \(0\) \\
Adiabatic sound speed \(C_a\) & \(337.33\,\mathrm{m\,s^{-1}}\) & \(1399.97\,\mathrm{m\,s^{-1}}\) \\
\bottomrule
\end{tabular}
\end{table}

\begin{figure}[!htb]
\centering
\figB{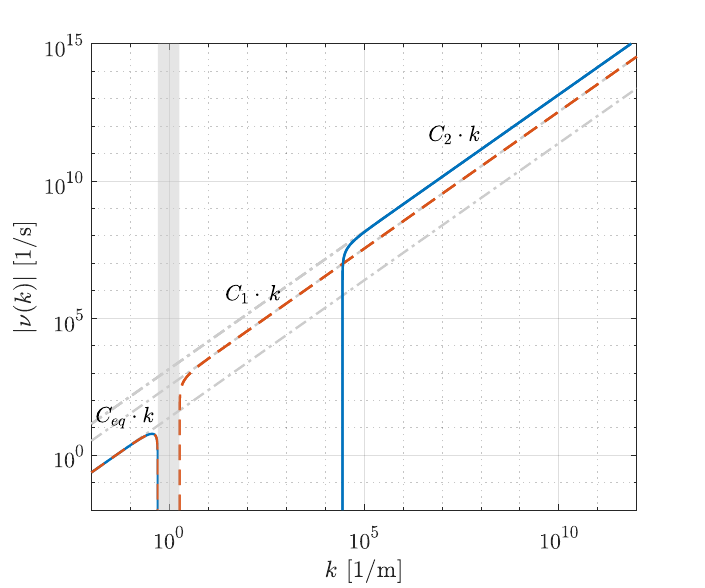}\hspace{3mm}%
\figB{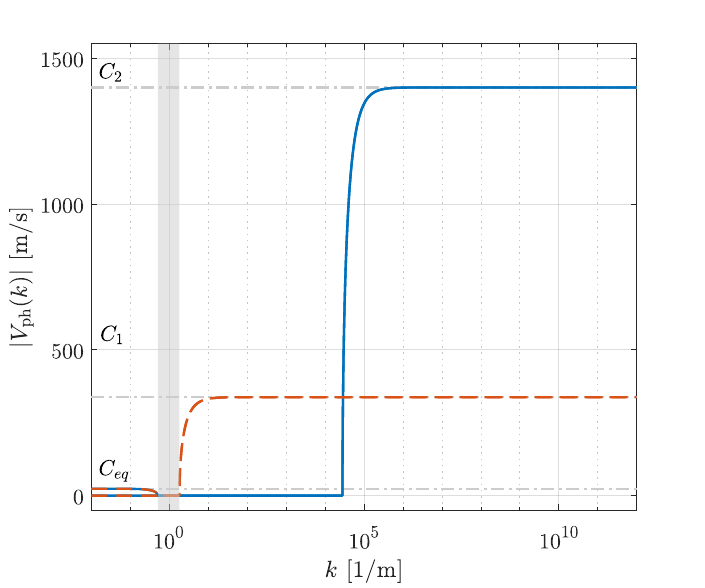}\hspace{3mm}%
\figB{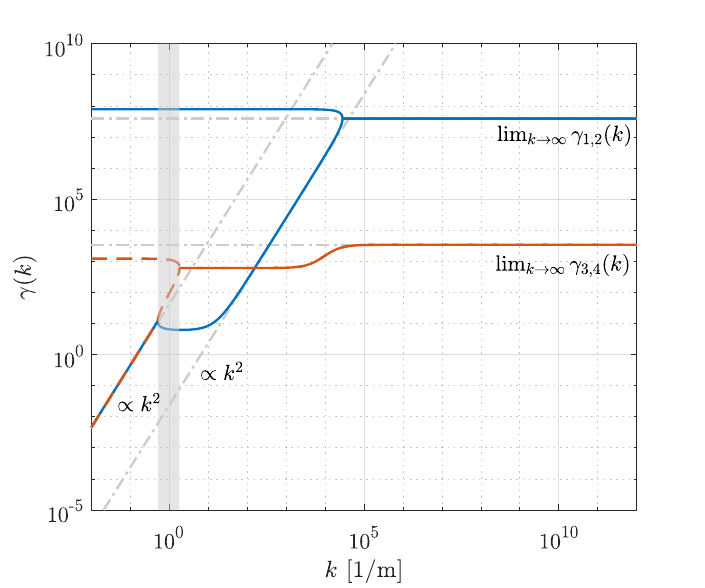}
\caption{Dispersion spectrum for \(\alpha_1=0.5\): oscillation frequencies \(|\nu_j(k)|\) (left), phase velocities \(|V_{\mathrm{ph},j}(k)|\) (centre) and damping rates \(\gamma_j(k)\) (right). Solid: branch of the liquid, approaching \(C_2\); dashed: branch of the gas, approaching \(C_1\). In the shaded band both acoustic branches are overdamped. The horizontal lines on the right are the limits \eqref{eq:short_wave_limit}.}
\label{fig:dispersion_symmetric}
\end{figure}

\paragraph{A symmetric mixture}
Figure~\ref{fig:dispersion_symmetric} shows the spectrum for \(\alpha_1=0.5\). Three regimes appear. For long waves the two acoustic branches merge into a single mode, which propagates with Wood's speed, here about \(24\,\mathrm{m\,s^{-1}}\), and is damped as \(\gamma\propto k^2\): on scales much larger than those of the relaxation, friction and pressure relaxation act on the acoustic wave as an effective bulk viscosity. When the acoustic frequency becomes comparable to the relaxation rates, the acoustic modes stop propagating: in the shaded band both branches are overdamped and \(V_{\mathrm{ph}}=0\). For shorter waves the branches re-emerge one at a time, the gas branch first and the liquid branch at \(k\approx10^4\,\mathrm{m^{-1}}\), approach the frozen speeds \(C_1\) and \(C_2\), and their damping rates saturate at the values \eqref{eq:short_wave_limit}.

\begin{figure}[!htb]
\centering
\figB{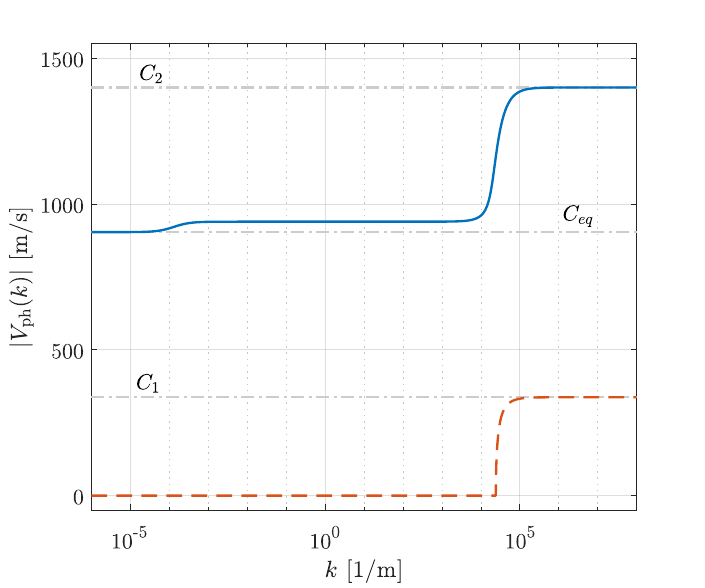}\hspace{3mm}%
\figB{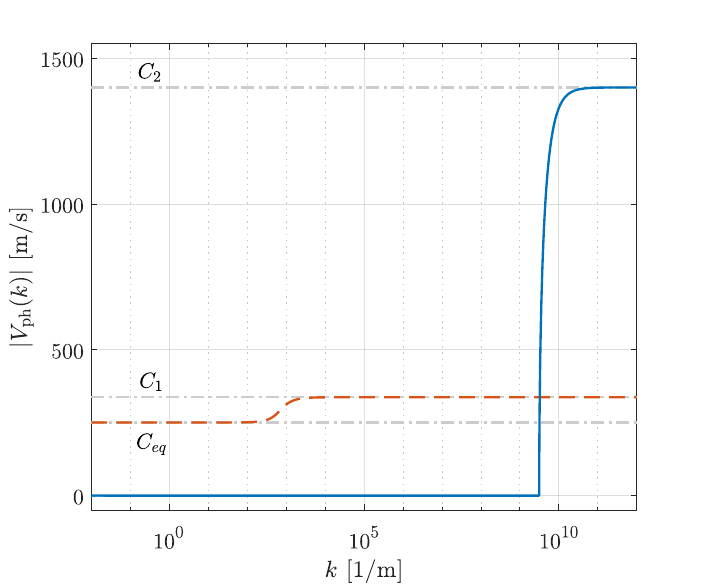}\\[2mm]%
\figB{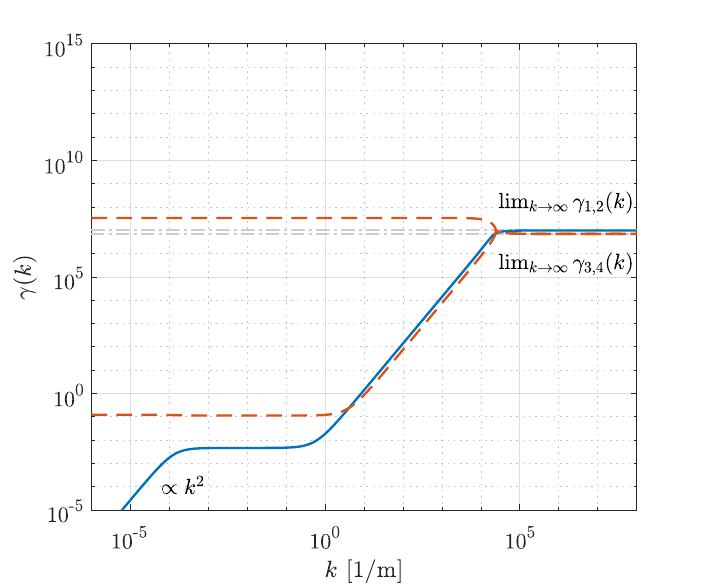}\hspace{3mm}%
\figB{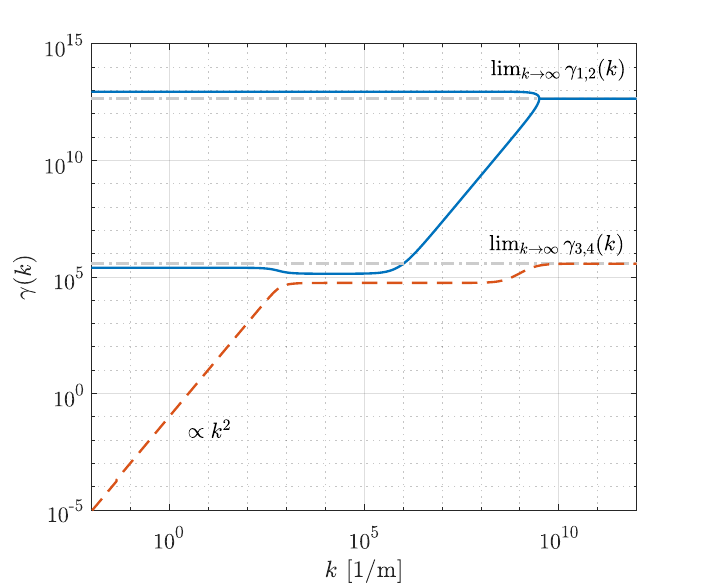}
\caption{Dilute mixtures: phase velocities (top) and damping rates (bottom) for a bubbly liquid, \(\alpha_1=10^{-4}\) (left), and a mist, \(\alpha_1=1-10^{-3}\) (right). Line styles as in Figure~\ref{fig:dispersion_symmetric}.}
\label{fig:dispersion_dilute}
\end{figure}

\paragraph{Dilute mixtures}
Near the pure-phase limits the common overdamped band disappears (Figure~\ref{fig:dispersion_dilute}). The branch of the continuous phase propagates at every wavenumber, from Wood's speed to its frozen speed, while the branch of the dispersed phase remains overdamped up to large wavenumbers. The extent of the equilibrium regime is strongly asymmetric: the regime \(\gamma\propto k^2\) ends near \(k\approx10^{-4}\,\mathrm{m^{-1}}\) in the bubbly liquid and extends to \(k\approx10^{2}\)--\(10^{3}\,\mathrm{m^{-1}}\) in the mist, about six orders of magnitude apart. The two dispersed phases play opposite roles. The bubbles add compressibility: to keep the mixture at its slow equilibrium speed they must expand and contract with the pressure wave carried by the fast liquid, which they can follow only for very long waves. The droplets add inertia: in the slow gas the drag keeps them moving with the carrier phase down to short wavelengths.

\begin{figure}[!htb]
\centering
\includegraphics[scale=0.369,trim=1.5bp 1.5bp 1.5bp 1bp,clip]{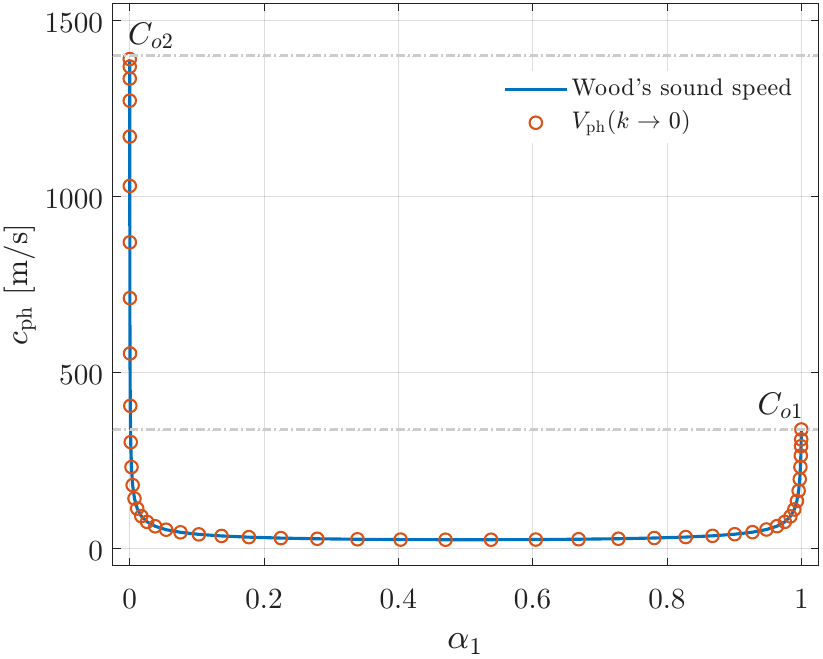}\hspace{3mm}%
\figB{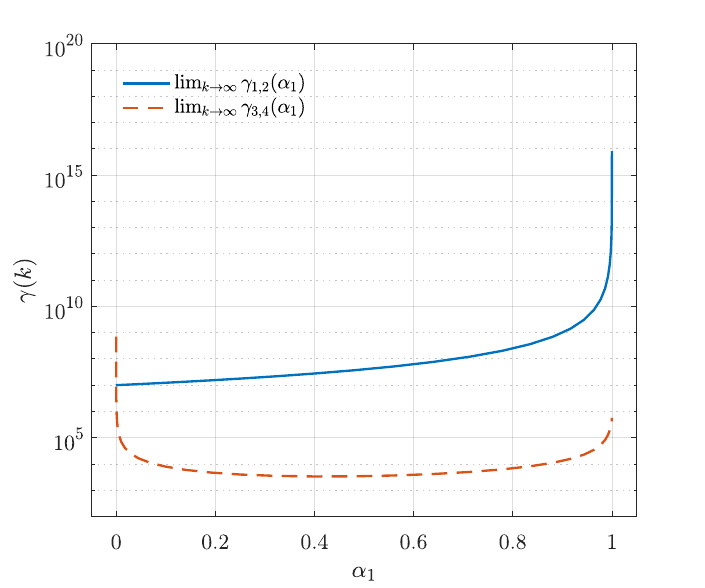}
\caption{Dependence on the composition. Left: phase velocity of the propagative root in the limit \(k\to0\) (circles) against Wood's speed \eqref{eq:wood_sound_speed} (solid). Right: high-frequency damping rates \eqref{eq:short_wave_limit} of the liquid pair (solid) and of the gas pair (dashed) as functions of \(\alpha_1\).}
\label{fig:wood_limit}
\end{figure}

\paragraph{Dependence on the composition}
Figure~\ref{fig:wood_limit} (left) verifies the long-wave limit over the whole composition range: for each \(\alpha_1\), the phase velocity of the propagative root in the \(k\to0\) regime coincides with Wood's speed \eqref{eq:wood_sound_speed}, which tends to \(C_2\) and \(C_1\) in the pure-phase limits. The relaxation times determine the wavenumber at which this regime is reached, not the limiting value. Figure~\ref{fig:wood_limit} (right) shows the high-frequency damping rates: the rate of the pair that belongs to the vanishing phase grows by several orders of magnitude as \(\alpha_1\to0\) or \(\alpha_1\to1\). This complements the frozen analysis of \ref{characteristic.str1}: the pair of a trace phase remains in the frozen spectrum, but relaxation damps it ever more strongly as the phase disappears.

The reduced HTC system therefore connects the frozen hyperbolic structure at short wavelengths with the mechanically relaxed Kapila--Wood response at long wavelengths; between the two, the composition of the mixture sets the scale at which the phases decouple.

\section{Constitutive relations and robust evaluation of the Bregman defects}\label{app:bregman}
We consider the two equations of state used in Section~\ref{sec:results}, the ideal gas and the stiffened gas. For each phase, dropping the phase index,
\begin{equation}\label{eq:eos}
	e(\rho,s)=e_{\mathrm{id}}(\rho,s)+\frac{B}{\rho},
	\qquad
	e_{\mathrm{id}}(\rho,s)=\frac{c_0^2}{\gamma(\gamma-1)}\Bigl(\frac{\rho}{\rho_0}\Bigr)^{\gamma-1}\exp\Bigl(\frac{s}{c_v}\Bigr),
	\qquad
	B=\frac{\rho_0c_0^2-\gamma p_0}{\gamma},
\end{equation}
with \(B=0\) for the ideal gas; \(B\) plays the role of the stiffness pressure \(p_\infty\). The pressure, the temperature and the specific Gibbs energy are
\[
	p=(\gamma-1)\,\rho\,e_{\mathrm{id}}-B,
	\qquad
	T=\frac{e_{\mathrm{id}}}{c_v},
	\qquad
	g=e+\frac p\rho-Ts .
\]
In terms of the densities, \(\varepsilon(\rho,\sigma)=\rho\,e_{\mathrm{id}}(\rho,\sigma/\rho)+B\): the stiffness enters as a constant and does not contribute to the divergence \eqref{eq:bregman_fluid}.

\paragraph{Closed form}
Let \(Y=(\rho,\bm m,\sigma)\) and \(Y'=(\rho',\bm m',\sigma')\) be two intrinsic states, \(e_{\mathrm{id}}\) be evaluated at \(Y\), and
\[
	t=\frac{\rho'}{\rho},
	\qquad
	z=(\gamma-1)\ln t+\frac{s'-s}{c_v},
	\qquad\text{so that}\qquad
	e_{\mathrm{id}}(\rho',s')=e_{\mathrm{id}}\,e^{z}.
\]
With \(\partial_\rho\varepsilon=e_{\mathrm{id}}(\gamma-s/c_v)\) and \(\partial_\sigma\varepsilon=e_{\mathrm{id}}/c_v\), the thermal part of \eqref{eq:bregman_fluid} is
\[
	\rho'e_{\mathrm{id}}e^{z}-\rho e_{\mathrm{id}}-e_{\mathrm{id}}\Bigl(\gamma-\frac{s}{c_v}\Bigr)(\rho'-\rho)-\frac{e_{\mathrm{id}}}{c_v}(\rho's'-\rho s)
	=\rho'e_{\mathrm{id}}\Bigl[e^{z}-\gamma+(\gamma-1)\frac{\rho}{\rho'}-\frac{s'-s}{c_v}\Bigr],
\]
and, using \((s'-s)/c_v=z-(\gamma-1)\ln t\),
\begin{equation}\label{eq:bregman_closed}
	D_E(Y'\,|\,Y)
	=\tfrac12\rho'|\vvv'-\vvv|^2
	+\rho'\,e_{\mathrm{id}}\,\phi(z)
	+(\gamma-1)\,\rho\,e_{\mathrm{id}}\,\psi(t),
	\qquad
	\phi(z)=e^{z}-1-z,
	\qquad
	\psi(t)=t\ln t-t+1 .
\end{equation}
Both functions are non-negative, \(\phi(z)\ge0\) for every \(z\) and \(\psi(t)\ge0\) for \(t>0\), and both vanish quadratically, \(\phi(z)\simeq z^2/2\) and \(\psi(t)\simeq(t-1)^2/2\). The three terms of \eqref{eq:bregman_closed} are therefore separately non-negative: the divergence is a sum of non-negative quantities, and its sign does not depend on any cancellation.

\paragraph{Evaluation}
For small jumps, \(\phi\) and \(\psi\) are differences of nearly equal numbers and are evaluated by their series: for \(|z|<1/8\), \(\phi(z)=\sum_{k\ge2}z^k/k!\); for \(|t-1|<1/8\), with \(x=t-1\), \(\psi(t)=\sum_{k\ge2}(-1)^k x^k/(k(k-1))\), and \(\ln t=2\,\mathrm{artanh}\bigl(x/(2+x)\bigr)\) is summed in the same way; otherwise the closed forms are used. The divergences, and with them the defects \eqref{eq:defects}, the productions \eqref{eq:production} and the denominator \(\delta_{ij}\) of \eqref{eq:abgrall_theta}, are thus computed to full relative accuracy, also when the jump across a face is small. Evaluating \(\delta_{ij}\) as \(\sum_a\Delta\pp^Z_a\cdot\Delta Z_a\) would instead subtract nearly equal products. A quadrature of the integral forms \eqref{eq:bregman_integral}, as done for the Roe matrix of the Hessian in \cite{HTC2022,ABGRALL2023127629}, is also avoided: near the traces the integrand of \eqref{eq:perspective_hessian} grows like \(1/\alpha_a(s)\), and the energy identities would hold only up to the quadrature error.

\bibliographystyle{plain}
\bibliography{references}

\end{document}